\documentclass[11pt]{article}

\usepackage{etex}

\usepackage{hyperref}

\usepackage{tikz-cd}
  
  \usepackage{rotating}
  
\newenvironment{tikzar}[1][]{{}\kern-4pt\begin{tikzcd}[ampersand replacement=\&,#1]}%
{\end{tikzcd}\kern-4pt{}}

\usepackage{amssymb}
\usepackage{amsmath}
\usepackage{amsthm}
\usepackage{newtxtext,newtxmath}
\RequirePackage{textcomp}
\usepackage{stmaryrd}
\usepackage{anyfontsize}
\SetSymbolFont{stmry}{bold}{U}{stmry}{m}{n}

\usepackage{bbm}
\usepackage{enumerate}
\usepackage{fullpage}
\usepackage{multicol}

\RequirePackage[T1]{fontenc}
\usepackage[utf8]{inputenc}
\RequirePackage{authblk}

\theoremstyle{plain}
\newtheorem{theorem}{Theorem}[section]
\newtheorem{proposition}[theorem]{Proposition}
\newtheorem{lemma}[theorem]{Lemma}
\newtheorem{corollary}[theorem]{Corollary}

\theoremstyle{definition}
\newtheorem{definition}[theorem]{Definition}
\theoremstyle{remark}

\newdimen\proofrulebreadth \proofrulebreadth=.05em
\newdimen\proofdotseparation \proofdotseparation=1.25ex
\newdimen\proofrulebaseline \proofrulebaseline=2ex
\newcount\proofdotnumber \proofdotnumber=3
\let\then\relax
\def\hfi{\hskip0pt plus.0001fil}
\mathchardef\squigto="3A3B
\newif\ifinsideprooftree\insideprooftreefalse
\newif\ifonleftofproofrule\onleftofproofrulefalse
\newif\ifproofdots\proofdotsfalse
\newif\ifdoubleproof\doubleprooffalse
\let\wereinproofbit\relax
\newdimen\shortenproofleft
\newdimen\shortenproofright
\newdimen\proofbelowshift
\newbox\proofabove
\newbox\proofbelow
\newbox\proofrulename
\def\shiftproofbelow{\let\next\relax\afterassignment\setshiftproofbelow\dimen0 }
\def\shiftproofbelowneg{\def\next{\multiply\dimen0 by-1 }%
\afterassignment\setshiftproofbelow\dimen0 }
\def\setshiftproofbelow{\next\proofbelowshift=\dimen0 }
\def\setproofrulebreadth{\proofrulebreadth}

\def\prooftree{
\ifnum  \lastpenalty=1
\then   \unpenalty
\else   \onleftofproofrulefalse
\fi
\ifonleftofproofrule
\else   \ifinsideprooftree
        \then   \hskip.5em plus1fil
        \fi
\fi
\bgroup
\setbox\proofbelow=\hbox{}\setbox\proofrulename=\hbox{}%
\let\justifies\proofover\let\leadsto\proofoverdots\let\Justifies\proofoverdbl
\let\using\proofusing\let\[\prooftree
\ifinsideprooftree\let\]\endprooftree\fi
\proofdotsfalse\doubleprooffalse
\let\thickness\setproofrulebreadth
\let\shiftright\shiftproofbelow \let\shift\shiftproofbelow
\let\shiftleft\shiftproofbelowneg
\let\ifwasinsideprooftree\ifinsideprooftree
\insideprooftreetrue
\setbox\proofabove=\hbox\bgroup$\displaystyle 
\let\wereinproofbit\prooftree
\shortenproofleft=0pt \shortenproofright=0pt \proofbelowshift=0pt
\onleftofproofruletrue\penalty1
}

\def\eproofbit{
\ifx    \wereinproofbit\prooftree
\then   \ifcase \lastpenalty
        \then   \shortenproofright=0pt  
        \or     \unpenalty\hfil         
        \or     \unpenalty\unskip       
        \else   \shortenproofright=0pt  
        \fi
\fi
\global\dimen0=\shortenproofleft
\global\dimen1=\shortenproofright
\global\dimen2=\proofrulebreadth
\global\dimen3=\proofbelowshift
\global\dimen4=\proofdotseparation
\global\count255=\proofdotnumber
$\egroup  
\shortenproofleft=\dimen0
\shortenproofright=\dimen1
\proofrulebreadth=\dimen2
\proofbelowshift=\dimen3
\proofdotseparation=\dimen4
\proofdotnumber=\count255
}

\def\proofover{
\eproofbit 
\setbox\proofbelow=\hbox\bgroup 
\let\wereinproofbit\proofover
$\displaystyle
}%
\def\proofoverdbl{
\eproofbit 
\doubleprooftrue
\setbox\proofbelow=\hbox\bgroup 
\let\wereinproofbit\proofoverdbl
$\displaystyle
}%
\def\proofoverdots{
\eproofbit 
\proofdotstrue
\setbox\proofbelow=\hbox\bgroup 
\let\wereinproofbit\proofoverdots
$\displaystyle
}%
\def\proofusing{
\eproofbit 
\setbox\proofrulename=\hbox\bgroup 
\let\wereinproofbit\proofusing
\kern0.3em$
}

\def\endprooftree{
\eproofbit 
  \dimen5 =0pt
\dimen0=\wd\proofabove \advance\dimen0-\shortenproofleft
\advance\dimen0-\shortenproofright
\dimen1=.5\dimen0 \advance\dimen1-.5\wd\proofbelow
\dimen4=\dimen1
\advance\dimen1\proofbelowshift \advance\dimen4-\proofbelowshift
\ifdim  \dimen1<0pt
\then   \advance\shortenproofleft\dimen1
        \advance\dimen0-\dimen1
        \dimen1=0pt
        \ifdim  \shortenproofleft<0pt
        \then   \setbox\proofabove=\hbox{%
                        \kern-\shortenproofleft\unhbox\proofabove}%
                \shortenproofleft=0pt
        \fi
\fi
\ifdim  \dimen4<0pt
\then   \advance\shortenproofright\dimen4
        \advance\dimen0-\dimen4
        \dimen4=0pt
\fi
\ifdim  \shortenproofright<\wd\proofrulename
\then   \shortenproofright=\wd\proofrulename
\fi
\dimen2=\shortenproofleft \advance\dimen2 by\dimen1
\dimen3=\shortenproofright\advance\dimen3 by\dimen4
\ifproofdots
\then
        \dimen6=\shortenproofleft \advance\dimen6 .5\dimen0
        \setbox1=\vbox to\proofdotseparation{\vss\hbox{$\cdot$}\vss}%
        \setbox0=\hbox{%
                \advance\dimen6-.5\wd1
                \kern\dimen6
                $\vcenter to\proofdotnumber\proofdotseparation
                        {\leaders\box1\vfill}$%
                \unhbox\proofrulename}%
\else   \dimen6=\fontdimen22\the\textfont2 
        \dimen7=\dimen6
        \advance\dimen6by.5\proofrulebreadth
        \advance\dimen7by-.5\proofrulebreadth
        \setbox0=\hbox{%
                \kern\shortenproofleft
                \ifdoubleproof
                \then   \hbox to\dimen0{%
                        $\mathsurround0pt\mathord=\mkern-6mu%
                        \cleaders\hbox{$\mkern-2mu=\mkern-2mu$}\hfill
                        \mkern-6mu\mathord=$}%
                \else   \vrule height\dimen6 depth-\dimen7 width\dimen0
                \fi
                \unhbox\proofrulename}%
        \ht0=\dimen6 \dp0=-\dimen7
\fi
\let\doll\relax
\ifwasinsideprooftree
\then   \let\VBOX\vbox
\else   \ifmmode\else$\let\doll=$\fi
        \let\VBOX\vcenter
\fi
\VBOX   {\baselineskip\proofrulebaseline \lineskip.2ex
        \expandafter\lineskiplimit\ifproofdots0ex\else-0.6ex\fi
        \hbox   spread\dimen5   {\hfi\unhbox\proofabove\hfi}%
        \hbox{\box0}%
        \hbox   {\kern\dimen2 \box\proofbelow}}\doll%
\global\dimen2=\dimen2
\global\dimen3=\dimen3
\egroup 
\ifonleftofproofrule
\then   \shortenproofleft=\dimen2
\fi
\shortenproofright=\dimen3
\onleftofproofrulefalse
\ifinsideprooftree
\then   \hskip.5em plus 1fil \penalty2
\fi
}

\newcommand{\op}{{o}}

\renewcommand{\to}{\xrightarrow{}}
\newcommand{\ot}{\xleftarrow{}}
\newcommand{\tto}[1]{\xrightarrow{#1}}
\newcommand{\oot}[1]{\xleftarrow{#1}}

\newcommand{\mono}{\rightarrowtail}

\newcommand{\epi}{\twoheadrightarrow}

\newcommand{\inclusion}{\hookrightarrow}
\newcommand{\mmono}[1]{\stackrel{#1}\rightarrowtail}
\newcommand{\oonom}[1]{\stackrel{#1}\leftarrowtail}
\newcommand{\eepi}[1]{\stackrel{#1}\twoheadrightarrow}
\newcommand{\iipe}[1]{\stackrel{#1}\twoheadleftarrow}

\renewcommand{\mapsfrom}{\mathrel{\reflectbox{\ensuremath{\mapsto}}}}

\newcommand{\TTo}[3]{\begin{tikzar}[sep=small]
{#1}\arrow[Rightarrow]{r}{#2}\& {#3}
\end{tikzar}}

\newcommand{\Set}{\mathsf{Set}}

\newcommand{\Cat}{\mathsf{Cat}}

\newcommand{\CAT}{\mathsf{CAT}}

\newcommand{\id}{\mathrm{id}}
\newcommand{\Dom}{\mathrm{Dom}}
\newcommand{\Cod}{\mathrm{Cod}}
\newcommand{\Ids}{\mathrm{Ids}}
\newcommand{\Id}{\mathrm{Id}}

\newcommand{\AAA}{{\cal A}}
\newcommand{\BBB}{{\cal B}}

\newcommand{\DDD}{{\cal D}}

\newcommand{\GGG}{{\cal G}}

\newcommand{\JJJ}{{\cal J}}
\newcommand{\KKK}{{\cal K}}
\newcommand{\LLL}{{\cal L}}

\newcommand{\NNN}{{\cal N}}

\newcommand{\RRR}{{\cal R}}
\newcommand{\SSS}{{\cal S}}

\newcommand{\XXX}{{\cal X}}
\newcommand{\YYY}{{\cal Y}}

\renewcommand{\Bbb}{\mathbb}
\newcommand{\AAa}{{\Bbb A}}
\newcommand{\BBb}{{\Bbb B}}

\newcommand{\CCc}{{\cal C}}
\newcommand{\CCcc}{{\Bbb C}}

\newcommand{\DDd}{{\Bbb D}}
\newcommand{\EEe}{{\Bbb E}}

\newcommand{\GGg}{{\Bbb G}}

\newcommand{\LLl}{{\Bbb L}}
\newcommand{\MMm}{{\Bbb M}}

\newcommand{\PPp}{{\Bbb P}}
\newcommand{\QQq}{{\Bbb Q}}
\newcommand{\RRr}{{\Bbb R}}

\newcommand{\TTt}{{\Bbb T}}

\newcommand{\XXx}{{\Bbb X}}
\newcommand{\YYy}{{\Bbb Y}}
\newcommand{\ZZz}{{\Bbb Z}}

\newcommand{\LlL}{{\mathfrak L}}

\newcommand{\mathbold}[1]{\mbox{\boldmath $#1$}}

\renewcommand{\ggg}{{\mathbold g}}
\newcommand{\hhh}{{\mathbold h}}

\mathcode`\<="4268 
\mathcode`\>="5269 
\mathchardef\gt="313E 
\mathchardef\lt="313C 

  \def\pushright#1{{
     \parfillskip=0pt            
     \widowpenalty=10000         
     \displaywidowpenalty=10000  
     \finalhyphendemerits=0      
     \leavevmode                 
     \unskip                     
     \nobreak                    
     \hfil                       
     \penalty50                  
     \hskip.2em                  
     \null                       
     \hfill                      
     {#1}                        
     \par}}                      
 
  \def\qed{\pushright{$\square$}\penalty-700 \smallskip}

 \newenvironment{prf}[1]{\begin{trivlist} \item[{\bf \;\,Proof}#1.]}%
 {\qed\end{trivlist}}

\newcommand{\beq}{\begin{equation}}
\newcommand{\eeq}{\end{equation}}
\newcommand{\ba}[1]{\begin{array}{#1}}
\newcommand{\ea}{\end{array}}
\newcommand{\bea}{\begin{eqnarray}}
\newcommand{\eea}{\end{eqnarray}}
\newcommand{\bear}{\begin{eqnarray*}}
\newcommand{\eear}{\end{eqnarray*}}
\newcommand{\bpr}{\begin{prf}{}}
\newcommand{\epr}{\end{prf}}
\newcommand{\bprf}[1]{\begin{prf}{#1}}
\newcommand{\eprf}{\end{prf}}

\newcommand{\Isb}[1]{\JJJ{#1}}
\newcommand{\Dfib}[1]{\comm{\mathsf{Dfib}}{#1}}
\newcommand{\Ofib}[1]{\comm{\mathsf{Ofib}}{#1}}
\newcommand{\Labs}{\UD}

\newcommand{\lft}[1]{\overleftarrow{#1}}
\newcommand{\rgt}[1]{\overrightarrow{#1}}

\newcommand{\supp}{\mathop{\underrightarrow{\mathrm{lim}}}}
\newcommand{\inff}{\mathop{\underleftarrow{\mathrm{lim}}}}
\newcommand{\liminff}{\mathop{\overleftarrow{\mathrm{lim}}}}
\newcommand{\limsupp}{\mathop{\overrightarrow{\mathrm{lim}}}}

\newcommand{\Dashv}{\mathrel{\mbox{$={\kern-1.1ex}|$}}}
\renewcommand{\vDash}{\mathrel{\mbox{$|{\kern-1.1ex}=$}}}

\newcommand{\mmult}{\odot}

\renewcommand{\Diamond}{\mbox{\large$\lozenge$}}

\newcommand{\Dabla}{\mbox{\begin{sideways} $\bowtie$\end{sideways}}}
\newcommand{\cmn}{\Delta}
\newcommand{\mnd}{\nabla}
\newcommand{\cmnu}{\blacktriangle}\newcommand{\mndu}{\blacktriangledown}

\newcommand{\Up}{{\Uparrow}}
\newcommand{\pU}[1]{{#1}{\Uparrow}}
\newcommand{\Do}{{\Downarrow}}
\newcommand{\UD}{{\Updownarrow}}

\newcommand{\Con}{{\sf Con}}
\newcommand{\Coc}{{\sf Coc}}

\newcommand{\Lan}{\underline{\cmn}}
\newcommand{\Ran}{\overline{\mnd}}
\newcommand{\dLan}{\underline{\vartriangle}}
\newcommand{\dRan}{\overline{\triangledown}}
\newcommand{\nRan}{\underline{\mnd}}
\newcommand{\nLan}{\overline{\cmn}}
\newcommand{\nnLan}{\overline{\underline{\cmn}}}
\newcommand{\nnRan}{\underline{\overline{\mnd}}}
\newcommand{\Kan}{\Lan\dashv\Ran}
\newcommand{\dKan}{\dLan\dashv\dRan}
\newcommand{\nKan}{\nRan\dashv\nLan}
\newcommand{\RLan}{\boxslash}
\newcommand{\LRan}{\boxbslash}
\newcommand{\RLAN}{\underline{\overline{\boxslash}}}
\newcommand{\LRAN}{\underline{\overline{\boxbslash}}}
\newcommand{\DRL}{\Do\CCc^\RLan}
\newcommand{\ULR}{\Up\CCc^\LRan}
\newcommand{\simDo}{\overline{\Do\CCc}}
\newcommand{\simsimDo}{\overline{\Do{\simDo}}}
\newcommand{\simUp}{\overline{\Up\CCc}}

\newcommand{\compson}{H}

\newcommand{\ladj}[1]{{#1}^\ast}
\newcommand{\radj}[1]{{#1}_\ast}
\newcommand{\adj}[1]{\ladj {#1} \dashv \radj {#1}}

\newcommand{\arrow}[1]{\comm{#1}{#1}}
\newcommand{\tarrow}[1]{\comm{#1}{#1^\op}}

\newcommand{\presheaf}{left {\action}}
\newcommand{\presheaves}{left {\action}s}

\newcommand{\action}{action}
\newcommand{\actions}{actions}
\newcommand{\Actions}{Actions}
\newcommand{\postsheaf}{right {\action}}
\newcommand{\postsheaves}{right {\actions}}

\newcommand{\aggg}{\alpha}

\newcommand{\bggg}{\beta}

\newcommand{\comm}[2]{{#1}\diagup{#2}}
\newcommand{\commm}[1]{\big[{#1}\big]}

\newcommand{\pbdown}{\mathbin{\rotatebox[origin=c]{135}{\Large$\ulcorner$}}}

  \newcommand{\isbell}[1]{\JJJ{#1}}   
   
  \newcommand{\gap}{\mathsf{g}}  
   \newcommand{\gapp}{\mathsf{J}}   
  \newcommand{\intv}{\jmath}
  \newcommand{\intvl}{\mathsf{j}}
  
  \newcommand{\Gaps}{\mathsf{Gaps}} 
   \newcommand{\Intvls}{\mathsf{Ints}} 
\newcommand{\Cuts}{\mathsf{Cuts}} 
\newcommand{\cuts}{\mathsf{cuts}} 
\newcommand{\ida}{\mathsf{a}} 
\newcommand{\idb}{\mathsf{b}} 
\newcommand{\idc}{\mathsf{c}} 
\newcommand{\idd}{\mathsf{d}}

\newcommand{\idepsilon}{\mathsf{e}} 
\newcommand{\ideta}{\mathsf{h}}

 \newcommand{\rrrr}{r}
 \newcommand{\dddd}{d}

\newcommand{\GAP}[3]{\TTo{{#1}\Big)}{#2}{\Big({#3}}}
\newcommand{\INTVL}[3]{\TTo{\Big({#1}}{#2}{{#3\Big)}}}

\newcommand{\sseq}[1]{\left(\, #1\, \right)}

\newcommand{\pseq}[1]{\big[\, #1\, \big]}
\newcommand{\psseq}[1]{\left[\, #1\, \right]}

\graphicspath{{PIC/}}

\title{Tight limits and completions\\
 from Dedekind-MacNeille to Lambek-Isbell
}
\author[1]{Dusko Pavlovic}
\author[2]{Dominic J.\,D.\,Hughes}
\affil[1]{University of Hawaii, Honolulu HI and Radboud University, Nijmegen, \ 
  \texttt{dusko@hawaii.edu}\thanks{Supported by NSF and AFOSR.}}
  \affil[2]{Apple Inc., Cupertino CA  and UC Berkeley, Berkeley CA\thanks{Visiting scholar, Logic Group. Many thanks to Wes Holliday and Dana Scott for inviting me.}}
\begin{document}
\date{}
\maketitle

\begin{flushright}
\parbox{12cm}{\small \it It is an open problem whether there exists a sup- and inf-complete category $\AAa''''$ with a sup- and inf-dense embedding $\AAa \to  \AAa''''$ in analogy to the Dedekind completion of an ordered set.}\\[1ex] 
\footnotesize Joachim Lambek \cite[Introduction]{LambekJ:completions}
\\[2ex]
\parbox{12cm}{\small \it No Lambek extension [i.e., a sup- and inf-dense embedding]  of the one-object category $\ZZz_4$ has finite limits [and cannot be inf-complete].}
\\[1ex]
{\footnotesize John Isbell \cite[Thm. 3.1]{IsbellJ:no-lambek}}
\end{flushright}

\begin{abstract}
While any infimum in a poset can also be computed as a supremum, and vice versa, categorical limits and colimits do not always approximate each other. 
If I approach a point from below, and you approach it from above, then we will surely meet if we live in a poset, but we may miss each other in a category. Can we characterize the limits and the colimits that approximate each other, and guarantee that we will meet?  Such limits and colimits are called \emph{tight}. Some critically important network applications depend on them. This paper characterizes tight limits and colimits, and describes tight completions, derived by applying the nucleus construction to adjunctions between loose completions. Just as the Dedekind-MacNeille completion of a poset preserves any existing infima and suprema, the tight completion of a category preserves any existing tight limits and colimits and is therefore idempotent.  
\end{abstract}

\newpage
\tableofcontents
\clearpage

\section*{Preliminaries}
\addcontentsline{toc}{section}{Preliminaries}

We made an effort to make this paper accessible not only to readers familiar with categories, but also to those who may not be, but come with enough interest, patience, and confidence. For the latter group, we include some clarifications of the basic concepts in the Appendix. For the former group, the Appendix should clarify any notation or concepts that may seem standard to some and nonstandard to others, and in any case do not lie on the main path of the narrative. Since the main path remains lengthy even with such deferrals, some of the routine introductory material is postponed.

\section{Background: From real numbers to tight poset completions}\label{Sec:Back}
\subsection{Loose vs tight completions}\label{Sec:pos}

A real number is an infinite object approximated arbitrarily closely by finite objects, the rational numbers. An infinite ordinal has a finitary normal form which expresses a finite descent down the tree of operations that generate it. A finite set of grammatical rules generates infinitely many sentences from a finite set of words. Finite programs describe  infinite computations. Infinite objects or processes have finite descriptions. 

A real number can be approximated by rational numbers from below and from above: it has a \emph{tight}\/ approximation. The length of the shortest program for a given algorithm can be approximated from above, but not from below. An interval can be approximated from below and from above, and that approximation is not tight.

The elements of a complete lattice can always be approximated tightly, as every supremum is also an infimum, and vice versa. A lattice has all suprema if and only if it has all infima. A complete category has all limits, but may not have all colimits.

\paragraph{Completions} pervade logics and mathematics. There are completions of partially ordered sets, metric spaces, measures, uniformities, categories, logical theories, knowledge bases, and many other structures. Intuitively, a completion removes some sort of incompleteness, or "holes" left behind a construction or specification. For example, we specify the natural numbers to count things, and the rational numbers to be able compare and measure them. The comparisons impose the order of rational numbers, usually denoted by $\QQq$. But some comparisons point to "holes" in this order. For example, the comparison of the length of a circle with its diameter comes above $3\frac{1}{8}$, but below $3\frac{1}{7}$. Archimedes' measurements and calculations brought the lower bound up to $3\frac{10}{71}$, and the upper bound down to $3\frac{177}{1250}$ \cite{Archimedes:circle}.  But no matter how much we refine the fractions, a hole between the lower and the upper bounds of the ratio  remains. Recognizing that the number $\pi$, expressing the length of a circle in terms of its diameter, is a persistent hole in the order of rational numbers, took many centuries, from the Pythagoreans to J.H.~Lambert in the XVIII century \cite{Pi-book}. We can approximate $\pi$ as closely as we wish, for example, by summing up an alternating series like
\beq\label{eq:pi} 4-\frac 4 3 + \frac 4 5 - \frac 4 7 +\frac 4 9 - \cdots\eeq 
If we stop after a positive summand, we end up above $\pi$, while if we stop at a negative summand, we end up below. There is always a hole for $\pi$ in the order $\QQq$ of rational numbers. Filling all such holes produces the continuous order of real numbers $\RRr$. There are many ways to fill the holes  \cite{PavlovicD:retracing,PavlovicD:CRN,PavlovicD:CRN1}. One way is by imagining that we can go infinitely long along the sequences of sums like \eqref{eq:pi} such that the distance between the partial sums become infinitely small. That way the continuum $\RRr$ arises as the \emph{metric completion}\/ of $\QQq$. Another way to fill the holes in $\QQq$ is to notice that each hole is the least upper bound, or \emph{supremum}, of the rationals below it, and the greatest lower bound, or \emph{infimum}, of the rationals above it, and to construct $\RRr$ by adjoining to $\QQq$ all suprema and infima. This way, the real continuum $\RRr$ arises as the \emph{lattice completion}\/ of the order $\QQq$.  This is how Richard Dedekind reconstructed it in his seminal memoir \cite{DedekindR:zahlen}. 

\paragraph{Formalizing the incompleteness and deriving completions.} More formally, the incompleteness gaps or "holes" in a structure arise when the outputs of some operations, canonically definable on a structure, have not been defined, and the "holes" in the structure that would be filled by their outputs are left empty. For example, the infima (meets, the least upper bounds) and the suprema (joins, the greatest lower bounds) are canonically definable in any poset. A poset is incomplete when it does not contain the suprema or the infima of all of its subsets. The poset thus has "holes" or gaps in the places where some suprema or infima would be. 

It is precisely definable where the suprema or the infima would be, but they are not there. To fill the "holes" for the suprema, we take the largest of sets that would have the same least upper bound and adjoin it to the structure. The largest of all sets that would have the same join (least upper bound) must be downward closed. The join-completion is thus the set of all lower-closed sets.

\subsection{Loose completions of posets}\label{Sec:pos-loose}

A loose completion of a poset $\PPp$ adjoins 
\begin{itemize}
\item the least upper bounds or \emph{suprema} $\bigvee$, or
\item the greatest lower bounds or \emph{infima} $\bigwedge$.
\end{itemize}
The constructions of the loose completions of posets are based on the observation that every set $S$ of elements of a poset $\PPp$ has the same upper bounds, and therefore the same supremum, as its lower closure $\lft S$ and the same lower bounds, and therefore the same infimum, as its upper closure $\rgt S$. In summary, the sets
\beq\label{eq:lowcl}
\begin{split}
\prooftree
S\subseteq \PPp
\justifies
\lft S = \{x\in \PPp|\exists s\in S.\ x\leq s\}\qquad \qquad \qquad \rgt S = \{x\in \PPp|\exists s\in S.\ s\leq x\}
\endprooftree
\end{split}
\eeq
satisfy  
\beq\label{eq:sup-reduction}
\bigvee S = \bigvee \lft S \qquad \qquad \qquad \qquad  \bigwedge S = \bigwedge \rgt S
\eeq
Therefore, the lower sets form the sup-completion, with each lower set as its own supremum, whereas the upper sets form the inf-completion. More precisely, consider the sets
\beq\label{eq:douppos}
\Do \PPp  =  \{\lft S\subseteq \PPp\ |\ x\leq y \wedge  \lft S(y) \implies \lft S(x)\}\qquad\quad
\Up \PPp  =  \{\rgt S\subseteq \PPp\ |\ \rgt S(x) \wedge x\leq y \implies \rgt S(y)\}
\eeq
where we write $S(x)$ instead of $x\in S$.  The suprema and the infima of $\XXX\subseteq \Do\PPp$ are 
\beq\label{eq:supremum}
\bigvee \XXX  =  \bigcup_{\lft X \in \XXX} \lft X\qquad\qquad\qquad\qquad  \bigwedge \XXX  =  \bigcap_{\lft X \in \XXX} \lft X
\eeq 
whereas the suprema and the infima of $\YYY\subseteq \Up\PPp$ are
\beq\label{eq:infimum}
\bigvee \YYY  =  \bigcap_{\rgt Y \in \YYY} \rgt Y\qquad\qquad\qquad\qquad  \bigwedge \YYY  =  \bigcup_{\rgt Y \in \YYY} \rgt Y
\eeq 
The lower sets in $\Do\PPp$ are thus ordered by set inclusion $\subseteq$, whereas the upper sets in $\Up\PPp$ are ordered by the \emph{opposite}\/ set inclusion $\supseteq$. The maps
\begin{align*}
\mnd: \PPp & \to \Do \PPp & \cmn :\PPp & \to \Up \PPp\\
x & \mapsto \{y\leq x\} & x &\mapsto \{y\geq x\}
\end{align*}
are thus monotone embeddings, in the sense that
\beq
\mnd x\subseteq \mnd y \ \ \iff \ \ x\leq y \ \ \iff\ \ \cmn x \supseteq \cmn y
\eeq  
They formally make $\Do\PPp$ into the sup-completion and $\Up\PPp$ into the inf-completion because every monotone map $f:\PPp\to \LLl$ into a complete lattice $\LLl$ induces a unique sup-preserving map $f_\vee$ and a unique inf-preserving map $f_\wedge$ making the following diagrams commute
\beq\label{eq:univ-pos}
\begin{tikzar}[row sep = 2ex]
\& \Do\PPp\ar{dd}[description]{\exists ! f_\vee}\\
\PPp\ar{ur}{\mnd}\ar{dr}[description]{\forall f}\\
\& \LLl
\end{tikzar}\qquad \qquad
\begin{tikzar}[row sep = 2ex]
\& \Up\PPp\ar{dd}[description]{\exists! f_\wedge}\\
\PPp\ar{ur}{\cmn}\ar{dr}[description]{\forall f}\\ 
\& \LLl
\end{tikzar}
\eeq
This follows from the fact that $\Do\PPp$ is $\vee$-generated by $\PPp$ through $\mnd$, whereas $\Up\PPp$ is $\wedge$-generated by $\PPp$ through $\cmn$, in the sense that every $\lft A\in \Do\PPp$ and every $\rgt B \in \Up\PPp$ satisfies
\beq\label{eq:infimum:supremum}
\lft A  = \bigcup_{x\in\lft A} \mnd x = \bigvee_{x\in\lft A} \mnd x \qquad\qquad\qquad\qquad  \rgt B = \bigcup_{y\in\rgt B} \cmn y  = \bigwedge_{y\in\rgt B} \cmn y 
\eeq 
While any infimum $x\wedge y$ that may exist in $\PPp$ is preserved under $\mnd$ and any supremum $x\wedge y$ is preserved under $\cmn$
\beq
\mnd(x\wedge y) = \mnd x \cap \mnd y \qquad \qquad \qquad \cmn(x\vee y) = \cmn x \cap \cmn y
\eeq
the sup-completion $\mnd$ does not preserve the suprema, and the inf-completion $\cmn$ does not preserve the infima, except of $x\leq y$, i.e. in the form $x\vee y = y$ and $x\wedge y = x$.


\subsection{Loose completions create gaps}
\label{Sec:loose-gaps}

The trouble with this simple and elegant approach is that adds the new suprema but destroys all of the old ones. Applying the above approach to the poset $\QQq$ or rational number gives the sup-completion $\Do\QQq$ containing the intervals 
\bea\label{eq:intervals}
\left(-\infty, \frac 1 2\right] & = &  \left\{x\in \QQq\ |\ x\leq\frac 1 2\right\}\ \  =\ \   \mnd \frac 1 2 \\ 
\left(-\infty, \frac 1 2\right) &= & \left\{x\in \QQq\ |\ x\lt \frac 1 2\right\}\ \ =\ \ \bigvee_{y\lt \frac 1 2} \mnd y\notag
\eea
The first one is the image of $\frac 1 2 \in \QQq$ along $\mnd$. The second one is the freshly adjoined supremum of the rationals below $\frac 1 2$. In $\QQq$, the supremum of the set of rationals below $\frac 1 2$ was $\frac 1 2$. In $\Do\QQq$, the supremum of the set of rationals below $\frac 1 2$ is the set $(-\infty, \frac 1 2)$ itself, whereas the image of $\frac 1 2\in\QQq$ along $\mnd$ is $(-\infty, \frac 1 2]$. But the two intervals in \eqref{eq:intervals} are different, and $\frac 1 2$ as the supremum of the numbers below $\frac 1 2$ is not preserved, but it is replaced by a new supremum. We are thus left with two copies of $\frac 1 2$: the old one $\mnd \frac 1 2$, and the new one $\bigvee_{y\lt \frac 1 2} \mnd y$. In-between the two there is a hole, because the set of upper bounds of both intervals in \eqref{eq:intervals} is $[\frac 1 2, \infty)$, and only the upper one intersects it. The same holds for any other rational number. Together with the irrational numbers, adjoined as suprema of rationals, the completion $\Do \QQq$ has thus adjoined a new copy of the rationals, as suprema of the rationals below them. It contains two copies of $\QQq$, with a new hole in-between each couple.  The same happens in $\Up \QQq$, except that the new rationals and the holes are just above the old rationals. They are an example of the \emph{loose}\/ completions posets, where all new suprema or infima are created, and the old ones are destroyed.

\subsection{Tight completions of posets} 

The crucial observation that enabled Dedekind to fill the holes was that every point of the continuum is \emph{both}\/ a supremum of the  rational points below it \emph{and also}\/ an infimum of the rational points above it. Adjoining to the order of rational numbers just the points with that property fills all holes in it. That was achieved by \emph{Dedekind's cuts}. The completions that adjoin only the suprema that are also infima, and only the infima that are also suprema, are called the \emph{tight}\/ completions. 

For simplicity, we describe Dedekind's tight completion in modern terms and for an abstract poset $\PPp$ rather than just the rational order $\QQq$. The observation that Dedekind's construction goes through for arbitrary posets and its general development are due to Holbrook MacNeille \cite{MacNeille}. The following diagram provides an overview of the construction.
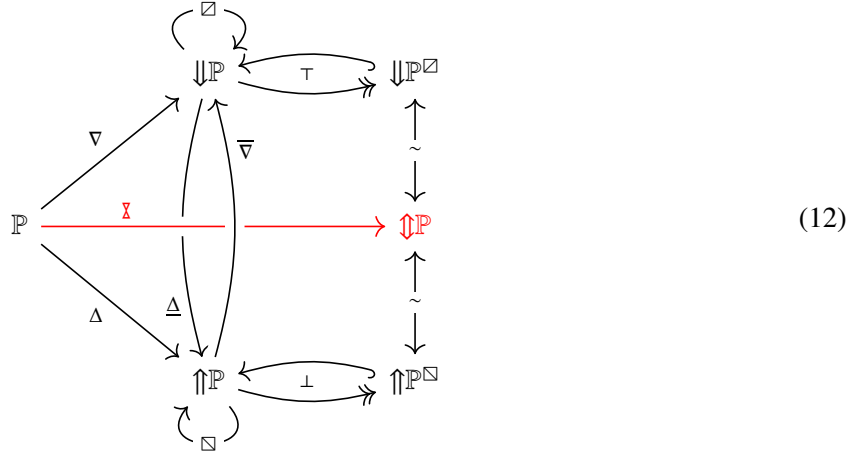
\beq\label{eq:DM}
\begin{tikzar}{}
\&\& \Do \PPp 
\arrow[loop, out = 135, in = 45, looseness = 4]{}[description]{\RLan}
\arrow[bend right = 15,near end]{dddd}[swap]{\Lan} \arrow[phantom]{dddd}{\dashv} 
\arrow[two heads,bend right = 15]{rr}
\arrow[phantom]{rr}{\scriptstyle\top}
 \&\& \Do\PPp^{\RLan} 
\arrow[leftrightarrow]{dd}[description]{\sim}
\arrow[hookrightarrow,bend right = 15]{ll} 
\\
\\
\PPp \arrow{uurr}{\mnd} \arrow{ddrr}[swap]{\cmn} 
\arrow[crossing over,color=red]{rrrr}[near start]{\color{red}\Dabla} 
\&\&\&\&{\color{red}\UD\PPp} \arrow[leftrightarrow]{dd}[description]{\sim}
\\
\\ 
\&\& \Up \PPp 
\arrow[loop, out = -50, in=-130, looseness = 4]{}[description]{\LRan}
\arrow[bend right = 15,near end,crossing over]{uuuu}[swap]{\Ran}  \arrow[two heads,bend right = 15]{rr} \arrow[phantom]{rr}{\scriptstyle\bot} \&\& \Up\PPp^{\LRan} 
\arrow[hookrightarrow,bend right = 15]{ll} 
\end{tikzar}
\eeq
The tight completion $\Dabla:\PPp\to \UD\PPp$ is the universal quotient of the two loose completions, $\mnd:\PPp\to \Do\PPp$ and $\cmn:\PPp\to \Up\PPp$, in which the suprema from $\Do\PPp$ and the infima $\Up\PPp$ are expressible in $\UD\PPp$ in terms of each other, and they generate each element of $\UD\PPp$ both from below and from above, tightly. The common quotient of $\Do\PPp$ and $\Up\PPp$ is constructed using their universal properties \eqref{eq:univ-pos}. Instantiating $f$ in the left-hand triangle of \eqref{eq:univ-pos} to $\cmn:\PPp \to \Up\PPp$ gives a sup-preserving map $\Lan:\Do\PPp \to \Up \PPp$. Instantiating it in the right-hand triangle of \eqref{eq:univ-pos} to $\mnd:\PPp \to \Do\PPp$ gives a meet-preserving map $\Ran:\Do\PPp \to \Up \PPp$. Since $\Up\PPp$ is ordered by reverse inclusion $\supseteq$, both maps are antitone with respect to the inclusion, and the unions are mapped to the intersections:
\beq\label{eq:pKan}
\begin{tikzar}[column sep=4em]
\lft S = \displaystyle \bigcup_{x\in \lft S} \mnd x  \arrow[thin,mapsto]{dd}
\& \Do \PPp\arrow[bend right = 15]{dd}[swap]{\Lan} \arrow[phantom]{dd}{\dashv} \& \displaystyle\Ran \rgt S = \bigcap_{x\in \rgt S} \mnd x  
\\
\\ 
\displaystyle\Lan\lft S = \bigcap_{x\in \lft S}\cmn x \&\Up \PPp 
\arrow[bend right = 15]{uu}[swap]{\Ran}  \& \displaystyle\rgt S = \bigcup_{x\in \rgt S} \cmn x \arrow[thin,mapsto]{uu}
\end{tikzar}
\eeq
A closer look shows that $\Lan \lft S$ is the set of upper bounds of the lower set $\lft S$, whereas $\Ran \rgt S$ is the set of lower bounds of the upper set $\rgt S$. The adjunction (Galois connection) $\Kan$, i.e., the logical equivalence 
\bea
\Lan \lft S \subseteq \rgt S & \iff & \lft S \supseteq \Ran \rgt S
\eea
states that the upper set $\rgt S$ contains all upper bounds of the lower set $\lft S$ if and only if $\lft S$ contains all lower bounds of $\rgt S$. The composites 
\beq\label{eq:closureops}
\RLan = \Ran\Lan:\Do\PPp\to \Do\PPp \qquad \mbox{ and }\qquad \LRan = \Lan\Ran:\Up\PPp\to \Up\PPp
\eeq
are closure operators as they satisfy  $\lft S \subseteq \RLan \lft S =\RLan\RLan \lft S$ and $\rgt S\subseteq \LRan \rgt S =\LRan\LRan \rgt S$ for all  $\lft S\in \Do\PPp$ and $\rgt S\in \Up\PPp$. While $\RLan \lft S$ is the set of all lower bounds of the set of upper bounds of $\lft S$, $\LRan \rgt S$ is the set of all upper bounds of the set of lower bounds of $\rgt S$. The lattices of closed sets 
\beq \label{eq:closed}
\Do\PPp^\RLan  =  \left\{\lft S\in \Do\PPp\ |\ \lft S = \RLan \lft S\right\}\qquad\qquad\quad
\Up\PPp^\LRan  =  \left\{\rgt S\in \Up\PPp\ |\ \rgt S = \LRan \rgt S\right\}
\eeq
are isomorphic. The easiest way to see this is to consider the lattice
\bea\label{eq:cuts}
\UD \PPp & = & \left\{<\lft S, \rgt S> \in \Do\PPp\times \Up \PPp\ |\ \lft S = \Ran \rgt S\  \wedge\  \Lan \lft S = \rgt S\right\}
\eea
The pairs $<\lft S, \rgt S>$ of sets of each other's bounds are the general form of \emph{Dedekind cuts}. The lower and upper sets $\lft S$ and $\rgt S$ obviously determine each other, and the isomorphisms
\beq\label{eq:three}
\Do\PPp^\RLan\ \cong \ \UD \PPp\ \cong\  \Up\PPp^\LRan 
\eeq
are straightforward. Instantiating $\PPp$ to the rational numbers $\QQq$ yields Dedekind's reals $\RRr = \UD\QQq$. As retracts of the complete lattices $\Do\PPp$ and $\Up\PPp$, the lattices in \eqref{eq:three} are complete, but their suprema and infima are different. In contrast with the loose suprema and infima 
\beq\label{eq:loose-infsup}
\bigvee \XXX  = \bigcup_{\lft S \in \lft \XXX} \lft S\qquad\qquad\qquad\qquad \bigwedge \YYY  =  \bigcup_{\rgt S \in \rgt \YYY} \rgt S
\eeq
from \eqref{eq:supremum} and \eqref{eq:infimum}, taking $\lft \AAA\subseteq \Do \PPp^\RLan$ and $\rgt \BBB\subseteq \Up\PPp^\LRan$ we have
\beq\label{eq:tight-infsup}
\bigvee \lft \AAA  = \Ran\bigcap_{\lft A \in \lft \AAA} \Lan \lft A\qquad\qquad\qquad\qquad \bigwedge \rgt \BBB  =  \Lan\bigcap_{\rgt B\in \rgt \BBB} \Ran \rgt B
\eeq
The suprema and the infima of a set of cuts $\Upsilon\subseteq \UD \PPp$ are now simply
\beq\label{eq:tight-cuts}
\bigvee \Upsilon  = \left< \Ran\bigcap_{<\lft C, \rgt C> \in \Upsilon} \rgt C\, , \bigcap_{<\lft C, \rgt C> \in \Upsilon} \rgt C \right> \qquad\quad \bigwedge \Upsilon  =  \left<\bigcap_{<\lft C,\rgt C>\in \Upsilon} \lft C\, ,\  \Lan\bigcap_{<\lft C,\rgt C>\in \Upsilon} \lft C\right>
\eeq
since each $<\lft C,\rgt C>\in \Upsilon$ satisfies $\lft C = \Ran\rgt C$ and $\rgt C =\Lan \lft C$. 
$\UD\PPp$-suprema are thus constructed as intersections in the inf-completion $\Up\PPp$, while the $\UD\PPp$-infima are constructed as intersections in the inf-completion $\Up\PPp$. 

\section{Problem: Tight completions of categories}\label{Sec:Prob}

How categories generalize posets and preorders and how categorical limits and colimits generalize posetal infima and suprema are summarized in Appendices~\ref{Appendix:basic} and \ref{Appendix:lim}. Here we only introduce notation.

\subsection{Loose completions of categories}\label{Sec:cat-loose}

\subsubsection{From lower and upper sets to left and right \actions}
To lift \eqref{eq:lowcl} from posets to categories, the lower and the upper closures of an arbitrary subset are first generalized to the comma categories $\comm \CCc D$ and $\comm D \CCc$ (as described in Appendix~\ref{Appendix:comma} or \cite[Sec.~II.6]{MacLaneS:CWM}) induced by an arbitrary functor $D:\DDd\to \CCc$. These comma categories are then reduced to the categories $\lft \DDd$ and $\rgt \DDd$ of their connected components (as described in Appendix~\ref{Appendix:diags}). The categorical version of \eqref{eq:lowcl} is thus 
\beq\label{eq:commaD}
\begin{split}
\prooftree
\prooftree
D:\DDd\to \CCc
\justifies
\Dom: \comm \CCc D\to \CCc\qquad \qquad \qquad \qquad \Cod: \comm D \CCc\to \CCc 
\endprooftree
\justifies
\rule{0ex}{2.6ex}
\lft D: \lft \DDd\to \CCc\qquad \qquad \qquad \qquad \rgt D: \rgt\DDd\to \CCc
\endprooftree
\end{split}
\eeq
The reductions \eqref{eq:sup-reduction} now lift to
\beq\label{eq:lim-reduction}
\supp D = \supp \lft D \qquad \qquad \qquad \qquad  \inff D = \inff \rgt D
\eeq
%
%
Just as the lattices of lower and of upper sets $\Do\PPp$ and $\Up\PPp$ from \eqref{eq:douppos} are respectively the $\vee$-completion and the $\wedge$-completion of a preorder $\PPp$, the categories of left \actions\ and of right \actions
\beq\label{eq:doupcat}
\Do \CCc  =  \left\{A:|\CCc|\to \Set\ |\ \CCc(x,y) \times   A_y \tto\ast A_x\right\}\quad
\Up \CCc  =  \left\{B:|\CCc|\to \Set\ |\ B_x \times  \CCc(x,y)  \tto ! B_y\right\}
\eeq
are respectively the $\supp$-completion and the $\inff$-completion of a category $\CCc$. The morphisms in these categories are the equivariant indexed functions, as described in Appendix~\ref{Appendix:basic}. The equivariance conditions assure that the indexed sets\footnote{Categorical conditions for indexing (essentially small) sets over (possibly large) categories are discussed in \cite{Freyd-Street:size,StreetR:size}, for example.} $A$ and $B$ can be viewed as functors $A:\CCc^\op\to \Set$ and $B:\CCc\to \Set$, or comprehended as functors $\lft A:\lft \AAa\to\CCc$ and $\rgt B:\rgt \BBb\to \CCc$ from the categories defined by
\begin{align}\label{eq:groth}
\Big\lvert\lft \AAa\Big\rvert &= \coprod_{x\in \CCc} A_x & \lft\AAa\left(s_x, t_y\right) & = \left\{f\in \CCc(x,y)\ |\ s_x=f\ast t_y\right\}\\
\Big\lvert\rgt \BBb\Big\rvert &= \coprod_{x\in \CCc} B_x & \rgt\BBb\left(u_x, v_y\right) & = \left\{f\in \CCc(x,y)\ |\ u_x!f = v_y\right\}\label{eq:groth-cov}
\end{align}
with the obvious projections $\lft A$ and $\lft B$. Such projections are recognized by the fact that they induce bijections $\comm{\AAa}{s} \cong \comm{\CCc}{\lft As}$ and $\comm{u}{\BBb} \cong \comm{\rgt Bu}{\CCc}$ for all $s\in \AAa$ and $u\in \BBb$. The functors that induce such bijections are called \emph{discrete fibrations} and \emph{discrete opfibrations} \cite[Vol.~2, Ch.~8]{BorceuxF:handbook}. Respectively, they form categories $\Dfib\CCc$ and $\Ofib\CCc$, with commutative triangles of functors as the morphisms. Constructions (\ref{eq:groth}--\ref{eq:groth-cov}) induce the equivalences
\beq\label{eq:DopUUUp}
\Dfib \CCc \simeq \Set^{\CCc^\op}\qquad \qquad\qquad\qquad \Ofib\CCc \simeq \Set^\CCc
\eeq
For the present work, the difference between the indexed families on the right of $\simeq$ or the total categories on the left of $\simeq$ is a matter of convenience. Unless stated otherwise, either standpoint will do and we use the same notations to refer to both sides:
\beq\label{eq:DopUUp}
 \Dfib \CCc \simeq \Do\CCc \simeq \Set^{\CCc^\op}\qquad \quad \Big(\Ofib\CCc\Big)^\op \simeq  \Up\CCc \simeq \left(\Set^\CCc\right)^\op \qquad \quad \Ofib\CCc \simeq \pU\CCc \simeq  \Set^\CCc 
\eeq
 
\paragraph{Arrows that point left are in the opposite category.} $\Up\CCc$-morphisms are written in the form
\beq\label{eq:arrowback}
\left(\rgt X\oot f\rgt Y\right)\  \in\   \Up\CCc\left(\rgt X, \rgt Y\right)\  =\  \pU\CCc\left(\rgt Y, \rgt X\right)
\eeq	
reminding us that $\Up\CCc$-morphisms are $\pU\CCc$-morphisms taken backwards.


\paragraph{A terminological conundrum: \emph{Pre}\/sheaves vs \emph{what}?} The contravariant functors into $\Set$, corresponding to discrete fibrations, are well-known as \emph{presheaves}. The covariant functors into $\Set$, corresponding to discrete opfibrations, are \emph{not}\/ well-known the postsheaves. The former got their name because they arise as an intermediary step in presenting sheaves of sets as the contravariant functors from the lattices of opens into $\Set$, but just those among such functors that satisfy the \emph{sheaf condition}\/  \cite{Godement,GrayW:sheaf-hist,MacLane-Moerdijk}. Following a similar convention, we could just as well call the integers the \emph{preprimes}, because the primes are just those integers that satisfy the condition that they are prime, and we have to construct the integers before we construct the primes. But presheaves have been called that for a long time and sticking with their name saves time. In the present work, however, the adjunction between $\Do\CCc$ and $\Up\CCc$, lifting \eqref{eq:DM} from posets to categories, plays a central role, and we need simple and symmetric names for both. We tried calling them \emph{discrete fibrations}\/ and \emph{discrete opfibrations}. In  \cite{PavlovicD:CALCO15}, we tried calling $\Do\CCc$ presheaves and $\Up\CCc$ postsheaves. Here we lift the standard terminology from group and monoid representations:
\begin{itemize}
\item \  \textbf{left\  \/ \actions}\ \/   $\mapsfrom$\ \  \emph{presheaves / discrete fibrations}, and 
\item \textbf{right \actions}\/ $\mapsfrom$ \emph{postsheaves / discrete opfibrations}, 
\end{itemize}
The added bonus is the reminder that categories are a special case of monoids, as explainned in Appendix~\ref{Appendix:catspan}.

\subsubsection{Loose completions as representations} 
While $\pU\CCc$ is a concrete category, $\CCc$ embeds into $\Up\CCc = \left(\pU\CCc\right)^\op$. The Yoneda embeddings, i.e., the representations
\begin{align}\label{eq:Yoneda}
\mnd : \CCc &\    \tto{\ \ \ \ } \   \Do \CCc & \cmn : \CCc &\   \tto{\ \ \ \ }\    \Up \CCc 
\\
 x & \ \mapsto  \left(\begin{tikzar}[row sep=1.11em,column sep=1.8em]
{\scriptstyle \CCc/x}\arrow{d}[pos = .3]{\scriptscriptstyle \Dom}  \\ 
\scriptstyle\CCc \end{tikzar}\right) &\  x 
& \ \mapsto  \left(\begin{tikzar}[row sep=1.11em,column sep=1.8em]
{\scriptstyle x/\CCc}\arrow{d}[pos = .3]{\scriptscriptstyle \Cod}  \\ 
{\scriptstyle\CCc} \end{tikzar}\right)\notag
\end{align}
are respectively the $\supp$-completion and the $\inff$-completion of $\CCc$. They generalize the $\wedge$-completions and the $\vee$-completions of posets and preorders in the sense of Appendix~\ref{Appendix:lim}. Every functor $K:\CCc\to \KKK$ into a $\supp$-complete category $\KKK$ induces a unique $\supp$-preserving map $K^\ast$ and every functor $L:\CCc\to \LLL$ into a $\inff$-complete category $\LLL$ induces  a unique $\inff$-preserving functor $L_\ast$ making the following diagrams commute
\beq\label{eq:univ-cat}
\begin{tikzar}[row sep = 2ex]
\& \Do\CCc\ar{dd}[description]{\exists ! K^\ast}\\
\CCc\ar{ur}{\mnd}\ar{dr}[description]{\forall K}\\
\& \KKK
\end{tikzar}\qquad \qquad
\begin{tikzar}[row sep = 2ex]
\& \Up\CCc\ar{dd}[description]{\exists! L_\ast}\\
\CCc\ar{ur}{\cmn}\ar{dr}[description]{\forall L}\\ 
\& \LLL
\end{tikzar}
\eeq
This follows from the fact that $\Do\CCc$ is $\supp$-generated by $\CCc$ along $\mnd$ and that $\Up \CCc$ is $\inff$-generated by $\CCc$ along $\cmn$, in the sense that every $\lft A\in \Do\CCc$ and every $\rgt B \in \Up\CCc$ satisfy
\beq\label{eq:inff}
\lft A  = \supp \mnd \lft A  \qquad\qquad\qquad\qquad  \rgt B = \inff \cmn \rgt B 
\eeq 
since $\lft A(x) = \supp \Set(\lft A-, x)$ and $\rgt B(y) = \supp \Set(y, \rgt B-)$. 
%
%
%
%
%
%
%
%
%
%
%
%
%
It is then easy to show that for any finite $\DDd$ and $\EEe$ such that $\inff D$ exists for $D:\DDd\to \CCc$ and $\supp E$ exists for $E:\EEe\to \CCc$ holds
\beq
\mnd(\inff D) = \inff \mnd D \qquad \qquad \qquad \cmn(\supp E) = \supp \cmn  E
\eeq
As before, the $\supp$-completion $\mnd$ does not preserve $\supp$ that may exist in $\CCc$, and the $\inff$-completion $\cmn$ does not preserve the existing $\inff$, but they add fresh $\supp$ and $\inff$, respectively.

\subsection{Tight completions of categories}\label{Sec:cat-tight}

\subsubsection{Lambek's Problem}\label{Sec:Lambek-problem}
In the early days of category theory, as Alexander Grothendieck demonstrated how categorifying 
provided a new foundation for algebraic geometry, Daniel Kan for homotopy theory, and Bill Lawvere for universal algebra, Jim Lambek took up categorifying completions. He presented his work in a series of lectures at the ETH in Z\" urich, and in Vol.~24 of \emph{Springer Lecture Notes in Mathematics} \cite{LambekJ:completions}. In Sec.~2, he lifted the Galois connection $\Kan:\Do\PPp \to \Up\PPp$ over a poset $\PPp$, displayed on \eqref{eq:DM}, to the adjunction $\Kan:\Do\CCc\to \Up \CCc$ over a category $\CCc$, displayed on the next diagram.
\beq\label{eq:DMC}
\begin{tikzar}{}
\&\& \Do \CCc 
\arrow[loop, out = 135, in = 45, looseness = 4]{}[description]{\RLan}
\arrow[bend right = 15]{dddd}[swap]{\Lan} \arrow[phantom]{dddd}{\dashv} 
\arrow[bend right = 15]{rr}[swap]{\ladj U}
\arrow[phantom]{rr}{\scriptstyle\top}
 \&\& \DRL 
\arrow[bend right = 15]{ll}[swap]{\radj U} 
\\
\\
\CCc \arrow{uurr}{\mnd} \arrow{ddrr}[swap]{\cmn} 
\\
\\ 
\&\& \Up \CCc 
\arrow[loop, out = -50, in=-130, looseness = 4]{}[description]{\LRan}
\arrow[bend right = 15,crossing over]{uuuu}[swap]{\Ran}  \arrow[bend right = 15]{rr}[swap]{\ladj V} \
\arrow[phantom]{rr}{\scriptstyle\bot} \&\& \ULR 
\arrow[bend right = 15]{ll}[swap]{\radj V} 
\end{tikzar}
\eeq
The functor $\Lan$ is the left Kan extension of the embedding $\cmn$, while $\Ran$ is the right Kan extension of the embedding $\mnd$. Their constructions echo \eqref{eq:pKan}, replacing $\bigcup$ with $\supp$ and $\bigcap$ with $\inff$. 
\beq\label{eq:cKan}
\begin{tikzar}[column sep=4em]
\lft A = \displaystyle \supp \mnd \lft A  \arrow[thin,mapsto]{dd}
\& \Do \CCc\arrow[bend right = 15]{dd}[swap]{\Lan} \arrow[phantom]{dd}{\dashv} \& \displaystyle\Ran \rgt B = \inff \mnd \rgt B  
\\
\\ 
\displaystyle\Lan\lft A = \inff\cmn \lft A \&\left(\pU\CCc\right)^\op 
\arrow[bend right = 15]{uu}[swap]{\Ran}  \& 
\displaystyle \rgt B = \supp \cmn \rgt B
\arrow[thin,mapsto]{uu}
\end{tikzar}
\eeq
where we write $\left(\pU\CCc\right)^\op$ instead of $\Up\CCc$ to emphasize that $\supp$ and $\inff$ at the bottom are in the category $\pU\CCc$ of {\postsheaf}s over $\CCc$, and not in $\Up\CCc=\left(\pU{\CCc}\right)^\op$, seen in \eqref{eq:DopUUp}. Viewing 
\[\lft \AAa\tto{\lft A} \CCc\tto{\cmn}\left(\pU\CCc\right)^\op\qquad \mbox{ and }\qquad \rgt \BBb\tto{\rgt B} \CCc\tto{\mnd}\Do\CCc
\]
as diagrams and applying the Yoneda Lemma and the fact that the representable functors preserve limits gives
\beq\label{eq:cones-cocones}
\begin{split}
\left(\Lan\lft A\right)_c \cong \pU\CCc\left(\cmn c, \inff \cmn \lft A\right) \cong \inff\pU\CCc\left(\cmn c,\cmn\lft A\right) \cong \inff\CCc\left(\lft A, c\right) = \Coc_\CCc\left(\lft A, c\right)\\
\left(\Ran\rgt B\right)_c \cong \Do\CCc\left(\mnd c, \inff \mnd \rgt B\right) \cong \inff \Do\CCc\left(\mnd c, \mnd\rgt B\right) \cong \inff \CCc\left(c,\rgt B\right) \cong \Con_\CCc\left(c, \rgt B\right)
\end{split}
\eeq
where $\Coc_\CCc\left(\lft A, c\right) = \Do\CCc\left(\lft A, \mnd c\right)$ is the set of cones and $\Con_\CCc\left(c,\rgt A\right) = \pU\CCc\left(\rgt A, \cmn c\right)$ is the set of cocones. To re-state \eqref{eq:cones-cocones} in words,  
\begin{itemize}
\item $\Lan$ takes $\lft A$ to the {\postsheaf} $\Lan\lft A: \Coc(\lft A) \tto{\Cod} \CCc$ of cocones over $\lft A$, whereas 
\item $\Ran$ takes $\rgt B$ to the {\presheaf} $\Ran\rgt B: \Con(\rgt B) \tto{\Dom} \CCc$ of cones under $\rgt B$, \end{itemize}
where $\Coc(\lft A) = \coprod_{c\in \CCc} \Coc(\lft A, c)$ and $\Con(\rgt A) = \coprod_{c\in \CCc} \Con(c,\rgt A)$. Finally, the unit and the counit of the adjunction $\Kan$ from \eqref{eq:cKan} are
\begin{align}
\eta \colon \lft A & \longrightarrow \Ran\Lan\lft A & \varepsilon \colon \Lan \Ran \rgt A  \longleftarrow &\  \rgt A\\
x & \mapsto \lambda \delta.\ \delta_x & \lambda \varrho.\ \varrho_y \mapsfrom &\  y\notag
\end{align}
which means that for $x\in \lft \AAa$ and $y\in \rgt \AAa$ 
\begin{itemize}
\item  $\eta(x): \lft A x \to \Lan\lft A$ is the cone $\left\{{\delta_x}: \lft Ax \to\  |\delta|\ \Big|\ \delta\in \Coc\big(\lft A,|\delta|\big) \right\}$, whereas
\item $\varepsilon(y): \Ran\rgt A  \to \rgt A y$ is the cocone $\left\{\varrho_y:|\varrho| \to \rgt Ay\ \Big|\ \varrho\in \Con\big(|\varrho|, \rgt A\big) \right\}$.
\end{itemize}
\paragraph{Summary.} Going back to \eqref{eq:pKan}, we see that 
\begin{itemize}
\item the {\postsheaf} of cocones generalizes the upper set of upper bounds, whereas
\item the {\presheaf} of cones generalizes the lower set of lower bounds.
\end{itemize}
The adjunction $\Kan$ of the extensions from   \eqref{eq:cKan} was displayed in Lambek's \cite[Sec.~2]{LambekJ:completions}, studied in Isbell's \cite{IsbellJ:structure}, and analyzed ever since, as the \emph{Isbell envelope, conjugacy, duality}  \cite{Avery-Leinster:isbell,DiLiberti:isbell,GarnerR:isbell}. To contribute to the variety, we usually call it the \emph{Isbell adjunction}, and follow \cite{GarnerR:isbell} in calling  the  \emph{Isbell monad}\/ the construction of the isomorphic comma categories induced by an extension adjunction. While the Isbell monad is thus a monad on the category of categories, every image of a category $\CCc$ along the Isbell monad, i.e., every instance of an Isbell adjunction \eqref{eq:cKan} induces a monad $\RLan = \Ran\Lan$ on $\Do\CCc$ and a comonad $\LRan = \Lan\Ran$ on $\Up\CCc = \left(\pU\CCc\right)^\op$, which is a monad on $\pU\CCc$.  At the time when Isbell and Lambek studied completions, monads were nascent, and neither author considered the categories of algebras $\DRL$ and $\ULR$ as the liftings of the lattices of closed sets $\Do\PPp^\RLan$ and $\Up\PPp^\LRan$. Lambek thus succeeded in lifting the main posetal completion constructions to categories, but closed his series of lectures and lecture notes by the statement quoted as motto at the beginning of this paper, announcing the open problem of categorifying the Dedekind-MacNeille completion.

\subsubsection{Isbell's dissolution of Lambek's Problem}\label{Sec:dissolution}
A couple of years later, John Isbell closed Lambek's Problem with the theorem also displayed at the beginning of this paper \cite[Thm.~3.1]{IsbellJ:no-lambek}. Isbell showed that no category with a $\supp$- and $\inff$-dense embedding of the cyclic group $\ZZz_4$ can have a particular equalizer that Isbell constructed. It follows that $\ZZz_4$ cannot be $\inff$- and $\supp$-densely embedded into an $\inff$-complete category, and that Lambek's Problem is generally not solvable for standard categorical limits and colimits. Isbell's equalizer is a compelling no-go result, but it does not explain why the road is closed. To see what goes wrong with generalizing the Dedekind-MacNeille schema \eqref{eq:DM} from posets to categories, consider the schema \eqref{eq:DMC}, this time instantiated to $\CCc = \ZZz_4$, the cyclic group of integer addition modulo 4, where a group is viewed as a category with a single object\footnote{Any group would do here. We use $\ZZz_4$ in reference to Isbell, who used this particular group to construct a particular equalizer. It is the smallest group with a nontrivial subgroup.}. 
%
%
The category $\Do\ZZz_4$ of left actions $\lft X=\left(\ZZz_4\times X\to X\right)$ is a $\supp$-completion of $\ZZz_4$, whereas $\Up\ZZz_4$ of right actions $\rgt Y = \left(Y \times \ZZz_4\to Y\right)$ is a $\inff$-completion. They are special cases of \eqref{eq:doupcat}. For a group, the two completions are, however, dually equivalent: \bea\label{eq:eqv}
\Do\ZZz_4 & \simeq & \left(\Up\ZZz_4\right)^\op
\eea
This is because the isomorphism $\ZZz_4\cong \ZZz^\op_4$, realized by the group inverse, induces the equivalence
\[
\Do \ZZz_4\  \ \simeq\  \ \Set^{\ZZz_4^\op} \ \ \simeq\ \  \Set^{\ZZz_4}\ \ \simeq \ \  \pU{\ZZz_4}
\]
while \eqref{eq:DopUUp} instantiates to
\[
\pU{\ZZz_4} \ \ \simeq \ \ \left(\Up \ZZz_4\right)^\op
\]
The Cayley group representations $\mnd:\ZZz_4\to \Do\ZZz_4$ and $\cmn:\ZZz_4\to\Up \ZZz_4$  are a special case of the Yoneda embeddings \eqref{eq:Yoneda}. 
Their Kan extensions form the adjunction $\Kan$ which maps the free group actions to the powers of their orbit sets, and kills the actions that are not free:
\begin{align}
\Lan \lft X & = \begin{cases} \ZZz_4^{X_0}  & \mbox{ if } X= \ZZz_4\times X_0\\
\emptyset & \mbox{ otherwise}
\end{cases}
&
\Ran \rgt Y & = \begin{cases}
\ZZz_4^{Y_0}  & \mbox{ if } Y= Y_0\times \ZZz_4\\
\emptyset & \mbox{ otherwise}
\end{cases}
\end{align}
Group actions $\lft X$ and $\rgt Y$ are free when the underlying sets are in the form $X = \ZZz_4\times X_0$ and $Y = Y_0\times \ZZz_4$ for some $X_0$ and $Y_0$ which the  $\ZZz_4$-action leaves fixed, acting only on the $\ZZz_4$-component:
\begin{align*}
\lft X \ = \ \Big(\ZZz_4\times \ZZz_4\times X  & \tto{(+)\times X_0} \ZZz_4\times X_0\Big)
& \qquad\qquad \rgt Y\  =\  \Big(Y_0\times \ZZz_4\times \ZZz_4 & \tto{Y_0\times(+)} Y_0\times \ZZz_4\Big)
\\
<m,n,x> &\ \  \longmapsto\ \  <m+n, x> & <y,m,n> &\ \  \longmapsto\ \  <y,m+n>
\end{align*}
The elements of $X_0$ and $Y_0$ remain fixed under the actions, and their elements thus correspond to the orbits. The $\Lan$-image of a free $\lft X$ and the $\Ran$-image of a free $\rgt Y$ are the pointwise actions 
\begin{align*}
\ZZz_4^{X_0} \times \ZZz_4 &\tto{\Lan\lft X} \ZZz_4^{X_0} 
& \ZZz_4 \times \ZZz_4^{Y_0} &\tto{\Ran\rgt Y} \ZZz_4^{Y_0}\\
<\xi, i> & \longmapsto \lambda x.\ \xi(x)+i &  <i,\upsilon > & \longmapsto \lambda y.\ i+\upsilon(y)
\end{align*} 
On the other hand, the adjunction $\Kan:\Up\ZZz_4\to \Do\ZZz_4$ induces on $\Do\ZZz_4$ the monad $\RLan=\Ran\Lan$ and on $\Up \ZZz_4$ the comonad $\LRan = \Lan\Ran$. They play the role of closure operators from \eqref{eq:closureops}. Their algebras and coalgebras play the roles of the closed sets from \eqref{eq:closed}. But while the two lattices of closed sets turn out to be isomorphic, the categories of algebras and of coalgebras are quite different. To spell out $\RLan\lft X = \Ran\Lan \lft X$, note that the pointwise action $\Lan \lft X = \left(\ZZz_4^{X_0}\times \ZZz_4 \to  \ZZz_4^{X_0}\right)$ is free, and that its set of orbits is the quotient $\ZZz_4^{X_0}/\ZZz_4$. Similarly for $\Ran \rgt Y = \left(\ZZz_4\times \ZZz_4^{Y_0} \to  \ZZz_4^{Y_0}\right)$, where the set oif orbits is $\ZZz_4^{Y_0} /\ZZz_4$. Picking out an arbitrary orbit $\bullet \in X_0$ makes $X_0$ into a disjoint union of $X_0^{\setminus \bullet} = X_0\setminus \{\bullet\}$ and $\{\bullet\}$, and induces the decomposition 
\[\ZZz_4^{X_0}\ \cong\  \left(\ZZz_4^{X_0} /\ZZz_4\right)\times \ZZz_4\  \cong\  \ZZz_4^{X_0^{\setminus \bullet}} \times \ZZz_4\]
Picking out an arbitrary orbit from $\Ran \rgt Y$, denoted by abuse of notation also $\bullet\in Y_0$, induces the decomposition 
\[ \ZZz_4^{Y_0}\ \cong\ \ZZz_4\times \left(\ZZz_4^{Y_0}/\ZZz_4\right)\ \cong\  \ZZz_4\times \ZZz_4^{Y_0^{\setminus \bullet}}\]
The actions on the cocones over the free $\lft X$ and on the cones under $\rgt Y$ and can now be written in the form
\beq 
\Lan \lft X=\left(\ZZz_4^{X_0^{\setminus \bullet}}\times \ZZz_4\times \ZZz_4 \tto{\id\times (+)} \ZZz_4^{X_0^{\setminus \bullet}}\times \ZZz_4\right) \qquad\quad
\Ran \rgt Y= \left(\ZZz_4\times \ZZz_4 \times \ZZz_4^{Y_0^{\setminus \bullet}} \tto{(+)\times \id} \ZZz_4 \times \ZZz_4^{Y_0^{\setminus \bullet}} \right)\eeq
The monad $\RLan = \Ran\Lan :\Do\ZZz_4 \to \Do \ZZz_4$ and the comonad $\LRan = \Lan\Ran :\Up\ZZz_4 \to \Up \ZZz_4$ are thus
\begin{align}
\RLan \lft X & = \begin{cases} \ZZz_4^{\ZZz_4^{X_0^{\setminus \bullet}}}  & \mbox{ if } X= \ZZz_4\times X_0\\
1 & \mbox{ otherwise}
\end{cases}
&
\LRan \rgt Y & = \begin{cases}
\ZZz_4^{\ZZz_4^{Y_0^{\setminus \bullet}}}  & \mbox{ if } Y= Y_0\times \ZZz_4\\
1 & \mbox{ otherwise}
\end{cases}
\end{align}
The equivalence $\Do\ZZz_4 \simeq \left(\Up\ZZz_4\right)^\op$ from \eqref{eq:eqv} therefore lifts to
\bea\label{eq:eqv-mnd}
\Do\ZZz_4^\RLan & \simeq & \left(\Up\ZZz_4^\LRan\right)^\op
\eea
The categories $\Do\ZZz_4^\RLan$ and $\Up\ZZz_4^\LRan$ could thus only be dual if they were self-dual. But one is built from quotients of coproducts of $\ZZz_4$, and the other from subobjects of products. While infima and suprema are symmetric and dual, products and coproducts are not, e.g. because for sets or numbers the products distribute over the coproducts but not vice versa.

\section{Approach: The nucleus of an adjunction}\label{Sec:nuc}

The crucial insight that opens up a path towards completing a suitable categorical version of the Dedekind-MacNeille schema \eqref{eq:DM} is that swapping the categories of algebras $\DRL$ and and coalgebras $\ULR$ in \eqref{eq:DMC}, adding the comparison functors, and composing them with the forgetful functors turns out to yield an adjunction. This adjunction is the \emph{nucleus}\/ \cite{PavlovicD:nucleus} of the Isbell adjunction, which we spell out in this section. The Dedekind-MacNeille isomorphism from \eqref{eq:DM} boils down to this nucleus when posets are generalized to categories. This nucleus is generally not an equivalence of categories, unless the category from which it is derived is a poset, or preorder. Nevertheless, the categories $\DRL$ and $\ULR$ always completely determine one another, and yield the tight completion of $\CCc$. In this section we spell out the nucleus. Categories of algebras and coalgebras for monads and comonads with their resolutions, going back to Eilenberg and Moore \cite{Eilenberg-Moore}, are described in most textbooks in category theory, for example \cite{BarrM:ttt,Lambek-Scott:book,MacLaneS:CWM}.

\subsection{The algebra-coalgebra adjunction}\label{Sec:Cnuc}
\begin{proposition}\label{Prop:Cnuc}
For any category $\CCc$, the Isbell adjunction $\Kan:\Up\CCc\to \Do\CCc$ induces the following schema:
\beq\label{eq:Cnuc}
\begin{tikzar}[row sep=3.5cm,column sep=4.5cm]
\Do\CCc \arrow[phantom]{d}[description]{\dashv}  
\arrow[loop, out = 135, in = 45, looseness = 4,thin]{}[description]{\mbox{\large$\RLan$}} 
\arrow[bend right = 13]{d}[swap]{\mbox{\large$\Lan$}} 
\arrow[thin,dashed]{r}[description]{\mbox{\large$\compson^{0}$}} 
\& 
\ULR \arrow[bend right = 13,pos=0.65,thin]{dl}[description]{\mbox{\large$\ladj V$}} 
\arrow[phantom]{d}[description]{{\dashv}}   
\arrow[bend right = 13]{d}[swap]{\mbox{\large$\nRan$}}
\\
\Up\CCc  
\arrow[loop, out = -45, in=-135, looseness = 6,thin]{}[description]{\mbox{\large$\LRan$}} 
\arrow[thin,dashed]{r}[description]{\mbox{\large$\compson_{1}$}}  
\arrow[bend right = 13]{u}[swap]{\mbox{\large$\Ran$}} 
\& 
\DRL 
\arrow[bend right = 13]{u}[swap]{\mbox{\large$\nLan$}}  \arrow[bend right = 13,crossing over,thin]{ul}[description,pos=.75]{\mbox{\large$\radj U$}} 
\end{tikzar}\eeq
The following well-known structures are displayed:
\begin{itemize}
\item the monad $\RLan=\Ran\Lan :\Do\CCc\to\Do\CCc$ with the (Eilenberg-Moore) category of algebras $\DRL$,
\item the comonad $\LRan = \Lan\Ran:\Up\CCc\to \Up \CCc$ with the (Eilenberg-Moore) category of coalgebras $\ULR$,
\item the free-forgetful resolution $\adj U:\DRL \to \Do\CCc$ of $\RLan$,
\item the forgetful-cofree resolution $\adj V:\Up\CCc \to \ULR$ of $\LRan$,
\item the comparison functor $\compson_1:\Up\CCc\to \DRL$ from $\Kan$ to $\adj U$ as resolutions of the monad $\Ran$ on $\Do\CCc$, which means that it  satisfies  
\[ \Ran = \radj U \circ \compson_1\qquad\qquad\qquad\qquad \compson_1\circ \Lan= \ladj U\]
\item the \emph{comparison functor}\/ $\compson^0:\Do\CCc \to \ULR$ from $\Kan$ to $\adj V$ as resolutions of the comonad $\LRan$ on $\Up\CCc$, which means that it  satisfies 
\[ \Lan = \ladj V \circ \compson^0 \qquad\qquad\qquad\qquad\Ran \circ \compson^0 = \radj V\]
\end{itemize}
We furthermore define the following functors between the algebras in $\DRL$ and the coalgebras in $\ULR$
\[ \nLan = \compson_0\circ\radj U\qquad\qquad\qquad\qquad
\nRan = \compson_1\circ \ladj V\]
which means that they act as follows
\beq\label{eq:funcdefn}
\begin{tikzar}[row sep = 1.8em,column sep=1.5em]
\Do\CCc \arrow[thin]{d}[pos=.35]{\compson^0}  \& \lft A\arrow[mapsto,thin]{d}   \&\&\&\& \&\&
\ULR \arrow[thin,bend right = 50]{dd}[swap]{\ladj V}
 \arrow[thick]{d}[pos=.35]{\nRan} \& \big({\rgt B} \oot{\ \beta\ } \Lan \Ran {\rgt B}\big)  \arrow[mapsto]{d}  \\
\ULR \& 
\big(\Lan {\lft A} \oot{\Lan \eta} \Lan \Ran \Lan {\lft A}\big)
\&\&\&\&\&\&
\DRL \&
\big(\Ran \Lan\Ran {\rgt B} \tto{\Ran \varepsilon} \Ran {\rgt B}\big) \\
\DRL \arrow[thin,bend left = 50]{uu}{\radj U}
\arrow[thick]{u}[swap,pos=.35]{\nLan} \& \big(\Ran \Lan {\lft A} \tto{\ \alpha\ } {\lft A}\big) \arrow[mapsto]{u}
\&\&\&\&\&\&
\Up\CCc 
\arrow[thin]{u}[swap,pos=.4]{\compson_1}  \& \arrow[mapsto,thin]{u} \rgt B 
\end{tikzar}
\eeq
The claim is that they form the\/ \emph{algebra-coalgebra adjunction} $\nRan\dashv \nLan :\DRL \to \ULR$.
\end{proposition}

\bprf{ sketch}
For any $\RLan \lft A \tto\alpha \lft A$ in $\DRL$ and $\rgt B\oot\beta \LRan \rgt B$ in $\ULR$, the adjunction correspondence
\bea\label{eq:nbij}
\DRL(\nRan \beta, \alpha) & \cong & \ULR(\beta, \nLan \alpha)
\eea 
follows from the fact that each of the following squares commutes if and only if the other commutes
\bea\label{eq:bigadj}
\begin{tikzar}[row sep=1.8em,column sep=1.3em]
\Ran\Lan \Ran {\rgt B} \ar{rr}[description]{\Ran\Lan \lft \intv} \ar{dd}[swap]{\Ran \varepsilon} \& \& \Ran \Lan {\lft A} \ar{dd}{\alpha}\\ \& \hspace{1em}  \\
\Ran {\rgt B} \ar{rr}[description]{\lft \intv} \ar{uurr}[description]{\Ran\rgt \intv} \&\& {\lft A}
\end{tikzar}
\hspace{2em}& \iff &\hspace{2em}
\begin{tikzar}[row sep=1.8em,column sep=1.3em]
{\rgt B} \ar[leftarrow]{rr}[description]{\rgt \intv} \ar[leftarrow]{dd}[swap]{\beta} \&\& \Lan {\lft A} \ar[leftarrow]{dd}{\Lan \eta}\\ \& \hspace{1em} \\
\Lan \Ran {\rgt B} \ar[leftarrow]{uurr}[description]{\Lan \lft \intv} \ar[leftarrow]{rr}[description]{\Lan \Ran \rgt \intv}  \&\&  \Lan \Ran \Lan {\lft A} 
\end{tikzar}
\eea
where we use $\nRan \beta = \Ran\varepsilon$ and $\nLan \alpha = \Lan \eta$ from \eqref{eq:funcdefn}. The natural bijection in \eqref{eq:nbij} is realized using the diagonals:
\[
\lft \intv \ \  \mapsto \ \  \rgt \intv = \beta \circ \Lan \lft \intv\qquad\qquad\qquad\qquad
\rgt \intv \ \  \mapsto \ \  \lft \intv = \alpha \circ \Ran \rgt \intv
\]
See \cite[Sec.~7.1]{PavlovicD:nucleus} for a detailed verification that such maps induce the claimed adjunction between algebras and coalgebras induced by any given adjunction, here instantiated to $\Kan$.
\epr

\subsection{The algebra-coalgebra adjunction is nuclear}
\begin{proposition}\label{Prop:Cnuc-monadic}
The algebra-coalgebra adjunction $\nRan\dashv \nLan :\DRL \to \ULR$ established in Prop.~\ref{Prop:Cnuc} is nuclear \cite{PavlovicD:nucleus}:  the functor $\nLan$ is monadic and $\nRan$ is comonadic. This means that lifting the schema displayed in \eqref{eq:Cnuc} from the Isbell adjunction $\Kan$ 
to the algebra-coalgebra adjunction $\nRan\dashv \nLan$ leads to the following situation: 
\beq\label{eq:Cnuc-monadic}
\begin{tikzar}[row sep=3.5cm,column sep=4.5cm]
\ULR \arrow[phantom]{d}[description]{\dashv}  
\arrow[loop, out = 135, in = 45, looseness = 4,thin]{}[description]{\mbox{\large$\RLAN$}} 
\arrow[bend right = 13]{d}[swap]{\mbox{\large$\nRan$}} 
\arrow[thin,dashed]{r}[description]{\mbox{\LARGE $\sim$}} 
\& 
\left(\DRL\right)^{\RLAN} \arrow[bend right = 13,pos=0.8,thin]{dl}
\arrow[phantom]{d}[description]{{\dashv}}   
\arrow[bend right = 13]{d}[swap]{\mbox{\large$\nnLan$}}  
\\
\DRL  
\arrow[loop, out = -45, in=-135, looseness = 6,thin]{}[description]{\mbox{\large$\LRAN$}} 
\arrow[thin,dashed]{r}[description]{\mbox{\LARGE $\sim$}}  
\arrow[bend right = 13]{u}[swap]{\mbox{\large$\nLan$}} 
\& 
\left(\ULR\right)^{\LRAN} 
\arrow[bend right = 13]{u}[swap]{\mbox{\large$\nnRan$}}  \arrow[bend right = 13,crossing over,thin]{ul}
\end{tikzar}\eeq
where the comparison functors are equivalences. 
In other words,
\begin{itemize}
\item $\ULR$ and the monad $\LRAN = \nLan\nRan$  on it determines $\DRL\simeq \left(\ULR\right)^{\LRAN}$, and conversely,
\item $\DRL$ and the monad $\RLAN = \nRan\nLan$  on it determines $\ULR\simeq \left(\DRL\right)^{\RLAN}$.
\end{itemize}
\end{proposition}

A detailed proof of the claim of Prop.~\ref{Prop:Cnuc-monadic} and a discussion of its consequences can be found in \cite[Sec.~7.2]{PavlovicD:nucleus}.

\section{Construction: From gaps and intervals to cuts}\label{Sec:cut}
\subsection{Gaps and intervals in posets}\label{Sec:gapint}

The algebra-coalgebra adjunction $\nKan$, spelled out in \eqref{eq:Cnuc} for categories, is obviously also  supported by posets as categories with at most one morphism per hom-set. The isomorphism
$$\Up\PPp^\LRan\ \  \cong\ \  \UD \PPp\ \ \cong\ \  \Do\PPp^\RLan$$
in the Dedekind-MacNeille schema \eqref{eq:DM} can now be explained as an intersection of the adjunctions $\Kan$ and $\nKan$. The schema in \eqref{eq:DM} is first refined to  
\beq\label{eq:DM-refin}
\begin{tikzar}[row sep=3.3ex,column sep=1.5em]
\&\& \Do \PPp 
\arrow[loop, out = 135, in = 45, looseness = 4]{}[description]{\RLan} 
\arrow[bend right = 15]{dddd}[swap]{\Lan} \arrow[phantom]{dddd}{\dashv} 
\arrow[two heads,bend right = 10,thin]{rrrr}
\arrow[phantom]{rrrr}{\scriptstyle\top}
 \&\&\&\&  \Up\PPp^{\LRan}
\arrow[hookrightarrow,bend right = 10,thin]{llll} 
\arrow[bend right = 15]{dddd}[swap]{\nRan} \arrow[phantom]{dddd}{\dashv}
\\
\\
\PPp \arrow{uurr}{\mnd} \arrow{ddrr}[swap]{\cmn} 
\\
\\ 
\&\& \Up \PPp 
\arrow[loop, out = -50, in=-130, looseness = 4]{}[description]{\LRan} 
\arrow[bend right = 15]{uuuu}[swap]{\Ran}  
\arrow[two heads,bend right = 10,thin]{rrrr} 
\arrow[phantom]{rrrr}{\scriptstyle\bot} 
\&\&\&\& 
\Do\PPp^{\RLan} 
\arrow[hookrightarrow,bend right = 10,thin]{llll} 
\arrow[bend right = 15]{uuuu}[swap]{\nLan}  
\end{tikzar}
\eeq
The adjunctions $\Kan$ and $\nKan$ mean that for $\lft X\in \Do\PPp$ and $\rgt Y\in \Up\PPp$ holds
\bea
\lft X\subseteq \Ran\rgt Y\ \  \iff&\ \ \  \lft X\leq\rgt Y\ \ \ &\iff\ \   \Lan\lft X\supseteq \rgt Y \label{eq:LR}\\
\Ran\rgt Y \subseteq \lft X\ \ \ \  \iff&\  \ \Ran\rgt Y\leq \Lan \lft X\ \  &\iff\ \ \ \ \rgt Y \supseteq \Lan\lft X \label{eq:nLR}
\eea
where we write $\lft X\leq \rgt Y$ to abbreviate $\forall x\in \lft X\ \forall y\in \rgt Y.\ x\leq y$, and in \eqref{eq:nLR} use $\Ran\Lan \lft X \subseteq \lft X$ for $\lft X\in \Do\PPp^\RLan$ and $\rgt Y\supseteq \Lan\Ran \rgt Y$ for $\rgt Y\in \Up\PPp^\LRan$. For a more intuitive picture, writing
\bea
\commm{\adj F} & = & \left\{<x,y>\in \XXx\times \YYy\ |\ \ladj F x\leq y\right\} \ \ =\ \ \left\{<x,y>\in \XXx\times \YYy\ |\  x\leq \radj F y\right\} 
\eea
for the distributor induced by any adjunction $\adj F:\BBb\to \AAa$, we can factorize the adjunctions in \eqref{eq:DM-refin} as follows:
\begin{gather}\label{diag:DM-fac}
\begin{tikzar}[column sep = 2.5cm,row sep = 1.5cm]
\Do \PPp \arrow[bend right=35,thin]{dd}
\arrow[bend right=15]{d}
\arrow[phantom]{d}[description]{\dashv} 
\arrow[thin,two heads]{r} 
\& 
\Up\PPp^\LRan 
\arrow[bend right=35,thin]{dd}
\arrow[bend right=15]{d}
\arrow[phantom]{d}[description]{\dashv} 
\\
\commm{\Kan}
\arrow[bend right=15]{u}
\arrow[bend right=15]{d}
\& 
\commm{\nKan}
\arrow[bend right=15]{u}
\arrow[bend right=15]{d}
\\
\Up \PPp \arrow[bend right=35,thin]{uu}
\arrow[phantom]{u}[description]{\dashv}
\arrow[bend right=15]{u}
\arrow[thin,two heads]{r} 
\&
\Do\PPp^\RLan 
\arrow[bend right=35,thin]{uu}
\arrow[phantom]{u}[description]{\dashv}
\arrow[bend right=15]{u}
\end{tikzar}
\end{gather}
By (\ref{eq:LR}--\ref{eq:nLR}), the distributors can in this case be written:
\bea\label{eq:commm-Kan}
\commm{\Kan} & = & \left\{ \left<\lft X, \rgt Y \right>\in \Do\PPp\times \Up\PPp\ \ \ \ |\ \ \  \forall x\in \lft X\ \forall y\in \rgt Y.\ \ x\leq y\right\}\\
\commm{\nKan} & = & \left\{ \left<\rgt Y,\lft X \right>\in \Up\PPp^\LRan\times \Do\PPp^\RLan \ |\ \forall x\leq \rgt Y \ \forall y\geq \lft X.\ x\leq y\right\}\label{eq:commm-nKan}
\eea
\paragraph{Tightening gaps and intervals.} The idea is  that 
\begin{itemize}
\item $\commm{\Kan}$ is the lattice of \emph{\textbf{gaps}} $\lft X\Big)\cdots \Big(\rgt Y$, which gives $\bigvee \lft X \leq \bigwedge \rgt Y$, whereas
\item $\commm{\nKan}$ is the lattice of \emph{\textbf{intervals}}   $\mnd \rgt Y\Big)\Big(\rgt Y\cdots \lft X\Big)\Big( \cmn \lft X$, which gives $\bigwedge \rgt Y \leq \bigvee \lft X $.
\end{itemize}
The intersection of the two, constructed as the pullback  
\beq\label{diag:intersection-P}
\begin{tikzar}[column sep = 0.5cm,row sep = 1cm]
\&\UD \PPp \ar{dr} \ar{dl} \ar[phantom]{dd}[description,pos=0.1]{\pbdown}
\\
\commm{\Kan} \ar[hook]{dr} \&\& \commm{\nKan} \ar[tail]{dl}[description]{\sigma}
\\
\& \Do\PPp\times \Up\PPp
\end{tikzar}
\eeq 
with $\sigma\left(\rgt Y, \lft X\right) = \left<\lft X, \rgt Y\right>$, yields the Dedekind-MacNeille completion 
\begin{itemize}
\item $\UD\PPp$ as the lattice of \emph{\textbf{cuts}}\ \  $\Ran \rgt Y= \lft X\Big)\Big( \rgt Y = \Lan \lft X $, which gives $\bigvee \lft X=\bigwedge \rgt Y$.
\end{itemize}       
This view of the tight completion of posets will guide us to the tight completion of categories.

\subsection{Gaps and intervals in categories}\label{Sec:bicnuc}

\paragraph{A general lemma.} The factorizations of  adjunctions in \eqref{diag:DM-fac} are instances of the general approach to adjunctions through comma categories.

\begin{lemma}\label{Lemma:comma}
Every adjunction $\adj F$ factors through a reflection $\Cod\dashv\varepsilon$ and a coreflection $\eta\dashv \Dom$
\beq\label{eq:adjun-fac}
\begin{tikzar}[column sep = large]
\AAa \arrow[bend right = 30,shift right = 2]{rrrr}[swap]{\ladj F}  \arrow[bend right = 15,shift right = 1,tail]{rr}[description]{\eta} \arrow[phantom]{rr}[description]{\top}  \&\&\commm{\adj F} \arrow[bend right = 15,shift right = 1,two heads]{rr}[description]{\Cod} \arrow[bend right = 15,shift right = 1,two heads]{ll}[description]{\Dom}
\arrow[phantom]{rr}[description]{\top}
\&\& 
\BBb \arrow[bend right = 30,shift right = 2]{llll}[swap]{\radj F}  \arrow[bend right = 15,shift right = 1,tail]{ll}[description]{\varepsilon}
\end{tikzar}
\eeq
where
\bea\label{eq:commm-adj}
\Big\lvert \commm{\adj F} \Big\rvert & = &  \coprod_{\substack{a\in \AAa\\ b \in \BBb}}\  
\left\{\ \left<
\begin{matrix}u\\ \vspace{-.1ex}\\ v
\end{matrix}\right>\  \in\ \  
\begin{matrix}\AAa\left(a, \radj F b\right)\\ \times\\ \BBb\left(\ladj F a, b \right)\end{matrix}\ \ \ \Bigg|\ \ \  
\begin{tikzar}[row sep = 1ex]
a\ar{rr}{u} \ar{dd}[description]{\eta}\&\&  \radj F b\\
\&\&\ladj F\radj F b\ar{dd}[description]{\varepsilon}
\\
\radj F\ladj F a \ar{uurr}[description,pos=0.4]{\radj F v}
\\
\ladj F a \ar{uurr}[description,pos=0.45]{\ladj F u} \ar{rr}[pos=0.4,swap]{v}\&\& b
\end{tikzar}
\ \ 
\right\} 
\\
 \commm{\adj F}(x,y) & = &  \left\{
\left<\begin{matrix} f \\ \vspace{-.005ex}\\ g \end{matrix}\right>\  \in\ \  \begin{matrix}\AAa\left(a_x, a_y\right)\\ \times\\ \BBb\left(b_x, b_y\right)\end{matrix}  \ \Bigg|\ \ \ 
\begin{tikzar}[row sep = 3ex,column sep = 3ex]
a_x\ar{rr}{f} \ar{dd}[description]{u_x}\&\& a_y   \ar{dd}[description]{u_y}\\
\\
\radj F b_x \ar{rr}[swap]{\radj F  g}\&\&\radj F b_y
\end{tikzar}
\ \ \ \wedge \ \ \ 
\begin{tikzar}[row sep = 3ex,column sep = 3ex]
\ladj F a_x\ar{rr}{\ladj F f} \ar{dd}[description]{v_x }\&\&  \ladj F b_x \ar{dd}[description]{v_y}\\
\\
b_x  \ar{rr}[swap]{g}\&\&b_y
\end{tikzar}
\ \ \right\}\notag
\eea
The functors in \eqref{eq:adjun-fac} are defined by
\begin{align*}
\eta(a) = & \begin{pmatrix}
a & a\tto\eta \radj F\ladj F a\\
\ladj F a & \ladj F a\tto\id \ladj Fa\end{pmatrix}
 & 
\Dom\begin{pmatrix}
a & a\tto u \radj Fb\\
b & \radj F a\tto v  b \end{pmatrix} & = a\\
\varepsilon(b) = & \begin{pmatrix}
\radj F b & \radj F b \tto\id \radj Fb \\
b & \ladj F \radj F b \tto\varepsilon b\end{pmatrix}
&  
\Cod\begin{pmatrix}
a & a\tto u \radj Fb\\
b & \radj F a\tto v  b \end{pmatrix} &= b
\end{align*}
\end{lemma}
The projections $\commm{\adj F} \to \comm{\AAa}{\radj F}$ and $\commm{\adj F}\to \comm{\ladj F}{\BBb}$ are easily seen to induce the isomorphisms 
\beq \comm{\AAa}{\radj F}\ \ \cong\ \ \commm{\adj F}\ \ \cong\ \ \comm{\ladj F}{\BBb}\eeq
The category $\commm{\adj F}$ is thus highly redundant, as it  glues together the isomorphic categories $\comm{\AAa}{\radj F}$ and $\comm{\ladj F}{\BBb}$. This redundancy will simplify the coming constructions just as the poset $\UD\PPp$ of cuts in \eqref{eq:DM} simplified the tight completion of a poset.

\paragraph{Factoring lose and tight extensions.} The category of gaps and the category of intervals arise, just like the eponymous lattices \eqref{diag:DM-fac}, by factoring adjunctions.
\begin{equation}\label{diag:Kan}%
\begin{tikzar}[column sep = 6cm,row sep = 2cm]
\Do \CCc \arrow[bend right=45]{dd}[swap]{\mbox{\Large$\Lan$}}
\arrow[bend right=10]{d}[swap,pos=.7]{\eta}
\arrow[phantom]{d}[pos=0.66,description]{\dashv} 
\arrow[thin,two heads]{r} 
\& 
\Up\CCc^\LRan 
\arrow[bend right=45]{dd}[swap]{\mbox{\Large$\nRan$}}
\arrow[bend right=10]{d}[swap,pos=0.7]{\eta}
\arrow[phantom]{d}[pos=0.66,description]{\dashv} 
\\
\begin{array}{c}
 \Do\CCc/\Ran
\\
 \cong
\\
 \commm{\Kan}
\\
 \cong
\\
 \Lan{} / \Up \CCc
\end{array}
\arrow[bend right=10]{u}[swap,pos=0.3]{\Dom}
\arrow[bend right=10]{d}[swap,pos=0.35]{\Cod}
\&
\begin{array}{c}
 \Up\CCc^\LRan/\nLan
\\
 \cong
\\
 \commm{\nKan}
\\
 \cong
\\
 \nRan{} / \Do\CCc^\RLan
\\
\end{array} 
\arrow[bend right=10]{u}[swap,pos=0.3]{\Dom}
\arrow[bend right=10]{d}[swap,pos=0.4]{\Cod}
\\
\Up \CCc \arrow[bend right=45]{uu}[swap]{\mbox{\Large$\Ran$}}
\arrow[phantom]{u}[pos=0.6,description]{\dashv}
\arrow[bend right=10]{u}[swap,pos=0.63]{\varepsilon}
\arrow[thin,two heads]{r} 
\&
\Do\CCc^\RLan 
\arrow[bend right=45]{uu}[swap]{\mbox{\Large$\nLan$}}
\arrow[phantom]{u}[pos=0.6,description]{\dashv}
\arrow[bend right=10]{u}[swap,pos=0.6]{\varepsilon}
\end{tikzar}%
\end{equation}
The difference is that a  gap in a poset is a pair of a lower set and an upper set such that one is below the other \eqref{eq:commm-Kan}, whereas a gap in a category is not just a {\presheaf} and a {\postsheaf} with some property, but also a set of cones and cocones \emph{witnessing}\/ that the {\presheaf} is below the {\postsheaf}. The witnesses are certain assignments between cones as the lower bounds of {\postsheaves} and cocones as the upper bounds of {\presheaves}, tied together by the adjunctions and represented by two pairs of isomorphic comma categories displayed in \eqref{diag:Kan}. Each of the categories $\commm{\Kan}$ and $\commm{\nKan}$ conveniently conjoins an isomorphic pair representing an adjunction, just as the lattice of cuts $\UD\PPp$ in \eqref{eq:DM} conveniently conjoins the isomorphic lattices of $\RLan$-closed lower sets and $\LRan$-closed upper sets. We write them out for completeness: 
\bea\label{eq:commmKan-obj}
\Big\lvert \commm{\Kan} \Big\rvert & = &  \coprod_{\substack{\lft A\in \Do\CCc\\ \rgt A \in \Up\CCc}}\  
\left\{\ \left<
\begin{matrix}\lft g^A\\ \vspace{-.1ex}\\ \rgt g^A
\end{matrix}\right>\  \in\ \  
\begin{matrix}\Do\CCc\left(\lft A, \Ran\rgt A\right)\\ \times\\ \Up\CCc\left(\Lan\lft A, \rgt A \right)\end{matrix}\ \ \ \Bigg|\ \ \  
\begin{matrix}\lft g^A = \Ran\rgt g^A\circ \eta
\\ \wedge\\ 
\rgt g^A=\varepsilon\circ\Lan\lft g^A
\end{matrix}\ \ 
\right\} 
\\
\commm{\Kan}(A,B) & = &  \left\{
\left<\begin{matrix}\lft f\\ \vspace{-.005ex}\\ \rgt f\end{matrix}\right>\  \in\ \  \begin{matrix}\Do\CCc\left(\lft A, \lft B\right)\\ \times\\ \Up\CCc\left(\rgt A, \rgt B\right)\end{matrix}  \ \Bigg|\ \ \ \begin{matrix}\lft g^B\circ \lft f = \Ran\rgt f \circ \lft g^A
\\ \wedge\\ 
\rgt g^B\circ \Lan\lft f = \rgt f \circ \rgt g^A
\end{matrix}\ \ \right\}\notag
\\[4ex]
\label{eq:commmnKan-obj}
\Big\lvert \commm{\nKan} \Big\rvert & = &  \coprod_{\substack{\left(\rgt \alpha\colon\LRan\rgt A\ot \rgt A\right)\, \in\, \ULR\\
\left(\lft \alpha\colon \RLan\lft A\to \lft A\right)\, \in\, \DRL }} \  
\left\{\ \left<
\begin{matrix}\rgt \intv^\alpha\\ \vspace{-.005ex}\\ \lft \intv^\alpha
\end{matrix}\right>\  \in\ \  
\begin{matrix}\ULR\left(\rgt \alpha, \nLan\lft \alpha\right)\\ \times\\ \DRL\left(\nRan\rgt \alpha, \lft \alpha \right)\end{matrix}\ \ \ \Bigg|\ \ \  
\begin{matrix}\rgt \intv^\alpha = \nLan\lft \intv \circ \rgt\alpha
\\ \wedge\\ 
\lft\intv^\alpha=\lft \alpha \circ \nRan\rgt\intv^\alpha
\end{matrix}\ \ 
\right\}
\\[3ex]
\commm{\nKan}(\aggg,\bggg) & = &  \left\{
\left<\begin{matrix}\rgt f\\ \vspace{-.01em}\\ \lft f\end{matrix}\right>\  \in\ \  \begin{matrix} \ULR\left(\rgt \aggg, \rgt \bggg\right)\\ \times\\ \DRL\left(\lft \aggg, \lft \bggg\right) \end{matrix}\ \ \Bigg|\ \ \ 
\begin{matrix}\rgt f\circ  \rgt \intv^\beta = \lft \intv^\alpha \circ \nLan \lft f
\\ \wedge\\ 
\lft \intv^\beta\circ \nRan \rgt f = \lft f\circ \lft \intv^\alpha
\end{matrix}\ \
\right\} \notag
\eea
In all of their isomorphic avatars, the objects of these categories are the categorical versions of gaps and intervals, respectively. Since the intuitions vary, we describe categorical gaps in some detail in Sec.~\ref{Sec:catgap}, and categorical intervals in Sec.~\ref{Sec:catint}.

\subsubsection{The Isbell adjunction and the category of gaps}\label{Sec:catgap}

If we follow up from the idea of Sec.~\ref{Sec:Lambek-problem} and think of a given {\presheaf} $\lft A:\lft \AAa\to \CCc$ as a lower set, then the posetal notion of an upper bound lifts to the categorical notion of a cocone $\delta \in \Coc(\lft A, c) =  \Do\CCc(\lft A, \mnd c)$. Analogously, when a {\postsheaf} $\rgt A:\rgt \AAa\to \CCc$ is thought of as an upper set, then the lower bounds lift to the cones $\varrho \in\Con(c,\rgt A) = \Up\CCc(\cmn c, \rgt A) = \pU(\rgt A, \cmn c)$. If a posetal gap $$\lft X\Big)\cdots \Big(\rgt Y$$ is a lower set $\lft X$ and an upper set $\rgt Y$ where every $x\in \lft X$ is a lower bound of $\rgt Y$ and every $y\in \rgt Y$ is an upper bound of $\lft X$, then a categorical gap must comprise a {\presheaf} and a {\postsheaf} with a natural family of arrows witnessing that each element of the {\presheaf} is below each element of the {\postsheaf}.

\begin{definition}\label{Def:gap}
A \emph{gap}\/ in a category $\CCc$ is a triple $A=<\lft A, \rgt A, \gap^A>$, usually written in the form  $$\GAP{\lft A} {\gap} {\rgt A}$$ 
where $\lft A\in \Do\CCc$ and $\rgt A\in \Up\CCc$ are {\action}s and $\gap \in \Do(\CCc\times \CCc^\op)\left(\lft A\times \rgt A^\op\, , \widetilde \CCc\right)$ is an {\action} homomorphism
\beq\label{eq:gap-pic}
\begin{tikzar}[row sep = 3em, column sep = 1em]
\lft \AAa \times \rgt \AAa^\op \ar[bend left=10]{rr}{\gap} \ar{dr}[description]{\lft A\times \rgt A^\op}
\&\& \tarrow{\CCc}\ar{dl}[description]{\widetilde \CCc}\\
\&\CCc\times \CCc^\op
\end{tikzar}
\eeq
i.e., a family of maps $\left\{\gap_{xy}:\lft Ax \to \rgt Ay\ |\ x\in \lft \AAa, y\in \rgt \AAa\right\}$ such that
\begin{enumerate}[a)]
\item  $\left\{\gap_{ay}:\lft Aa \to \rgt Ay\ |\ y\in \rgt \AAa\right\}$ is a cone for each $a\in \lft \AAa$, whereas 
\item  $\left\{\gap_{xb}:\lft Ax \to \rgt Ab\ |\ x\in \lft \AAa\right\}$ is a cocone for each $b\in \rgt \AAa$.
\end{enumerate}
The set of gaps is written in the form
\beq\label{def:Gaps}
\Gaps\left(\lft A, \rgt A\right) \ \  = \ \  \left \{\ \ 
\begin{tikzar}[sep=small]
{\lft A\big)\!\!}\arrow[Rightarrow]{r}{\gap}\& {\!\!\big(\rgt A}
\end{tikzar}\ \ \right\}\ \ =\ \ \Do\left(\CCc\times \CCc^\op\right)\left(\lft A\times \rgt A^o\,,\  \widetilde{\CCc}\right) 
\eeq
where $\widetilde\CCc=\left(\tarrow\CCc\tto{<\Dom,\Cod>} \CCc\times \CCc^\op \right)$ is the distributor representing $\CCc$, as described in Appendix~\ref{Appendix:comma}.
\end{definition}

\begin{definition}\label{Def:isbellcat}
The \emph{category of gaps}\/ in $\CCc$ is
\bea\label{eq:isbmon}
\Big\lvert \isbell\CCc \Big\rvert & = &  \coprod_{\substack{\lft A\in \Do\CCc\\ \rgt A \in \Up\CCc}} 
\Gaps\left(\lft A, \rgt A\right) 
\\
\isbell\CCc(A,B) & = &  \left\{\ \ 
\left<\begin{matrix}\lft f\\ \hspace{1em}\\ \rgt f\end{matrix}\right>\  \in\ \  \begin{matrix}\Do\CCc\left(\lft A, \lft B\right)\\ \times\\ \Up\CCc\left(\rgt A, \rgt B\right)\end{matrix}  
\ \ \Bigg|\ \ 
\begin{tikzar}[row sep=3em]
\lft A\ar{r}{\lft f} \ar{d}[description]{\lft\gap^A} \&\lft B \ar{d}[description]{\lft \gap^B}\\
\Ran\rgt A \ar{r}[swap]{\Ran\rgt f} \& \Ran \rgt B
\end{tikzar}\ \ 
\right\}
\notag
\eea
\end{definition}

\paragraph{Remark.} The requirement on the $\isbell\CCc$-morphisms can be equivalently written in the dual form, and also as a chu-morphism \cite{BarrM:chu96,PavlovicD:chuI,PrattV:chu} over  categories.
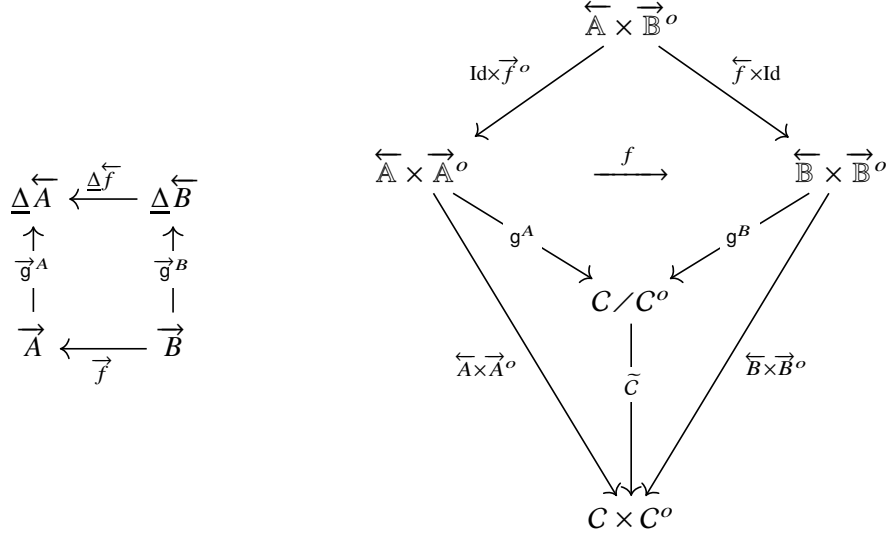
\[
\begin{tikzar}[row sep=3em]
\Lan \lft A\ar[leftarrow]{r}{\Lan\lft f} \ar[leftarrow]{d}[description]{\rgt\gap^A} \&\Lan\lft B \ar[leftarrow]{d}[description]{\rgt \gap^B}\\
\rgt A \ar[leftarrow]{r}[swap]{\rgt f} \&\rgt B
\end{tikzar}
\qquad\qquad\qquad
\begin{tikzar}[row sep=3em,column sep=3em]
\& \lft \AAa \times \rgt \BBb^\op \arrow{dl}[swap]{\Id \times \rgt f^o} \arrow{dr}{\lft f \times \Id} \\
\lft \AAa \times \rgt \AAa^o \arrow{dr}[description]{\gap^A}\arrow{dddr}[swap]{\lft A\times \rgt A^o} 
\&
\tto{\ \ \ f\ \ \ }
\& 
\lft \BBb \times \rgt \BBb^\op \arrow{dl}[description]{\gap^B} \arrow{dddl}{\lft B\times \rgt B^\op}
\\
\&\tarrow{\CCc}\arrow{dd}[description,pos=0.35]{\widetilde\CCc} \\ \\
\& \CCc\times \CCc^o
\end{tikzar}
\] 

\paragraph{Background and terminology.} The $\isbell{}$-construction induces the \emph{Isbell monad}\/ $\isbell:\Cat\to \Cat$, which was introduced and studied in detail by Richard Garner in \cite{GarnerR:isbell}. Garner calls the objects of $\isbell\CCc$ \emph{cylinders}. In earlier literature, for example \cite{Freyd-Kelly}, the term "cylinder" was used to refer to singly indexed families $Fx\to Gx$, rather than doubly indexed in the form $Fx\to Gy$. While the term "bicone" might be more descriptive, in the present context the "gap" intuition seems helpful.

The correspondences from \eqref{eq:LR} now lift to categories. 

\begin{lemma} There are natural bijections
\beq\label{eq:gapcomma}
\Do\CCc(\lft A, \Ran\rgt A) \ \ \cong \ \ \Gaps(\lft A, \rgt A) \ \ \cong \ \ \Up\CCc(\Lan\lft A, \rgt A)
\eeq
\end{lemma}

\bpr
By \eqref{eq:cones-cocones}, fibers of the {\presheaf}\ $\Ran\rgt A$ are cones under $\rgt A$, and fibers of the {\postsheaf} $\Lan\lft A$ are cocones over $\lft A$. The claim follows from Def.~\ref{Def:gap}. The naturality check is straightforward.
\epr

\begin{proposition}
The correspondences in \eqref{eq:gapcomma} induce the equivalence $\isbell\CCc\simeq \commm{\Kan:\Up\CCc\to\Do\CCc}$.
\end{proposition}

\paragraph{Summary.} The isomorphisms $\comm{\Do\CCc}{\Ran}\  \cong\  \commm{\Kan}\ \cong\  \comm{\Lan}{\Up\CCc}$ on the left in \eqref{diag:Kan}, and in refined form in \eqref{eq:gapcomma}, are realized by the bijective correspondences that can be depicted as follows
\beq\label{eq:corresp-gap}
\begin{split}
\prooftree
\prooftree
\lft A\Big) \tto{\lft \gap} \Ran\rgt A\Big)\Big(\rgt A
\Justifies
\GAP{\lft A}{\gap}{\rgt A}
\thickness = 2.5ex
\endprooftree
\Justifies
\lft A\Big)\Big(\Lan\lft A \oot{\rgt \gap} \Big(\rgt A
\thickness = 2.5ex
\endprooftree
\qquad\qquad\qquad\qquad
\begin{tikzar}[row sep = 1ex]
\Lan\rgt A \&\&  \rgt A\ar{ll}[swap,pos=0.4]{\rgt \gap} 
\ar{dd}[description]{\varepsilon}
\\
\Ran\Lan\lft A\ar{ddrr}[description,pos=0.4]{\Ran\lft \gap}
\\
\&\& \Lan\Ran\rgt A \ar{uull}[description,pos=0.4]{\Lan\rgt \gap}
\\
\lft A\ar{uu}[description,pos=0.45]{\eta} \ar{rr}[pos=0.4]{\lft \gap}\&\& \Ran\rgt A
\end{tikzar}
\end{split}
\eeq 
A morphism $f\in\isbell\CCc(A,B)$ shows how the gap $A$ contains the gap $B$, roughly in the form 
$$\lft A\Big)\tto{\lft f}\GAP{\lft B}{\gap^B}{\rgt B}\oot{\rgt f} \Big(\rgt A$$

\subsubsection{The algebra-coalgebra adjunction and the category of intervals}\label{Sec:catint}

An interval in a category can be construed as a pair of gaps
\bear 
\INTVL{\rgt A}{\gapp}{\lft A} & = & \left<\ \begin{matrix}\GAP{\Ran\rgt A}{\gapp_\bullet}{\Lan\lft A}\\
\GAP{\Ran\rgt A}{\gapp^\bullet}{\Lan\lft A}\end{matrix}\ \ \right>
\eear 
related by the following bijective correspondence
\beq\label{eq:corresp-intvl}
\begin{split}
\prooftree
\prooftree
\Ran\rgt A\Big) \tto{\lft \gapp_\bullet} \Ran\Lan\lft A\Big)\tto{\lft\alpha}\lft A\Big)\Big(\Lan \lft A
\Justifies
\INTVL{\rgt A}{\gapp}{\lft A}
\thickness = 2ex
\endprooftree
\Justifies
\Ran\rgt A\Big)\Big(\rgt A \oot{\rgt \alpha}\Big(\Lan\Ran\rgt A \oot{\rgt\gapp^\bullet}\Big(\Lan\lft A
\thickness = 2ex
\endprooftree
\qquad\qquad\qquad
\begin{tikzar}[row sep = 1ex]
\Lan\Ran\lft A \&\&  \Lan \lft A\ar{ll}[swap,pos=0.4]{\rgt \gapp^\bullet} 
\ar{dd}[description]{\Lan\lft\alpha}
\\
\Ran\Lan\Ran\rgt A\ar{ddrr}[description,pos=0.4]{\Ran\rgt \gapp^\bullet}
\\
\&\& \Lan\Ran\Lan \lft A \ar{uull}[description,pos=0.4]{\Lan\lft \gapp_\bullet}
\\
\Ran \rgt A\ar{uu}[description,pos=0.45]{\Ran\rgt\alpha} \ar{rr}[pos=0.4,swap]{\lft \gapp_\bullet}\&\& \Ran\Lan\lft A
\end{tikzar}
\end{split}
\eeq 
Note that the correspondence is dual to \eqref{eq:corresp-gap}, with the role of the adjunction unit $\eta$ and counit $\varepsilon$ played, respectively by $\Ran\rgt \alpha$ and $\Lan\lft \alpha$. The fact that they realize a bijective correspondence was suggested by Prop.~\ref{Prop:Cnuc}. 

\paragraph{Explanation.} Recall that an interval in a poset was presented in Sec.~\ref{Sec:gapint} as a gap between
\begin{itemize}
\item  the lower bounds $\Ran\rgt Y$ under a closed upper set $\rgt Y = \Lan\Ran\rgt Y$ and
\item the upper bounds $\Lan\lft X$ over a closed lower set $\lft X = \Ran\Lan\lft X$.
\end{itemize}
While a posetal interval between closed $\rgt Y$ and $\lft X$ was a gap between $\Ran \rgt Y$ and $\Lan\lft X$
\[
\Ran \rgt Y \Big)\Big(\rgt Y\cdots \lft X\Big)\Big( \Lan \lft X
\]
with $\bigwedge \rgt Y = \bigvee\Ran \rgt Y$ and $\bigwedge \lft X = \bigwedge \Lan\lft X$, in categories
\begin{itemize}
\item lower sets lift to {\presheaves}\ $\lft A$ with closure witnessed by algebra structures $\Ran\Lan\lft A\tto{\lft\alpha}\lft A$, while
\item upper sets lift to {\postsheaves}\ $\rgt A$, with closure witnessed by coalgebra structures $\rgt A\oot{\rgt\alpha}\Lan\Ran\rgt A$. 
\end{itemize}
A categorical interval between a closed {\postsheaf} $\rgt A$ and closed {\presheaf} $\lft A$ is therefore not just a gap between $\Ran\rgt A$ and $\Lan\lft A$, but depends on a coalgebra $\rgt \alpha$ on $\rgt A$ and an algebra $\lft \alpha$ on $\lft A$, as specified in \eqref{eq:corresp-intvl} on the right. This dependency is spelled out in the following definition.

\begin{definition}\label{Def:interval}
An \emph{interval}\/ in a category $\CCc$ is a tuple $\alpha = \left< \lft A,\rgt A, \lft\alpha, \rgt\alpha,\gapp_\bullet,\gapp^\bullet\right>$, usually written in the form
\[
\INTVL{\rgt \alpha}{\gapp} {\lft \alpha}
\]
where $\left(\Ran\Lan\lft A\tto{\lft \alpha}\lft A\right)\in \DRL$ is a $\RLan$-algebra, $\left(\rgt A\oot{\rgt\alpha}\Lan\Ran\rgt A\right)\in \ULR$ is a $\LRan$-coalgebra, and $\gapp_\bullet, \gapp^\bullet \in \Do\left(\CCc\times \CCc^\op\right)\left( \Ran\rgt A\times\Lan\lft A, \widetilde \CCc\right)$ are gaps as in Def.~\ref{Def:gap}, such that for all $\varrho \in \Con(x,\rgt A)$ and $\delta \in \Coc(\lft A,y)$
\beq\label{eq:interval}
\gapp_\bullet(\varrho, \delta) = \left(x\tto{\varrho\lft\alpha\left(\rgt\gapp^\bullet \delta\right)}y\right) 
\qquad\qquad
\gapp^\bullet(\varrho, \delta) = \left(x\tto{\delta\lft\alpha\left(\lft \gapp_\bullet \varrho\right)}y\right)
\eeq
The set of intervals  in a category $\CCc$ is thus
\beq\label{def:Intervals}
\Intvls\left(\lft \alpha, \rgt \alpha\right)  \  = \   \left \{\  
\begin{tikzar}[sep=small]
{\big(\lft \alpha\!\!}\arrow[Rightarrow]{r}{\gapp}\& {\!\!\rgt \alpha\big)}
\end{tikzar}\ 
\right\}\  =\  \left\{\gapp = \left<\gapp_\bullet,\gapp^\bullet\right>\in \Gaps\left(\Ran\rgt A, \Lan\lft A\right)\\ \Big|\ \eqref{eq:interval}\ \right\} 
\eeq
\end{definition}

\begin{lemma}\label{Lemma:interval-one}
Condition \eqref{eq:interval} is equivalent to
\beq\label{eq:lemma-interval-one}
\lft \gapp_\bullet(\varrho) = \left(\Lan\lft A \tto{\rgt \gapp^\bullet} \Lan\Ran\rgt A \tto{\rgt\alpha} \rgt A\tto{\varrho} \Lan x\right) 
\qquad\qquad
\rgt \gapp^\bullet(\delta) = \left(\Ran\rgt A \tto{\lft \gapp_\bullet} \Ran\Lan\lft A \tto{\lft\alpha} \lft A\tto{\delta} \Ran y\right)
\eeq
These equations are respectively equivalent to the commutativity of the second and of the first triangle in \eqref{eq:corresp-intvl}.
\end{lemma}

\begin{proposition}\label{lemma:interval-two}
The gaps $\gapp_\bullet, \gapp^\bullet \in \Do\left(\CCc\times \CCc^\op\right)\left( \Ran\rgt A\times\Lan\lft A, \widetilde \CCc\right)$ that satisfy
\beq\label{eq:lemma-interval-two}
\gapp_\bullet(\varrho, \delta) = \varrho\lft\alpha\left(\lambda \widetilde\varrho.\ \delta\rgt\alpha\left(\lft\gapp_\bullet \widetilde\varrho\right)\right)
\qquad\qquad
\gapp^\bullet(\varrho, \delta) = \delta\rgt\alpha\left(\lambda \widetilde\delta.\ \varrho\lft\alpha\left(\rgt\gapp^\bullet \widetilde\delta\right)\right)
\eeq
also satisfy \eqref{eq:lemma-interval-one} and thus form an interval.
\end{proposition}

\bpr
Substituting in \eqref{eq:corresp-intvl} the $\Ran$-image of the first triangle into the second one, and the $\Lan$-image of the second triangle into the first one, yields the following commutative squares
\beq\label{eq:proof-interval-two}
\begin{tikzar}[row sep = 4em, column sep = 3.5em]
\Ran\Lan\Ran\rgt A \ar{r}{\Ran\Lan\lft\gapp_\bullet} \ar[leftarrow]{d}[description]{\Ran\rgt\alpha}\ar{dr}[description]{\Ran\rgt\gapp^\bullet}\& \Ran\Lan\Ran\Lan\lft A\ar{d}[description]{\Ran\Lan\lft\alpha}
\\
\Ran\rgt A \ar{r}[swap]{\lft \gapp_\bullet} \& \Ran\Lan\lft A
\end{tikzar}
\hspace{8em}
\begin{tikzar}[row sep = 4em, column sep = 3.5em]
\Lan\Ran\rgt A  \& 
\Lan\lft A \ar{l}[swap]{\rgt \gapp^\bullet} 
\\
\Lan\Ran\Lan\Ran\rgt A
\ar{u}[description]{\Lan\Ran\rgt\alpha}
\& 
\Lan\Ran\Lan\lft A \ar[leftarrow]{u}[description]{\Lan\lft \alpha}
\ar{l}{\Lan\Ran\rgt\gapp^\bullet} \ar{ul}[description]{\Lan\lft\gapp_\bullet}
\end{tikzar}
\eeq
Going around each of the squares corresponds to one of the equations in \eqref{eq:lemma-interval-two}. This shows that a pair of gaps satisfying \eqref{eq:corresp-intvl} and thus forming an interval also satisfy \eqref{eq:lemma-interval-two}. The diagonals of the squares in  \eqref{eq:proof-interval-two} show how either of the gaps one of the conditions in \eqref{eq:lemma-interval-two} determines the other gap with which it forms an interval. 
\epr

%

\paragraph{When are the gaps $\gapp_\bullet$ and $\gapp^\bullet$ equal?} It is easy to see that $\gapp_\bullet = \gapp^\bullet$ if and only if $\lft\gapp_\bullet$ and $\rgt \gapp^\bullet$ are each other's transposes along $\Kan$, i.e. if they satisfy not only  \eqref{eq:corresp-intvl} but also \eqref{eq:corresp-gap}. We will see in Sec.~\ref{Sec:catcut} that this only happens for very special {\action}s. 

\begin{proposition}\label{Prop:jay}
For any $\RLan$-algebra $\Ran\Lan\lft A\tto{\lft \alpha}\lft A$  and any $\LRan$-coalgebra $\rgt A\oot{\rgt \alpha}\Lan\Ran\rgt A$ there are bijections
\bea\label{eq:intvlcomma}
\DRL\left(\nRan \rgt \alpha, \lft \alpha\right) \ \   \cong &  \Intvls\left(\rgt \alpha, \lft \alpha\right) & \cong\ \ \ULR\left(\rgt \alpha, \nLan\lft \alpha\right)\\
\left(\Ran\rgt A\tto{\lft\gapp_\bullet} \Ran\Lan\lft A\tto{\lft \alpha}\lft A \right)=\lft \intv\ \  \mapsfrom & \left<\gapp_\bullet, \gapp^\bullet\right>  & \mapsto \ \ \rgt  \intv = \left(\rgt A\oot{\rgt \alpha}\Lan\Ran\rgt A \oot{\rgt\gapp^\bullet}\Lan \lft A\right)\notag\\
\lft \intv\ \  \mapsto & \left<\Ran\rgt \intv, \Lan\lft \intv\right>  & \mapsfrom \ \ \rgt  \intv\notag
\eea
\end{proposition}

\bprf{ sketch} The correspondence between $\lft \intv$ and $\rgt \intv$ as transposes along $\nKan$ was displayed in \eqref{eq:bigadj}, modulo some renamings\footnote{The $\lft \alpha$ and $\rgt \alpha$ in \eqref{eq:proof-interval-two} correspond to $\alpha$ and $\beta$ in \eqref{eq:bigadj}. In Sec.~\ref{Sec:Cnuc}, there was no justification to call both an algebra and a coalgebra $\alpha$ and to distinguish them by the arrows on top. In the current section, such  algebra-coalgebra pairs will appear as components in multi-component objects. Instead of carrying around a name for each particular component, we shall write $\lft \alpha$ and $\rgt \alpha$ for the components of an object $\alpha$.}. Extending the left-hand square in \eqref{eq:proof-interval-two} by $\lft \alpha$ and the right-hand square by $\rgt \alpha$ shows that the correspondence of $\gapp_\bullet$ and $\gapp^\bullet$ in \eqref{eq:intvlcomma} induces the correspondence of $\lft \intv$ and $\rgt \intv$ in \eqref{eq:bigadj}. The squares in \eqref{eq:bigadj} have $\Ran\varepsilon$ and $\Lan\eta$ where the squares in \eqref{eq:proof-interval-two} have $\Ran\rgt \alpha$ and $\Lan\lft \alpha$, but the chase showing that one version commutes if and only if the other version commutes is routine, and can be found in \cite[Lemma~7.3]{PavlovicD:nucleus}.
\epr

\begin{definition}\label{Def:intercat}
The \emph{category of intervals}\/ in $\CCc$ is  
\bea\label{eq:lambek-obj}
\lvert \NNN\CCc\rvert & = &  \coprod_{\substack{\left(\rgt \alpha\colon\LRan\rgt A\ot \rgt A\right)\, \in\, \ULR\\
\left(\lft \alpha\colon \RLan\lft A\to \lft A\right)\, \in\, \DRL }} 
\Intvls\left(\rgt \alpha,\lft \alpha\right) 
\\[3ex]
\NNN\CCc(\aggg,\bggg) & = &  \left\{
\left<\begin{matrix}\rgt f\\ \hspace{1em}\\ \lft f\end{matrix}\right>\  \in\ \  \begin{matrix} \ULR\left(\rgt \aggg, \rgt \bggg\right)\\ \times\\ \DRL\left(\lft \aggg, \lft \bggg\right) \end{matrix}\ \ \Bigg|\ \ \ 
\begin{tikzar}[row sep=3em,column sep = 3.5em]
\Ran\rgt A \ar{r}{\Ran\rgt f} \ar{d}[description]{\lft  \gapp^\aggg}  
\& \Ran\rgt B\ar{d}[description]{\lft\gapp^\bggg}
\\
\Ran\Lan \lft A \ar{r}[swap]{\Ran\Lan \lft f}\& \Ran\Lan\lft B
\end{tikzar}\ \ 
\right\} \notag
\eea
\end{definition}

\paragraph{Remark.} The the morphism condition can again be equivalently written in a transposed and in a chu form
{\small\[
\begin{tikzar}[row sep=3em,column sep = 3.5em]
\Lan\Ran\rgt A \ar[leftarrow]{r}{\Lan\Ran\rgt f} \ar[leftarrow]{d}[description]{\rgt  \gapp^\aggg}  
\& \Lan\Ran\rgt B\ar[leftarrow]{d}[description]{\rgt\gapp^\bggg}
\\
\Lan \lft A \ar[leftarrow]{r}[swap]{\Lan \lft f}\& \Lan\lft B
\end{tikzar}
\qquad\qquad
\begin{tikzar}[row sep=3em,column sep=.2em]
\& \Con(\rgt A) \times \Coc(\lft B)^\op \arrow{dl}[swap]{\Id \times \Lan\lft f^\op} \arrow{dr}{\Ran\rgt f \times \Id} \\
\Con(\rgt A) \times \Coc(\lft A)^\op \arrow{dr}[description]{\gapp^A}\arrow{dddr}[swap]{\Ran\rgt A\times \Lan\lft A^\op} 
\&
\& 
\lft \Con(\rgt B) \times \Coc(\lft B)^\op \arrow{dl}[description]{\gapp^B} \arrow{dddl}{\Ran\rgt B\times \Lan\lft B^o}
\\
\&\tarrow{\CCc}\arrow{dd}[description,pos=0.35]{\widetilde\CCc} \\ \\
\& \CCc\times \CCc^o
\end{tikzar}
\]}%
where $\Con(\rgt A)$ and $\Coc(\lft A)$ are, respecrively, the categories of cones under the {\postsheaf} $\rgt A$ and of cocones over the {\presheaf} $\lft A$. They were previously used in Sec.~\ref{Sec:Lambek-problem}.

\paragraph{Background.} While the $\JJJ$-construction from Def.~\ref{Def:isbellcat} originates from the \emph{$\JJJ$sbell monad} \cite{GarnerR:isbell}, the $\NNN$-construction in Def.~\ref{Def:intercat} is closely related with the \emph{$\NNN$ucleus construction} and the \emph{Street monad} \cite{PavlovicD:nucleus}. 

\begin{proposition}
The bijections \eqref{eq:intvlcomma} induce the equivalence $\NNN\CCc \simeq \commm{\nKan:\DRL\to \ULR}$.\hspace{-1ex}
\end{proposition}

%

\subsection{Cuts in categories}\label{Sec:catcut}

A cut in a poset, going back to Dedekind's definition of the reals, is a gap which is at the same time an interval, as constructed in \eqref{diag:intersection-P}, which means that
\begin{itemize}
\item  a gap $\lft X\Big)\cdots \Big(\rgt Y$ assures $\bigvee \lft X \leq \bigwedge \rgt Y$, whereas 
\item an interval $\mnd \rgt Y\Big)\Big(\rgt Y\cdots \lft X\Big)\Big( \cmn \lft X$  assures $\bigwedge \rgt Y \leq \bigvee \lft X $, so that
\item a cut presents an element as a supremum and an infimum at once, since $\bigvee \lft X = \bigwedge \rgt Y$. 
\end{itemize}
Lambek's Problem was to generalize such constructions to categories, with suprema and infima replaced by colimits and limits. Isbell showed that this cannot be achieved in a universal way for all limits and colimits, and we are trying to find the maximal family for which it can be achieved. In this section we define categorical cuts: gaps that are also intervals.

\subsubsection{Full cuts}\label{Sec:full}
\begin{definition}
A \emph{full cut}\/ in a category $\CCc$ is a tuple $\alpha=\left<\lft A, \rgt A, \lft \alpha, \rgt \alpha, \gap^\alpha, \gapp^\alpha\right>$, where
\begin{itemize}
\item $\lft A\in \Do\CCc$ and $\rgt A\in \Up\CCc$ are {\action}s,
\item $\left(\Ran\Lan\lft A\tto{\lft \alpha}\lft A\right)\in \DRL$ is a $\RLan$-algebra, $\left(\rgt A\oot{\rgt\alpha}\Lan\Ran\rgt A\right)\in \ULR$ is a $\LRan$-coalgebra,
\item $\GAP{\lft A}{\gap^A}{\rgt A}$ is a gap,
\item $\INTVL{\rgt \alpha}{\gapp^\alpha}{\lft \alpha}$ is an interval
\end{itemize}
such that the following diagrams commute
\beq\label{eq:cut}
\begin{tikzar}[row sep = 3em,column sep = 3em]
\& \Ran\Lan\lft A \ar{dr}[description]{\lft\alpha} \ar{dl}[description]{\Ran \rgt\gap}\\
\Ran\rgt A\ar{r}[description]{\lft\gapp_\bullet} \& \Ran\Lan\lft A \ar{r}[description]{\lft \alpha}\& \lft A\\
\& \lft A\ar{ul}[description]{\lft\gap} \ar[equals]{ur}
\end{tikzar}
\qquad 
\qquad
\begin{tikzar}[row sep = 3em,column sep = 3em]
\& \Lan\Ran\rgt A \ar{dr}[description]{\Lan\lft\gap} \ar{dl}[description]{\rgt \alpha}\\
\rgt A\ar[leftarrow]{r}[description]{\rgt\alpha} \& \Lan\Ran\rgt A \ar[leftarrow]{r}[description]{\rgt\gapp^\bullet}\& \Lan\rgt A\\
\& \rgt A\ar[equals]{ul} \ar{ur}[description]{\rgt\gap}
\end{tikzar}
\eeq
\end{definition}

\paragraph{Remark.} The two requirements in \eqref{eq:cut} are not independent: the bottom triangle of one and the top triangle of the other imply the complementary pair of triangles. To see this, first note that the middle lines of \eqref{eq:cut} are the homomorphisms from \eqref{eq:intvlcomma}
{\small
\[ \lft \intv = \left(\Ran\rgt A\tto{\lft\gapp_\bullet} \Ran\Lan\lft A\tto{\lft \alpha}\lft A \right)\,  \in\,  \DRL(\nRan\rgt \alpha,\lft \alpha) \ \ \ \ \ \mbox{and}\ \ \  \ \  
\rgt \intv = \left(\rgt A\oot{\rgt \alpha}\Lan\Ran\rgt A \oot{\rgt\gapp^\bullet}\Lan \lft A\right)\,  \in\,  ^\LRan\hspace{-0.1ex}\pU\CCc(\nLan\lft \alpha, \rgt\alpha)
\]}%
Prop.\ref{Prop:jay} established that these homomorphisms are each other's transposes along $\nKan:\left(^\LRan\pU\CCc\right)^\op\to \DRL$, where we write $^\LRan\hspace{-0.07ex}\pU\CCc$ instead of\footnote{$\pU\CCc$ denotes the category of {\postsheaf}s and equivariant homomorphisms, so that $\Up\CCc = \left(\pU\CCc\right)^\op$. See \eqref{eq:DopUUp}.} $\left(\pU\CCc\right)^{\LRan^\op}$. The left and the right square in \eqref{eq:cut} now correspond respectively to the left and right pair of equations:
\begin{align}
\lft \intv \circ \lft\gap &=\   \id & \rgt \intv \circ \Lan\lft \gap &= \rgt \alpha \label{eq:rgt}\\
\lft \intv \circ \Ran\rgt \gap &= \lft \alpha & \rgt \intv \circ \rgt \gap &= \id \label{eq:lft}
\end{align}
The top pair of these equations turns out to be equivalent to the bottom pair:
 
\begin{proposition}
The pair of equations in \eqref{eq:rgt} is equivalent to the pair \eqref{eq:lft}.
\end{proposition}

\bpr
Equations \eqref{eq:rgt} say that the large rectangle and the right-hand triangle commute in the following diagram.
\beq\label{eq:seventy}
\begin{tikzar}[row sep = 4em,column sep = 5em]
\Ran\Lan\Ran\Lan\lft A \ar[two heads]{d}[description]{\Ran\varepsilon}\ar[two heads,pos=0.4]{r}[description]{\Ran\Lan\Ran\rgt \gap}\& \Ran\Lan\Ran\rgt A \ar[two heads]{r}[description]{\Ran\Lan\lft \intvl}
 \ar[two heads]{d}[description]{\Ran\varepsilon}
 \& \Ran\Lan\lft A  \ar[two heads]{d}[description]{\lft \alpha}
 \\
 \Ran\Lan \lft A \ar[two heads]{r}[description]{\Ran\rgt\gap}
 \&
 \Ran \rgt A 
 \ar[two heads]{r}[description]{\lft \intvl} \ar[tail]{ur}[description]{\Ran \rgt \intvl}
 \& \lft A\\
 \& \lft A \ar[tail]{ul}[description]{\eta} \ar[tail]{u}[description]{\lft\gap} \ar[equals]{ur}
\end{tikzar}
\eeq
Equations \eqref{eq:lft} say that the left-hand triangle and the large rectangle  commute in the following diagram.
\beq
\label{eq:seventyone}
\begin{tikzar}[row sep = 4em,column sep = 5em]
\Lan\Ran\rgt A  \ar[two heads]{d}[description]{\rgt \alpha}
 \& 
 \ar[two heads]{l}[description]{\Lan\Ran\rgt \intvl}\Lan\Ran\Lan\lft A  \ar[two heads]{d}[description]{\Lan\eta}\& \Lan\Ran\Lan\Ran\rgt A
 \ar[two heads]{d}[description]{\Lan\eta}  \ar[two heads,pos=0.4]{l}[description]{\Lan\Ran\Lan\lft \gap} 
 \\
\rgt A 
 \&
 \Lan \lft A 
 \ar[two heads]{l}[description]{\rgt \intvl} \ar[tail]{ul}[description]{\Lan \lft \intvl}
 \& 
 \Lan\Ran\rgt A \ar[two heads]{l}[description]{\Lan\lft \gap}
 \\
 \& \rgt A \ar[tail]{ur}[description]{\varepsilon} \ar[tail]{u}[description]{\rgt\gap} \ar[equals]{ul}
\end{tikzar}
\eeq 
The adjunction correspondence between $\lft \intv$ and $\rgt \intv$, established in Prop.\ref{Prop:jay}, and displayed in the general form in \eqref{eq:bigadj}, ensures that the square with the diagonal in \eqref{eq:seventyone} commutes if and only if the square with the diagonal in \eqref{eq:seventy} commutes. The adjunction correspondence between $\lft \gap$ and $\rgt \gap$, as the adjoint views of the gap, ensures that the bottom left triangle in \eqref{eq:seventy} commutes if and only if the bottom right triangle commutes in \eqref{eq:seventyone}. The squares without the diagonal commute in both diagrams by the naturality. It follows that each of the diagrams commutes if and only if the other does. Since \eqref{eq:seventy} commutes if and only if the equations in \eqref{eq:rgt} hold, and \eqref{eq:seventyone} commutes if and only if the equations in \eqref{eq:lft} hold, this proves the claim.
\epr

\begin{corollary}\label{Corollary:eyes} 
A cut $\alpha=\left<\lft A, \rgt A, \lft \alpha, \rgt \alpha, \gap^\alpha, \gapp^\alpha\right>$ makes $\lft A$ into a retract of $\Ran \rgt A$ and $\rgt A$ into a retract of $\Lan \lft A$. The (co)algebra structures $\lft\alpha$ and $\rgt\alpha$ factor through these retractions:
\beq\label{eq:eyes}
\begin{tikzar}[row sep = 3em,column sep = 4.5em]
 \Ran\Lan \lft A \ar[two heads,bend left=15]{r}[description]
{\Ran\rgt\gap}
\ar[two heads,bend left=40]{rr}[description]{\lft \alpha}
 \&
 \Ran \rgt A \ar[two heads,bend left=15]{r}[description]{\lft \intv} \ar[tail,bend left = 15]{l}[description]{\Ran \rgt \intv}
 \& 
 \lft A \ar[tail,bend left=15]{l}[description]
{\lft\gap} \ar[tail,bend left=40]{ll}[description]
{\ideta} 
\& 
\rgt A \ar[tail,bend left=15]{r}[description]{\rgt \gap} \ar[tail,bend left=40]{rr}[description]{\idepsilon}
 \&
 \Lan \lft A \ar[two heads,bend left=15]{l}[description]{\rgt \intv} \ar[tail,bend left=15]{r}[description]{\Lan\lft \intv}
 \& 
 \Lan\Ran\rgt A \ar[two heads,bend left=15]{l}[description]{\Lan\lft \gap} \ar[two heads,bend left=40]{ll}[description]{\rgt \alpha}
\end{tikzar} 
\eeq
where $\ideta$ is the transpose of the idempotent $\Lan\lft A\eepi{\rgt \intv} \rgt A \mmono{\rgt \gap} \Lan\lft A$, whereas $\idepsilon$ is the transpose of $\Ran\rgt A\eepi{\lft \intv} \lft A \mmono{\lft \gap} \Ran\rgt A$.
\end{corollary}

\begin{definition}\label{Def:cutcat}
The \emph{category of full cuts}\/ in $\CCc$ is 
{\small 
\bea\label{eq:cut-full-def}
\lvert \LlL\CCc\rvert & = &  \coprod_{\substack{\left(\rgt \alpha\colon\LRan\rgt A\ot \rgt A\right)\, \in\, \ULR\\
\left(\lft \alpha\colon \RLan\lft A\to \lft A\right)\, \in\, \DRL }} 
\left\{\left<\begin{matrix}\gap \\ \hspace{1em}\\ \gapp\end{matrix}\right>\  \in\ \  \begin{matrix} \Gaps\left(\lft A, \rgt A\right)\\ \times\\ \Intvls\left(\rgt \alpha,\lft \alpha\right) \end{matrix}\ \ \Bigg|\ \ \ 
\begin{tikzar}[row sep=2em,column sep=1em]
\Ran\Lan\lft A \ar{rr}[description]{\Ran\rgt\gap}\ar{dd}[description]{\lft\alpha} \&\&
\Ran\rgt A \ar{dl}[description]{\lft\gapp_\bullet} 
\&
\rgt A \ar[leftarrow]{dl}[description]{\rgt\alpha} \ar[equals]{dd}
\\
\& \Ran\Lan\rgt A \ar{dl}[description]{\lft \alpha}
\& \Lan\Ran\lft A \ar[leftarrow]{dl}[description]{\rgt \gapp^\bullet} 
\\
\lft A  
\&
\Lan\lft A \ar[leftarrow]{rr}[description]{\lft\gap}  \&\&\rgt A
 \end{tikzar}
\right\}
\\[3ex]
\LlL\CCc(\aggg,\bggg) & = &  \left\{
\left<\begin{matrix}\rgt f\\ \hspace{1em}\\ \lft f\end{matrix}\right>\  \in\ \  \begin{matrix} \ULR\left(\rgt \aggg, \rgt \bggg\right)\\ \times\\ \DRL\left(\lft \aggg, \lft \bggg\right) \end{matrix}\ \ \ \Bigg|\ \ \ \ 
\begin{tikzar}[row sep=3em,column sep = 3.5em]
\lft A \ar{d}[description]{\lft f} \ar{r}[description]{\lft\gap^A} \&\Ran\rgt A \ar{d}[description]{\Ran\rgt f} \ar{r}[description]{\lft  \gapp_\bullet^A}  
\& \Ran\Lan \lft A \ar{d}[description]{\Ran\Lan \lft f} 
\\
\lft B  \ar{r}[description]{\lft\gap^B}\& \Ran\rgt B\ar{r}[description]{\lft\gapp_\bullet^B}\& \Ran\Lan\lft B
\end{tikzar}
\ \ 
\right\} \notag
\eea
}
\end{definition}

\paragraph{Dual condition for morphisms.} As noted after Definitions \ref{Def:isbellcat} and \ref{Def:intercat}, the morphism requirements can be equivalently stated in dual form, and both dual versions can be expressed together in the chu form. This time, moreover, each of the morphism components determines the other:
\beq\label{eq:twomorph}
\begin{tikzar}[row sep = 8ex,column sep = large]
\Ran\rgt A \ar{d}[description]{\Ran\rgt f}\ar[two heads]{r}[description]{\lft \intv^A} \& \lft A \ar[dashed]{d}[description]{\lft f} \ar[tail]{r}[description]{\lft\gap^A} \& \Ran\rgt A \ar{d}[description]{\Ran\rgt f}
\\
\Ran\rgt B \ar[two heads]{r}[description]{\lft \intv^B} \& \lft B \ar[tail]{r}[description]{\lft\gap^B} \& \Ran\rgt B
\end{tikzar}
\qquad\qquad
\begin{tikzar}[row sep = 8ex,column sep = large]
\Lan\lft A \ar{d}[description]{\Lan\lft f}\ar[two heads]{r}[description]{\rgt \intv^A} \& \rgt A \ar[dashed]{d}[description]{\rgt f} \ar[tail]{r}[description]{\rgt\gap^A} \& \Lan\lft A \ar{d}[description]{\Lan\lft f}
\\
\Lan\lft B \ar[two heads]{r}[description]{\rgt \intv^B} \& \rgt B \ar[tail]{r}[description]{\rgt\gap^B} \& \Lan\lft B
\end{tikzar}
\eeq
%
%

\noindent\textbf{If $\gapp_\bullet=\gapp^\bullet$ then $\Ran\rgt A\cong \lft A$  and $\rgt A\cong \Lan\lft A$.} As mentioned in Sec.~\ref{Sec:bicnuc}, an interval $\gapp = <\gapp_\bullet, \gapp^\bullet>$ satisfies $\gapp_\bullet= \gapp^\bullet$ if and only if the {\presheaf} version $\lft\gapp_\bullet$ and the {\postsheaf} version $\rgt \gapp^\bullet$ are each other's transposes along $\Kan$, and thus satisfy not only \eqref{eq:corresp-intvl} but also \eqref{eq:corresp-gap}. One of the conditions in \eqref{eq:corresp-gap} says that the northeast triangle in the following diagram commutes.
\beq\begin{tikzar}[column sep = 8em, row sep = 5em]
\Ran\rgt A\ar[tail]{r}{\eta} \ar[tail]{dr}[description]{\lft\gapp_\bullet = \Ran\rgt\intv} 
\ar[two heads]{d}[description]{\lft \intv}
\& \Ran\Lan\Ran\rgt A\ar[two heads]{d}[description]{\Ran\rgt\gapp^\bullet = \Ran\Lan\lft \intv}\\
\lft A \ar[tail]{r}[swap]{\eta} \& \Ran\Lan\lft A
\end{tikzar}
\eeq
Since the square commutes by the naturality of $\eta$, the southwest triangle must commute as well. If the interval is a component of a cut, then $\rgt\intv$ is a split epi in $\pU\CCc$ and $\Ran\rgt \intv$ is thus a split monic in $\Do\CCc$. Since $\Ran\rgt\intv = \eta\circ\lft \intv$, it follows that $\lft \intv$ must also be a split monic. Since by \eqref{eq:rgt} $\lft \intv$ is a split epi in any cut, it follows that it must be an isomorphism. An interval in a cut thus satisfies $\gapp_\bullet= \gapp^\bullet$ if and only if $\lft \intv$ is an isomorphism. A symmetric argument proves the same for $\rgt \intv$. Corollary~\ref{Corollary:eyes} then yields the isomorphisms $\Ran\Lan\lft A\cong \lft A$ and $\rgt A\cong \Lan\Ran\rgt A$. An exercise, for example, using \eqref{eq:cones-cocones}, shows that this is satisfied when $\lft A$ and $\rgt A$ are retracts of representables. The converse seems to be true by \emph{reductio ad absurdum}, but it may not be true constructively.

\subsubsection{Simple cuts}\label{Sec:simple}

\begin{definition}
A \emph{simple cut}\/ in a category $\CCc$ is a tuple $A=\left<\lft A, \rgt A, \lft\gap, \rgt\gap, \lft\intv, \rgt\intv\right>$, where
\begin{itemize}
\item $\lft A\in \Do\CCc$ and $\rgt A\in \Up\CCc$ are {\action}s,
\item $\lft \gap\in \Do\CCc\left(\lft A, \Ran\rgt A\right)$ and $\rgt \gap\in \pU\CCc\left(\rgt A, \Lan\lft A\right)$ are transposes along $\Kan$, and thus form a gap,
\item $\lft\intv \in \Do\CCc\left(\Ran\rgt A, \lft A\right)$ and $\rgt \intv\in \pU\CCc\left(\Lan\lft A, \rgt A\right)$ are transposes along $\nKan$ (and thus form an interval, see Prop.~\ref{Prop:cutcut}), and they make the following diagrams commute 
\bea\label{eq:simple-jay}
\begin{tikzar}[column sep=3em]
\Ran\Lan\Ran\rgt A \ar[shift left=2]{r}{\Ran\varepsilon}
\ar[shift right=2]{r}[swap]{\Ran \idepsilon}
\&
 \Ran {\rgt A} \ar{r}[description]{\lft \intv} \&  {\lft A}  \end{tikzar}
& \qquad &
\begin{tikzar}[column sep=3em]
\rgt A \&
 \Lan {\lft A} \ar{l}[description]{\rgt \intv} \& \ar[shift left=2]{l}{\Lan\ideta} \ar[shift right=2]{l}[swap]{\Lan \eta}
 \Lan\Ran\Lan {\lft A}  \end{tikzar}
\eea
where
\beq\label{eq:id-def}\idepsilon = \left(\Lan\Ran\rgt A\oot{\Lan\lft\intv} \Lan\lft A\oot{\rgt\gap} \rgt A\right)
\qquad\qquad\qquad\ideta = \left(\lft A \tto{\lft \gap} \Ran\rgt A \tto{\Ran\rgt \intv} \Ran\Lan\lft A\right) 
\eeq
\end{itemize}
such that $\gap$s and $\intv$s form retractions: 
\beq\label{eq:cut-simple}
\begin{tikzar}[row sep = .25em,column sep = 5em]
\& \lft A \ar[equals]{dd} \ar[tail]{dl}[description]{\lft\gap}\& \& \rgt A \ar[equals]{dd} \ar[tail]{dr}[description]{\rgt\gap}\\
\Ran\rgt A\ar[two heads]{dr}[description]{\lft\intv} \&\&\&\& \Lan\lft A\ar[two heads]{dl}[description]{\rgt\intv}\\
\& \lft A \&\& \rgt A
\end{tikzar}
\eeq
The set of simple cuts between $\lft A$ and $\rgt A$ is written
\bear
\Cuts\left(\rgt A, \lft A\right) & = & \left\{\left<\lft\gap^A, \rgt \gap^A, \lft \intv^A, \rgt \intv^A\right> \right\}
\eear
\end{definition}

\begin{proposition}\label{Prop:cutcut}
Every simple cut determines a unique full cut, and vice versa.
\end{proposition}

\bpr
A simple cut $A=\left<\lft A, \rgt A, \lft\gap, \rgt\gap, \lft\intv, \rgt\intv\right>$ induces a full cut $\alpha=\left<\lft A, \rgt A, \lft \alpha, \rgt \alpha, \gap, \gapp\right>$ comprising the same-named components\footnote{This means that $\lft A$ and $\rgt A$ are the same and $\lft\gap$ and $\rgt\gap$ are folded into $\gap$.} and moreover
\begin{align} 
\lft \alpha =& \left(\Ran\Lan\lft A \tto{\Ran\rgt\gap}\Ran\rgt A \tto{\lft \intv}\lft A\right) & 
\rgt \alpha = & \left(\rgt A \oot{\rgt \intv} \Lan\lft A \oot{\Lan\lft \gap}\Lan\Ran\rgt A\right)\label{eq:alphadef}
\\
\lft \gapp_\bullet =& \left(\Ran\rgt A \tto{\Ran\rgt \intv} \Ran\Lan\lft A\right) & 
\rgt \gapp^\bullet =& \left(\Lan\Ran\rgt A \oot{\Lan\lft \intv} \Lan\lft A\right)\label{eq:Jay-bullet-def}
\end{align}
We must prove that $\lft \alpha$ and $\rgt \alpha$ satisfy the requirements for, respectively, a $\RLan$-algebra and a $\LRan$-coalgebra and that, together with $\lft\gapp_\bullet$ and $\rgt\gapp^\bullet$, they also satisfy \eqref{eq:cut}. The conditions that make $\lft \alpha$ from \eqref{eq:alphadef} into a $\RLan$-algebra are in the following diagram.
\beq\label{eq:proof-simple}
\begin{tikzar}[column sep = 4.5em, row sep = 4em]
\Ran\Lan\Ran\Lan\lft A \ar{r}[description]{\Ran\Lan\Ran\rgt \gap} \ar[bend right = 30]{dd}[description]{\mu\ =\Ran\varepsilon} 
\ar[bend left = 20]{rr}{\Ran\Lan\lft \alpha}
\&
\Ran\Lan\Ran\rgt A\ar{r}[description]{\Ran\Lan\lft\intv} \ar[bend right = 30]{dd}{\Ran\varepsilon} 
\&
\Ran\Lan\lft A\ar{d}[description]{\Ran\rgt \gap}
\ar[bend right = 30]{dd}[swap]{\lft \alpha}
\\
\&\& \Ran\rgt A \ar{d}[description]{\lft \intv}
\&
\lft A\ar{ul}[swap]{\eta} \ar{l}[description]{\lft \gap}
\ar[equals]{dl}
\\
\Ran\Lan\lft A \ar{r}[description]{\Ran\rgt \gap} 
\ar[bend right = 20]{rr}[swap]{\lft \alpha}
\&
\Ran\rgt A \ar{r}[description]{\lft \intv} \& \lft A\end{tikzar}
\eeq 
Chasing the triangles on the right gives
\beq\label{eq:alphaeta}
\lft \alpha \circ \eta\  =\  \lft \intv \circ \Ran\rgt \gap \circ \eta \ =\  \lft \intv \circ \lft \gap\  =\  \id
\eeq
The requirement that $\lft \alpha \circ \Ran\Lan\lft \alpha = \lft \alpha \circ \mu$ corresponds to the commutativity of the outermost quadrilateral of arcs, consisting of two inner quadrilaterals. The left-hand one commutes by the naturality of $\varepsilon$. It remains to be proved that the right-hand quadrilateral also commutes:
\beq
\lft \alpha \circ \Ran\Lan\lft\intv\  \stackrel{\eqref{eq:alphadef}}=\  \lft\intv\circ \Ran\gap\circ \Ran\Lan\lft\intv \ \stackrel{\eqref{eq:id-def}}=\ \lft \intv \circ \Ran\idepsilon\ \stackrel{\eqref{eq:simple-jay}}=\  \ \lft\intv\circ\Ran\varepsilon
\eeq
The proof that $\rgt \alpha$ is a $\LRan$-coalgebra is dual. Towards a proof that $\lft \alpha$ and $\rgt \alpha$, together with $\lft \gapp_\bullet$ and $\rgt \gapp^\bullet$ defined in \eqref{eq:Jay-bullet-def}, satisfy \eqref{eq:cut}, note that Prop.~\ref{Prop:jay} implies that setting $\lft \gapp_\bullet$ and $\rgt \gapp^\bullet$ as in \eqref{eq:Jay-bullet-def} returns
\beq\label{eq:jayJay}
\lft \intv\ =\ 
\left(\Ran\rgt A\tto{\lft\gapp_\bullet} \Ran\Lan\lft A\tto{\lft \alpha}\lft A \right)\qquad \qquad 
\rgt  \intv\  =\  \left(\rgt A\oot{\rgt \alpha}\Lan\Ran\rgt A \oot{\rgt\gapp^\bullet}\Lan \lft A\right)
\eeq 
Conditions \eqref{eq:cut} now reduce to \eqref{eq:cut-simple} and \eqref{eq:alphadef}. Conversely, given a full cut $\alpha=\left<\lft A, \rgt A, \lft \alpha, \rgt \alpha, \gap, \gapp\right>$, we need to show that the simple cut $A=\left<\lft A, \rgt A, \lft\gap, \rgt\gap, \lft\intv, \rgt\intv\right>$, with $\lft \intv$ and $\rgt \intv$ defined by \eqref{eq:jayJay}, satisfies \eqref{eq:simple-jay} and \eqref{eq:cut-simple}. But the left-hand equation in \eqref{eq:simple-jay} is the right-hand rectangle in \eqref{eq:proof-simple}. This rectangle must commute because the outer square commutes for any algebra $\lft \alpha$, the left-hand rectangle by the naturality of $\varepsilon$, and $\Ran\Lan\Ran\rgt\gap$ is a split epimorphism.  Finally, substituting \eqref{eq:simple-jay} in the bottom triangles of \eqref{eq:cut} yields the triangles of \eqref{eq:cut-simple}. \epr

\begin{definition}\label{Def:cutsimple}
The \emph{category of simple cuts}\/ in $\CCc$ is 
{\small 
\bea\label{eq:cut-simple-def}
\lvert \LLL\CCc\rvert & = &  \coprod_{\substack{
\rgt A\in \Up\CCc\\
\lft A\in \Do\CCc}} 
\ \ \Cuts\left(\rgt A, \lft A\right)\\[3ex]
\LLL\CCc(A, B) & = &  \left\{
\left<\begin{matrix}\rgt f\\ \hspace{1em}\\ \lft f\end{matrix}\right>\  \in\ \  \begin{matrix} \Up\CCc\left(\rgt A, \rgt B\right)\\ \times\\ \Do\CCc\left(\lft A, \lft B\right) \end{matrix}\ \ \ \Bigg|\ \ \ \ 
\begin{tikzar}[row sep = 8ex,column sep = 7ex]
\Ran\rgt A \ar{d}[description]{\Ran\rgt f}\ar{r}[description]{\lft \intv^A} \& \lft A \ar{d}[description]{\lft f} \ar{r}[description]{\lft\gap^A} \& \Ran\rgt A \ar{d}[description]{\Ran\rgt f}
\\
\Ran\rgt B \ar{r}[description]{\lft \intv^B} \& \lft B \ar{r}[description]{\lft\gap^B} \& \Ran\rgt B
\end{tikzar}\ \ 
\right\} \notag
\eea
}
\end{definition}

\paragraph{Redundancies.} The fact that each of the morphism components determines the other is like in \eqref{eq:twomorph}. The equivalent dual condition on the $\LLL\CCc$-morphism is obtained as there. 

A straightforward consequence of Prop.~\ref{Prop:cutcut} is the following:

\begin{proposition}
The category $\LLL\CCc$ of simple cuts (Def.~\ref{Def:cutsimple}) is isomorphic to the category  $\LlL\CCc$ of full cuts (Def.~\ref{Def:cutcat}).
\end{proposition}

%
%

\subsubsection{Absolute cuts}\label{Sec:abs}

\begin{definition}\label{Def:absolute}
An \emph{absolute cut}\/ in a category $\CCc$ is a tuple $\ida =\left<\,\lft A, \rgt A, \lft \ida, \rgt \ida, \lft \varphi, \rgt \varphi\, \right>$ where
\begin{itemize}
\item $\lft A\in \Do\CCc$ and $\rgt A\in \Up\CCc$ are {\action}s,
\item $\lft \ida \in\Do\CCc\left(\Ran\rgt A, \Ran\rgt A\right)$, $\rgt \ida \in\Up\CCc\left(\Lan\lft A, \Lan\lft A\right)$, $\lft \varphi \in\Do\CCc\left(\RLan \lft A, \RLan \lft A\right)$ and $\rgt \varphi \in\Up\CCc\left(\LRan \rgt A, \LRan \rgt A\right)$ are idempotents with joint splittings: 
\end{itemize}
\bea\label{eq:ida}
\begin{tikzar}[column sep=2em]
\Ran\rgt A \ar[bend right = 15,shift right=2]{dd}[swap]{\lft\ida} \ar[two heads]{d}{\lft \intv} \& \Lan\Ran\rgt A \&\& \Lan\Ran\rgt A
\ar{ll}[swap]{\rgt \varphi}
\\
\lft A \ar[tail]{d}{\lft \gap} \&\Lan\lft A  \ar[tail]{u}{\Lan\lft\intv} \& \& \Lan\lft A \ar{ll}[description]{\rgt \ida} 
\ar[tail]{u}[swap]{\Lan\lft\intv}  
\\
\Ran\rgt A   \& \Lan\Ran\rgt A \ar[two heads]{u}{\Lan\lft\gap} \&\& \Lan\Ran\rgt A \ar{ll}{\rgt \varphi
} 
\ar[two heads]{u}[swap]{\Lan\lft\gap}
 \end{tikzar}
& \qquad &
\begin{tikzar}[column sep=2em]
\Lan\lft A   \& \Ran\Lan\lft A \ar{rr}{\lft\varphi} \ar[two heads]{d}[swap]{\Ran\rgt\gap} \&\& \Ran\Lan\rgt A   \ar[two heads]{d}{\Ran\rgt\gap}
\\
\rgt A \ar[tail]{u}[swap]{\rgt \gap}  \&\Ran\rgt A \ar{rr}[description]{\lft \ida}  \ar[tail]{d}[swap]{\Ran\rgt\intv} \& \& \Ran\rgt A  
  \ar[tail]{d}{\Ran\rgt\intv}
\\
\Lan\lft A \ar[bend left = 15,shift left=2]{uu}{\rgt\ida} \ar[two heads]{u}[swap]{\rgt \intv}  \& \Ran\Lan\lft A \ar{rr}[swap]{\lft\varphi}  \&\& \Ran\Lan\lft A  
 \end{tikzar}
\eea
We write $\cuts\left(\lft A, \rgt A\right)$ to denote the set of absolute cuts between $\lft A$ and $\rgt A$.
\end{definition}

\paragraph{Why "absolute"?} A categorical structure is called absolute when it is preserved by all functors.   Appendix~\ref{Appendix:split} explains the sense in which absolute properties are characterized by splittings.

\begin{definition}\label{Def:cutabsol}
The \emph{category of absolute cuts}\/ in $\CCc$ is 
{\small 
\bea\label{eq:cut-abs-def}
\lvert \Labs\CCc\rvert & = &  \coprod_{\substack{\lft A\in \Do\CCc\\
\rgt A\in \Up\CCc}} \
\cuts\left(\lft A, \rgt A\right)\\[3ex]
\Labs\CCc(\ida, \idb) & = &  \left\{\ \left<\begin{matrix}\rgt f\\ \hspace{1em}\\ \lft f\end{matrix}\right>\  \in\ \  \begin{matrix} \Up\CCc\left(\rgt A, \rgt B\right)\\ \times\\ \Do\CCc\left(\lft A, \lft B\right) \end{matrix}\ \ \ \Bigg|
\ \ \ \ 
\begin{tikzar}[row sep = 8ex,column sep = 7ex]
\Ran\rgt A \ar{d}[description]{\Ran\rgt f}\ar{r}[description]{\lft \ida} \& \Ran\rgt A \ar{d}[description]{\Ran\rgt f}
\\
\Ran\rgt B \ar{r}[description]{\lft \idb} \& \Ran\rgt B
\end{tikzar}\  \ \ \ \raisebox{-1ex}{\Large $\wedge$}\ \ \ \  
\begin{tikzar}[row sep = 8ex,column sep = 7ex]
\Lan\lft A \ar[leftarrow]{d}[description]{\Lan\lft  f}\ar[leftarrow]{r}[description]{\rgt \ida} \& \Lan\lft  A \ar[leftarrow]{d}[description]{\Lan\lft  f}
\\
\Lan\lft B \ar[leftarrow]{r}[description]{\rgt \idb} \& \Lan\lft  B
\end{tikzar}
\ \ 
\right\} \notag
\eea
}
\end{definition}

\begin{proposition}\label{Prop:sim-abs-equiv}
The categories of simple and absolute cuts are equivalent.
\end{proposition}

\bpr
For fixed {\action}s $\lft A$ and $\rgt A$, a simple cut induces an absolute cut
\beq\label{eq:simabseq}
\begin{tikzar}[row sep = 1ex]
\LLL\CCc \arrow[bend left = 7]{rrr}
\arrow[phantom]{rrr}[description]{\simeq} \&\&\& \Labs\CCc   \arrow[bend left = 7]
{lll}
\\
\left<\lft \gap, \rgt \gap, \lft \intv, \rgt \intv \right> \ar[mapsto]{rrr} \&\&\& 
\left<\lft \ida, \rgt \ida, \lft \varphi, \rgt \varphi\right>
\end{tikzar}
\eeq
by composing the gaps and the intervals into the idempotents as follows:
\begin{align*}
\lft \ida &= \lft \gap\circ\lft \intv & \lft \varphi &= \Ran\rgt \intv\circ\lft\gap\circ\lft\intv\circ\Ran\rgt \gap \\
\rgt \ida &= \rgt \gap\circ \rgt \intv &  \rgt \varphi &= \Lan\lft \intv\circ\rgt\gap\circ\rgt \intv\circ\Lan\lft \gap
\end{align*}
The fact that $\rgt \ida$ satisfies \eqref{eq:ida}  follows from \eqref{eq:simple-jay}. Conversely, given an absolute cut $\ida = \left<\lft \ida,\rgt\ida,\lft \varphi, \rgt \varphi\right>$, the corresponding simple cut $\alpha = \left<\lft\gap,\rgt \gap, \lft \intv, \rgt \intv\right>$ is defined by the splittings in \eqref{eq:ida} and the composites in \eqref{eq:eyes}. Since \eqref{eq:ida} directly ensures \eqref{eq:simple-jay}, the checks that the defined structure is a simple cut and that an $\Labs\CCc$-morphism from \eqref{eq:cut-abs-def} induces a unique $\LLL\CCc$-morphism from \eqref{eq:cut-simple-def} are straightforward.
\epr

\subsection{Cut algebras and coalgebras}\label{Sec:cutalg}

The adjunction factorization spelled out in Lemma~\ref{Lemma:comma} can for the case of $\nKan$ be tightened from the comma categories to the cut category.

\begin{proposition}\label{Prop:lemma}
For every category $\CCc$ the algebra-coalgebra adjunction $\nKan$ factors through a reflection $\Cod\dashv\varepsilon$ and a coreflection $\eta\dashv \Dom$
\beq\label{eq:nKan-fac}
\begin{tikzar}[column sep = large]
\ULR \arrow[bend right = 30,shift right = 2]{rrrr}[swap]{\nRan}  \arrow[bend right = 15,shift right = 1,tail]{rr}[description]{\eta} \arrow[phantom]{rr}[description]{\top}  \&\&\Labs\CCc\arrow[bend right = 15,shift right = 1,two heads]{rr}[description]{\Cod} \arrow[bend right = 15,shift right = 1,two heads]{ll}[description]{\Dom}
\arrow[phantom]{rr}[description]{\top}
\&\& 
\DRL \arrow[bend right = 30,shift right = 2]{llll}[swap]{\nLan}  \arrow[bend right = 15,shift right = 1,tail]{ll}[description]{\varepsilon}
\end{tikzar}
\eeq
where the functors are
\bear
\eta\left(\rgt A\oot{\rgt\alpha}\Lan\Ran\rgt A\right)  & = & \begin{pmatrix}
\Ran\rgt A & \Ran\Lan\Ran\rgt A \eepi{\Ran\varepsilon}\Ran\rgt A\mmono{\eta\Ran}\Ran\Lan\Ran\rgt A\\
\rgt A & \Lan\Ran\rgt A\iipe{\varepsilon}\rgt A \oonom{\rgt \alpha} \Lan\Ran\rgt A\end{pmatrix}  \\
\Dom\begin{pmatrix}
\lft A & \Ran\rgt A\eepi{\lft \intv} \lft A \mmono{\lft \gap}\Ran\rgt A\\
\rgt A & \Lan\lft A\oonom{\rgt \gap} \rgt A\iipe{\rgt\intv}\Lan\lft A\end{pmatrix} & = & \left(\rgt A\iipe{\rgt \intv}\Lan\lft A \iipe{\Lan\lft \gap} \Lan\Ran\rgt A \right)\\
\varepsilon\left(\Ran\Lan\lft A\tto{\lft\alpha}\lft A\right) & = & \begin{pmatrix}
\lft A & \Ran\Lan\lft A \eepi{\lft \alpha}\lft A\mmono{\eta}\Ran\Lan\lft A\\
\Lan\lft A & \Lan\Ran\Lan\lft A\iipe{\Lan\eta}\Lan\lft A \oonom{\varepsilon\Lan} \Lan\Ran\Lan \lft A\end{pmatrix}  
\\
\Cod\begin{pmatrix}
\lft A & \Ran\rgt A\eepi{\lft \intv} \lft A \mmono{\lft \gap}\Ran\rgt A\\
\rgt A & \Lan\lft A\oonom{\rgt \gap} \rgt A\iipe{\rgt\intv}\Lan\lft A\end{pmatrix} & = & \left(\Ran\Lan\lft A\eepi{\Ran\rgt\gap} \Ran\rgt A \eepi{\lft \intv} \lft A \right)
\eear
\end{proposition}

\paragraph{Remark.} The \emph{algebra-as-coalgebra}\/ \cite{PavlovicD:nucleus} and the \emph{simple nucleus}\/ \cite{PavlovicD:LICS17} constructions  follow from Prop.~\ref{Prop:lemma} as a corollary. The functors $\Dom$ and $\Cod$ in \eqref{eq:nKan-fac} are both faithful, since each of the two components of any $\Labs\CCc$-morphism determines the other one, by \eqref{eq:twomorph} and Prop.~\ref{Prop:sim-abs-equiv}. (They are not full in general.) 
The categories of $\RLan$-algebras and of $\LRan$-coalgebras can now be equivalently presented in terms of absolute-cut idempotents. 

\begin{corollary}
The category of algebras $\DRL$ is equivalent to the full subcategory of $\Labs\CCc$ spanned by objects of the form $\left<\lft A, \Lan\lft A, \left(\Ran\Lan\lft A\tto{\lft\ida}\Ran\Lan\lft A\right), \left(\Lan\Ran\Lan\lft A\oonom{\varepsilon \Lan}\Lan\lft A \iipe{\Lan\eta} \Lan\Ran\Lan\lft A\right)\right>$. 

Dually, the category of coalgebras $\ULR$ is equivalent to the full subcategory of $\Labs\CCc$ spanned by objects of the form $\left<\Ran\rgt A, \rgt A, \left(\Ran\Lan\Ran\rgt A\eepi{\Ran\varepsilon}\Ran\rgt A\mmono{\eta\Ran}\Ran\Lan\Ran\rgt A\right), \left(\Lan\Ran\rgt A\oot{\rgt \ida} \Lan\Ran\rgt A\right)\right>$.
\end{corollary}

\begin{definition}\label{Def:cut-alg}
A \emph{cut-algebra} is a pair $\kappa_{\rgt A}= \left<\rgt A, \kappa:\Ran\rgt A\to \Ran\rgt A\right>$ where $\kappa$ is idempotent and $\Lan\kappa$ splits on $\rgt A$. A \emph{cut-coalgebra} is a pair $\varphi_{\lft A} = \left<\lft A, \varphi:\Lan\lft A\to \Lan\lft A\right>$ where $\varphi$ is idempotent and $\Ran\varphi$ splits on $\lft A$.  The category of \emph{cut-algebras} is
\bea \label{eq:algComon}
|\simUp| & = & \coprod_{\rgt A\in \Up\CCc}\  \left\{\kappa \in \Do\CCc(\Ran\rgt A, \Ran\rgt A) \ \ \big|\ \ \begin{tikzar}[row sep=1.5em,column sep=.25em]
\Ran\rgt A  \ar{dd}[description]{\kappa } \ar{ddrr}[description]{\kappa } 
\&
\Lan \Ran\rgt A  \ar[leftarrow]{ddrr}[description]{\Lan \kappa}  \&\& \rgt A \ar[tail]{ll}
\\ \\
\Ran\rgt A  \ar{rr}[description]{\kappa} \&\& \Ran \rgt A \& \Lan \Ran\rgt A  \ar[two heads]{uu}
 \end{tikzar}\ \ \right\} 
\\[3ex]
\simUp \left(\kappa_{\rgt A}, \chi_{\rgt B} \right) \ & = &\hspace{2em}  \left\{t \in \Up\CCc(\rgt A,\rgt B)\ \Big|\ \ 
\begin{tikzar}[row sep=3em,column sep=3em]
\Ran\rgt A  \ar{r}[description]{\Ran  t} \ar{d}[description]{\kappa} \& \Ran  \rgt B \ar{d}[description]{\chi}  \\ 
\Ran\rgt A  \ar{r}[description]{\Ran  t}\& \Ran  \rgt B 
\end{tikzar}\ 
\right\} \notag
\eea
The category of \emph{cut-coalgebras} over $\CCc$ is
\bea \label{eq:coalgMon}
|\simDo| & = & \coprod_{\lft A\in \Do\CCc}\  \left\{\varphi \in \Up\CCc(\Lan\lft A, \Lan\lft A) \ \ \big|\ \ \begin{tikzar}[row sep=1.5em,column sep=.25em]
\Lan\lft A  \ar[leftarrow]{dd}[description]{\varphi } \ar[leftarrow]{ddrr}[description]{\varphi } 
\&
\Ran \Lan  \lft A  \ar{ddrr}[description]{\Ran \varphi} \ar[two heads]{rr}  \&\& \lft A \ar[tail]{dd}
\\ \\
\Lan\lft A  \ar[leftarrow]{rr}[description]{\varphi} \&\& \Lan\lft A \& \Ran \Lan\lft A  
 \end{tikzar}\ \ \right\} 
\\[3ex]
\simDo \left(\varphi_{\lft A}, \psi_{\lft B}\right) \  & = &\hspace{2em}  \left\{f \in \Do\CCc(\lft A,\lft B)\ \Big|\ \ 
\begin{tikzar}[row sep=3em,column sep=3em]
\Lan\lft A  \ar[leftarrow]{r}[description]{\Lan  f} \ar{d}[description]{\varphi} \& \Lan  \lft B \ar{d}[description]{\psi}  \\ 
\Lan\lft A  \ar[leftarrow]{r}[description]{\Lan  f}\& \Lan  \lft B 
\end{tikzar}\ 
\right\} \notag
\eea
\end{definition}

\begin{proposition}\label{Prop:alg-coalg-simple}
Cut-algebras over $\CCc$ correspond to $\RLan$-algebras, while cut-coalgebras correspond to $\LRan$-coalgebras. There are thus the following equivalences of categories
\beq\label{eq:alg-coalg-simple}
\simUp\  \simeq\  \DRL\qquad \qquad \qquad \qquad \simDo\  \simeq\  \ULR
\eeq
\end{proposition}

\paragraph{Comments.} The equivalences in \eqref{eq:alg-coalg-simple} are detailed in \cite[Prop.~8.1]{PavlovicD:nucleus} and \cite[Thm.~III.2]{PavlovicD:LICS17}, for general adjunctions. The requirement in \eqref{eq:algComon} that the $\Lan$-image of the idempotent $\kappa$ must split on $\rgt A$ implies on one hand that $\rgt A \cong \Lan\lft A$, where $\lft A$ is obtained by splitting $\kappa$ itself. Moreover, it was shown in the above references that this splitting comprises a $\RLan$-algebra and the corresponding unit, i.e.,
\bea 
\kappa & = & \left(\Ran\rgt A \cong \Ran\Lan \lft A\eepi{\lft\alpha}\lft A \mmono{\eta} \Ran\Lan\lft A\cong \Ran\rgt A \right) 
\eea
The requirement is thus that the $\kappa$-image of a cone $\RRR$ indexed over the cocones $\delta\in\Lan\lft A$ must comprise the components $\kappa(\RRR)_\delta \ = \  \delta_{\lft\alpha(\RRR)}$.

\section{Result: The completion}\label{Sec:limits}
\subsection{From loose to tight limits}\label{Sec:lim}

While the categorical limits and colimits are usually defined for diagrams as functors from arbitrary small categories \cite[Ch.~III, Sections 3--4]{MacLaneS:CWM}, it was mentioned in Sec.~\ref{Sec:cat-loose} that they can be defined equivalently in terms of {\action}s. The reductions between the two treatments are summarized in  Appendix~\ref{Appendix:diags}. The limit operations then arise as adjoints of the two Yoneda embeddings: 
\beq\label{eq:lims}
\begin{tikzar}[column sep = 3em]
\CCc \arrow[bend right=15]{rr}[swap]{\mnd} 
\arrow[phantom]{rr}[description]{\scriptstyle \bot} \&\& \Do \CCc \arrow[bend right=15]{ll}[swap]{\supp}
 \&\&\& \CCc \arrow[bend right=15]{rr}[swap]{\cmn}
 \arrow[phantom]{rr}[description]{\scriptstyle \top} \&\& \Up \CCc \arrow[bend right=15]{ll}[swap]{\inff}
\end{tikzar}
\eeq
This view has been used in foundational theories, such as Street-Walters' \emph{Yoneda structures} \cite{Street-Walters:yoneda} and Ross Steet's \emph{cosmoi}\/ \cite{StreetR:cosmoi}.  Recalling that the {\presheaf} homomorphisms in $\Do\CCc(\lft \AAa, \mnd c)$ can be construed as the cocones from $\lft A$ to $c$, while the cones from $c$ to $\rgt A$ are the {\postsheaf} homomorphisms in $\pU\CCc(\rgt A, \cmn c) = \Up\CCc(\cmn c, \rgt A)$, the adjunctions in \eqref{eq:lims} say that limits represent cones whereas colimits represent cocones, which means that there are natural bijections
\beq\label{eq:coc-con-loose}
\CCc\left(\supp \lft A, c\right) \cong \Coc\left(\lft A, c\right)\qquad\qquad\qquad \CCc\left(c, \inff \rgt A\right) \cong \Con\left(c,\rgt A\right)
\eeq
Tight limits yield to a similar presentation.

\begin{definition}\label{Def:tightlim} Tight colimits $\limsupp$ and tight limits $\liminff$ are the adjoints
\beq\label{eq:tight-def}
\begin{tikzar}[column sep = 3em]
\CCc \arrow[bend right=15]{rr}[swap]{\mnd} \arrow[phantom]{rr}[description]{\scriptstyle \bot} \&\& \simDo \arrow[bend right=15]{ll}[swap]{\limsupp}
 \&\&\& \CCc \arrow[bend right=15]{rr}[swap]{\cmn} 
 \arrow[phantom]{rr}[description]{\scriptstyle \top}\&\&\simUp \arrow[bend right=15]{ll}[swap]{\liminff}
\end{tikzar}
\eeq
where $\mnd$ and $\cmn$ map the objects of $\CCc$ to the corresponding representable {\action}s and their identities.
\end{definition}

\paragraph{Comments.} Since by Prop.~\ref{Prop:alg-coalg-simple} we have $\simDo\simeq \ULR$ and $\simUp \simeq \DRL$, tight limits and tight colimits can be equivalently defined in terms of $\LRan$-coalgebras and $\RLan$-algebras. While working with cut-algebras and coalgebras from Def.~\ref{Def:cut-alg} is more convenient for current purposes, the familiar presentations of algebras and coalgebras are often helpful in applications. The present concepts of tight limit and colimit were first studied in \cite{PavlovicD:CALCO15}, where they were called respectively \emph{limit supremum}\/ and \emph{limit infimum}. 

\begin{proposition}\label{Prop:coc-con-tight}
Tight (co)limits represent (co)cones fixed by the induced idempotents. More precisely, for a cut-coalgebra $\varphi_{\lft A}$ and a cut-algebra $\kappa_{\rgt A}$ there are natural bijections
\bea\label{eq:coc-con-tight}
\CCc\left(\limsupp_{\varphi}\lft A, c\right) & \cong & \left\{\delta\in \Coc(\lft A, c) \ |\ \varphi(\delta) = \delta\right\}\\
\CCc\left(c, \liminff_{\kappa}\rgt A \right) &\cong & \left\{\varrho\in \Con(c,\rgt A)\ |\ \kappa(\varrho) = \varrho\right\}
\eea
\end{proposition}

\bpr
By \eqref{eq:algComon}, a morphism $\delta\in \simDo\left(\varphi_{\lft A}, \mnd c\right)$ is of the form
\beq\label{eq:diagcutalg}
\begin{tikzar}[column sep = 3em,row sep = .5em]
\lft A \ar{r}[description]{\delta} \& \mnd c\\
\Lan\lft A\ar{dd}[description]{\varphi}  \& \cmn c \ar{l}[description]{\Lan \delta} \ar{dd}[description]{\id}\\ \\[3ex]
\Lan\lft A \& \cmn c \ar{l}[description]{\Lan \delta}
\end{tikzar}
\eeq
But $\delta\in \Do\CCc\left(\lft A, \mnd c\right)$ is clearly a cocone. By the Yoneda Lemma, its image $\cmn \delta \in \pU\CCc\left(\cmn c, \Lan\lft A\right)$ corresponds to an element of the fiber $\Lan\lft A_c$, and that element is the cocone $\delta$ itself. The commutativity of the square in \eqref{eq:diagcutalg} says that $\varphi (\delta) = \delta$. The story for $\varrho \in \simUp\left(\cmn c, \lft \ida\right)$ is analogous.
\epr

\begin{corollary}
A category that is complete under loose (co)limits is also complete under tight (co)limits. A functor that preserves loose (co)limits also preserves tight (co)limits.
\end{corollary}

\bpr
Prop.~\ref{Prop:coc-con-tight} says that tight colimits are retracts of {\presheaf}s in the form $\Ran\rgt A$. In a cocomplete category, the colimit of any {\presheaf} as a diagram is representable, and the idempotents split (when coequalized with identities). Tight colimits can thus be constructed using such loose colimits. A functors preserving loose colimits must also preserve tight colimits, since the latter can be reduced to the former.\epr

\paragraph{Loose (co)limits are not derivable from tight (co)limits.} We saw in Sec.~\ref{Sec:dissolution} that only the free group actions lie in the image of the extensions $\Ran$ and $\Lan$. Only the projective actions are retracts of the free ones. Any left group action can be obtained as a loose colimit of copies of the representable action of the group on itself, but only the projective group actions can be obtained as tight colimits. The right group actions can be obtained as loose limits, but only the injective ones are obtained as tight limits.

\subsubsection{Tight diagrams and their limits}\label{Sec:tightdiag} 
The standard (loose) limit and colimit operations take a diagram $D$ as input, producing (respectively) a limit cone out of $\inff D$ or colimit cocone into $\supp D$. Without loss of generality, we restrict diagrams to {\postsheaves} $\rgt D$ for $\inff\rgt D$ or to {\presheaves} $\lft D$ for $\supp\lft D$, as explained at the beginning of Sec.~\ref{Sec:lim}. \emph{Tight}\/ limit and colimit take as input not just a suitable diagram but also a cone idempotent for tight limit or a cocone idempotent for tight colimit. More precisely, tight limit $\liminff$ takes as input a pair $\left<\rgt D, \kappa\right>$, where $\rgt D$ is a {\postsheaf} and $\kappa$ is a {\presheaf} idempotent
\beq
\begin{tikzar}[row sep = 6ex]
\Con\left(\rgt D\right) \ar[bend left = 10]{rr}{\kappa}\ar{dr}[description]{\Ran\rgt D} \&\& \Con\left(\rgt D\right) \ar{dl}[description]{\Ran\rgt D}\\
\& \CCc
\end{tikzar}
\eeq
such that $\Lan\kappa$ splits on $\rgt D$, which means that the cocones $\GGG\in \Lan\Ran\rgt D$ with $\Lan\kappa(\GGG) = \GGG$, i.e., $\GGG_\varrho = \GGG_{\kappa(\varrho)}$, must be in a bijective correspondence with the domain $\rgt \DDd$ of the {\postsheaf} $\rgt D$. We call such a pair $\left<\rgt D,\kappa\right>$ a 
\emph{right tight diagram}.

Dually, tight colimit $\limsupp$ takes as input a \emph{left tight diagram}, i.e., a pair $\left<\lft D, \varphi\right>$, where $\lft D$ is a {\presheaf} and $\varphi$ is a {\postsheaf} idempotent 
\beq
\begin{tikzar}[row sep = 6ex]
\Coc\left(\lft D\right) \ar[bend left = 10,leftarrow]{rr}{\varphi}\ar{dr}[description]{\Lan\lft D} \&\& \Coc\left(\lft D\right) \ar{dl}[description]{\Lan\lft D}\\
\& \CCc
\end{tikzar}
\eeq
such that $\Ran\varphi$ splits on $\lft D$, which means that the cones $\RRR\in \Ran\Lan\lft D$ with $\Ran\varphi(\RRR) = \RRR$, i.e., $\RRR_\delta = \RRR_{\varphi(\delta)}$, are in a bijective correspondence with the domain $\lft \DDd$ of the {\presheaf} $\lft D$. 

\begin{proposition}
$\simDo$ is $\limsupp$-complete. $\simUp$ is $\liminff$-complete.
\end{proposition}

\bpr
We construct the left adjoint in 
\beq\label{eq:tight-left-adjoint}
\begin{tikzar}[column sep = 3em]
\simDo \arrow[bend right=15]{rr}[swap]{\mnd} \arrow[phantom]{rr}[description]{\scriptstyle \bot} \&\& \simsimDo \arrow[bend right=15]{ll}[swap]{\limsupp}
\end{tikzar}
\eeq
An object of $\simsimDo$ is a left tight diagram $\left<\lft D, \varphi\right>$ where $\lft D:\lft \DDd\to \simDo$ is a {\presheaf} in $\Do\simDo$ and $\varphi:\Lan\lft D\to\Lan\lft D$ is a {\postsheaf} idempotent in $\Up\simDo$. The nodes of the diagram $\lft D$ are cut-coalgebras $\left<\lft D_i, \varphi_i \right>$, indexed by $i\in \lft \DDd$, where $\lft D_i:\lft \DDd_i \to \CCc$ are {\presheaves} in $\Do\CCc$ and $\varphi_i$ are {\postsheaf} idempotents in $\Up\CCc$ splitting on $\varphi_i = \left(\Lan\lft D_i\epi \rgt D_i \mono \Lan\lft D_i\right)$ so that $\Ran\rgt D_i = \lft D_i$. Writing $\rgt D:\lft \DDd \to \pU\CCc$ for the diagram mapping $i\in \lft \DDd$ to $\rgt D_i$, we have $\rgt A = \inff \rgt D$ with a limit cone
\bea
\rrrr\colon \rgt A &\tto{ \ \ \ } & \rgt D
\eea
Each component $\rrrr_i: \rgt A\to \rgt D_i$ induces $\Ran\rrrr_i : \lft D_i \to \lft A$, where $\lft A = \Ran\rgt A$ and $\lft D_i = \Ran\rgt D_i$ as before. Hence the cocone $\Ran\rrrr:\lft D\to \lft A$ in $\Do\CCc$. The claim is now that the sought tight colimit in $\simDo$ is $\left<\lft A, \psi\right>  =  \liminff_{\varphi} \lft D$, with the cocone
\bea \dddd \ = \ \varphi\left( \Ran\rrrr\right) & : &\left<\lft D,\varphi\right> \to \left<\lft A, \psi\right>\eea
for $\psi:\Lan\lft A \to \Lan\lft A$ induced by $\varphi:\Lan\lft D\to \Lan \lft D$  by composing
\bear 
\psi_z & = & \Bigg( \Coc_{\Do\CCc}\left(\lft A, z\right)\cong \Coc_{\simDo}\left(\lft D, \mnd z\right) \tto{\varphi_{\mnd z}}\Coc_{\simDo}\left(\lft D, \mnd z\right) \cong \Coc_{\Do\CCc}\left(\lft A, z\right)
\Bigg)
\eear
where $\Coc_{\Do\CCc}\left(\lft A, z\right)\cong \Coc_{\simDo}\left(\lft D,\mnd z\right)$ remains to be checked. Since $\lft A = \Ran\rgt A$, any cocone $\delta:\lft A\cong \Ran\rgt A\tto{\ \ \ } z$ in $\CCc$ must be in the form $\delta = \Ran\varrho$ for a unique cone $\varrho: z \tto{\ \ \ } \rgt A$. Since $\rgt A = \inff \rgt D$, $\varrho$ corresponds to a unique cone $\rrrr\circ\varrho: \cmn z \tto{\ \ \ }\rgt D$ in $\pU\CCc$. By construction of $\rgt D$, the cocone $\Ran(\rrrr\circ\varrho): \lft D  \to \mnd z$ is unique, and lies in $\simDo$. To see that the described map from $\Coc_{\Do\CCc}\left(\lft A,z\right)$ to $\Coc_{\simDo}\left(\lft D,\mnd z\right)$ has an inverse, use $\lft D=\Ran\rgt D$ and $\lft A = \Ran\rgt A$ and the fact that $\rgt A$ is a limit again. Checking that the cut-coalgebra structures are preserved is straightforward, using the construction of $\rgt D$ by splitting the cut-coalgebra structures, and reasoning as in the proof of Prop.~\ref{Prop:coc-con-tight}.

\epr

\subsubsection{Representables generate tightly}\label{Sec:repres-gen}
\begin{proposition}
Any cut-coalgebra $\varphi_{\lft A} \in \simDo$ is a tight colimit of the representables in the image of $\mnd$ and any cut-algebra $\kappa_{\rgt B} \in \simUp$ is a tight limit of representables in the image of $\cmn$, in the form
\bea
\varphi_{\lft A} = \limsupp_{\mnd\varphi}\left(\lft \AAa\tto{\lft A} \CCc\tto{\mnd}\simDo\right)\qquad \qquad \qquad \qquad \kappa_{\rgt B} = \liminff_{\cmn\kappa}\left(\rgt \BBb\tto{\rgt B}\CCc\tto{\cmn}\simUp\right)
\eea
\end{proposition}

\bpr
Any {\presheaf} $\lft A \in \Do\CCc$ is a loose colimit of the form $\lft A = \supp\mnd \lft A$. Thus we have the cocone $r\in \Coc\left(\mnd \lft A, \lft A\right)$ in $\Do\CCc$. The cocone $\varrho = \varphi(r) \in \Coc\left(\mnd \lft A, \lft A\right)$ also satisfies $\varphi(\varrho) = \varrho$, and each of its components is a $\simDo$-morphism, since it makes a diagram like \eqref{eq:diagcutalg} commute. So we have a cocone $\varrho \in \Coc\left(\left<\mnd \lft A,\mnd\varphi\right>, \varphi_{\lft A}\right)$ in $\simDo$, where $\left<\mnd \lft A,\mnd\varphi\right>$ is a tight diagram.  
It remains to be proved that $\varrho$ is a colimit cocone  in $\simDo$. Let $\gamma \in \Coc\left(\left<\mnd \lft A,\mnd \varphi\right>, \xi_{\lft X}\right)$ be an arbitrary tight cocone in $\simDo$, i.e., a cocone $\gamma \in \Coc\left(\mnd \lft A,\lft X\right)$ such that for every cocone $\delta \in \Coc(\lft X)$ with $\xi(\delta) = \delta$ holds $\varphi(\delta\circ \gamma) = \delta\circ \gamma$. It follows that $\varphi(\delta\circ\gamma) = \xi(\delta) \circ \gamma$, which implies the commutativity of the square in the following diagram
\beq\label{eq:factoring}
\begin{tikzar}[row sep=small,column sep = small]
\lft A \ar{rrrrrr}{\widehat \gamma} \&\&\&\&\&\& \lft X\\
\Lan \lft A\&\&\&\&\&\& \Lan \lft X \ar{llllll}[description]{\Lan \widehat \gamma}\\
\& \xi(\delta)\circ\widehat\gamma \&\&\&\& \xi(\delta) \ar[mapsto]{llll}\\[-2.5ex]
\& =
\\[-3ex]
\&\varphi(\delta\circ\widehat\gamma)
\\
\\
\\
\& \delta\circ \widehat\gamma\ar[mapsto]{uuu} \&\&\&\& \delta \ar[mapsto]{llll}\ar[mapsto]{uuuuu}\\
\Lan \lft A\ar{uuuuuuu}[description]{\varphi}\&\&\&\&\&\& \Lan \lft X \ar{llllll}[description]{\Lan \widehat \gamma} \ar{uuuuuuu}[description]{\xi}
\end{tikzar}
\eeq
where $\widehat\gamma\in \Do\CCc(\lft A, \lft X)$ is the factoring of $\gamma \in \Coc\left(\mnd \lft A, \lft X\right)$ through the colimit $\varrho\in \Coc\left(\mnd\lft A,\lft A\right)$ in $\Do\CCc$, i.e.,
\bea\label{eq:gammafac}
\gamma & = &\left(\mnd \lft A \tto\varrho\lft A\tto{\widehat \gamma} \lft X\right)\eea
Diagram \eqref{eq:factoring} says that $\widehat\gamma$ is also a $\simDo$-morphism, and \eqref{eq:gammafac} lifts to $\simDo$ as
\bea\label{eq:gammafac-sim}
\gamma & = &\left(\left<\mnd \lft A,\mnd\varphi\right> \tto\varrho\varphi_{\lft A}\tto{\widehat \gamma} \xi_{\lft X}\right)\eea
This proves $\varphi_{\lft A} = \limsupp_{\mnd\varphi} \mnd\lft A$. The proof for $\kappa_{\rgt B} = \liminff_{\cmn\kappa} \cmn\rgt B$ is dual.
\epr
%
%
%

\subsection{Loose limit completion: the Isbell monad}\label{Sec:isbellmon}

The category of gaps construction from Sec.~\ref{Sec:catgap} can be construed as the answer to the question of whether the Yoneda embeddings $\mnd:\CCc\to \Do\CCc$ and $\cmn:\CCc\to\Up\CCc$, giving rise to the colimit and the limit operations as their adjoints \eqref{eq:lims}, can be conjoined so that both the colimits and the limits arise from the same operation. Garner formulated this answer as his \emph{Isbell (pseudo)monad} $\Isb{}:\CAT\to \CAT$ in the eponymous paper \cite{GarnerR:isbell}. If the Yoneda embeddings are viewed as the units of the two extensively used completion monads \cite{Street-Walters:yoneda,StreetR:cosmoi}, the unit of the Isbell monad follows.
\begin{gather}\label{eq:isbell-unit}
\begin{tikzar}[column sep = 3.cm,row sep = 2cm]
\& \Do \CCc \arrow[bend right = 15,dotted]{dl}[swap]{\scriptstyle\supp}
\arrow[bend right=10,tail]{dr}[swap,pos=0.52]{\eta}
\arrow[phantom]{dr}[sloped]{\scriptstyle\top}
\\
\CCc\arrow{rr}[description,pos=.3]{\Diamond}[swap]{\scriptstyle\bot}{\scriptstyle\bot} 
\arrow{ur}[sloped]{\scriptstyle \bot}[swap]{\mnd} 
\arrow{dr}{\cmn}[swap,sloped]{\scriptstyle \bot} 
\&\& 
\Isb\CCc \arrow[bend right = 10,dotted]{ll}[swap]{\scriptstyle\supp} \arrow[bend left = 10,dotted]{ll}{\scriptstyle\inff}
\arrow[bend right=10,two heads]{ul}[swap]{\scriptstyle\Dom}
\arrow[bend right=10,two heads]{dl}[swap,pos=0.52]{\scriptstyle\Cod}
\\
\& \Up \CCc 
\arrow[bend right=10,tail]{ur}[swap]{\varepsilon}
\arrow[phantom]{ur}[sloped]{\scriptstyle\bot}
\arrow[bend left = 15,dotted]{ul}{\scriptstyle\inff}
\end{tikzar}
\end{gather}
The dotted arrows hint at the algebras for the monads $\Do, \Up\colon\CAT\to\CAT$, which are of Kock-Z\"oberlein type \cite{Kelly-Lack:property,KockA:zoeberlein,PavlovicD:chuI,StreetR:fib-bicat,Zoeberlein}, and thus have to be adjoint, and essentially unique. In terms of Sec.~\ref{Sec:catgap}, an algebra $\Xi\colon \Isb\CCc\to \CCc$ for the Isbell monad maps a gap $\GAP{\lft A} {\gap} {\rgt A}$ from $\Isb\CCc$ to an object $x$ in $\CCc$ which lies between $\lft A$ and $\rgt A$, in the sense of the following diagram.
\beq
\begin{tikzar}[row sep=4em,column sep=5em]
\lft \AAa\times \rgt{\AAa}^\op \arrow{dr}[swap]{\lft A\times \rgt A^o}  \arrow{r}[description]{\lft \Xi\times \rgt\Xi^\op} \arrow[bend left=30]{rr}[description]{\gap} \& \CCc/x\times x/\CCc^o \arrow{r}[description]{\bullet} \arrow{d}[description]{\hspace{.5em}\mnd x \times \cmn x^o} \& \tarrow \CCc \arrow{dl}{\widetilde \CCc}\\
\& \CCc\times \CCc^o
\end{tikzar} 
\eeq
With the adjunction equipment described in \cite{GarnerR:isbell}, the $\Isb{}$-algebra $\Xi$ thus determines a cocone $\lft \Xi$ from $\lft A$ to $x$ and a cone $\rgt \Xi$ from $x$ to $\rgt A$ such that $\lft \Xi\bullet\rgt \Xi = \gap$ \eqref{eq:gap-pic}. Garner proves in \cite[Thm.~12]{GarnerR:isbell} that such algebras equip $\CCc$ with his  \emph{cylindric factorisation systems}. The two extremal cases are when $\lft \Xi\in \Coc\left(\lft A, x\right)$ is a colimit cocone, i.e., $x= \supp \lft A$, and when $\rgt \Xi\in \Con\left(x,\rgt A\right)$, i.e., $x= \inff \rgt A$. When they exist, a factorization chooses an $x$ to decompose the gap through
\[\begin{tikzar}
\supp \lft A \ar{r}\& x \ar{r}  \& \inff \rgt A
\end{tikzar}
\] 
When they do not exist, an $\Isb{}$-algebra chooses $x$ as their joint approximation. The number of such choices can be thought of as the length of the gap. That is the sense in which factoring the gap is a loose approximation of the limit and the colimit. The idea of the tight limit and colimit is that they should tightly approximate each other, which will be achieved in the next section.

\subsection{Tight limit completion: the Lambek monad}
\label{Sec:lambekmon}

In this section, we reach the goal and prove that the cut category $\UD\CCc$ is the tight completion of $\CCc$. The embedding $\between\colon\CCc\to \UD\CCc$ is induced by the Yoneda embeddings again, extended to cut-algebras and cut-coalgebras.
\begin{gather}\label{eq:lambek-unit}
\begin{tikzar}[column sep = 3cm,row sep = 2cm]
\& \simDo \arrow[bend right = 15,dotted]{dl}[swap]{\scriptstyle\limsupp}
\arrow[bend right=10,two heads]{dr}[swap]{\scriptstyle H}
\arrow[phantom]{dr}[sloped,description]{\scriptstyle\top}
\\
\CCc\arrow{rr}[description,pos=.3]{\between}[pos=.55]{\scriptstyle\simeq} 
\arrow{ur}[sloped]{\scriptstyle \bot}[swap]{\mnd} 
\arrow[shift right = .5ex]{dr}[swap]{\cmn}
\arrow[phantom,bend left=4]{dr}[description,sloped]{\scriptstyle \top} 
\& \& \Labs\CCc \arrow[bend right = 10,dotted]{ll}[swap,pos=.45]{\scriptstyle\liminff \ \cong\  \limsupp} 
\arrow[bend right=13,tail]{ul}[swap]{\scriptstyle\nLan\Cod}
\arrow[bend right=13,tail]{dl}[swap]{\scriptstyle\nRan\Dom}
\\
\& \simUp 
\arrow[bend right=10,two heads]{ur}[swap]{\scriptstyle E}
\arrow[phantom]{ur}[sloped,description]{\scriptstyle\bot}
\arrow[shift left = .5ex,bend right = 15,dotted]{ul}[swap]{\scriptstyle\liminff}
\end{tikzar}
\end{gather}
The functors $H$ and $E$ are induced by \eqref{eq:ida}. The functors $\mnd$ and $\cmn$ map $x\in \CCc$ to $\left<\mnd x, \id_{\cmn x}\right>$ and $\left<\cmn x, \id_{\mnd x}\right>$, respectively, and induce
\bea
\between x & = & \begin{pmatrix}
\mnd x & \cmn x\oot{\id}\cmn x\\
\cmn x & \mnd x \tto\id \mnd x\end{pmatrix}
\eea
Putting together the results from Sections ~\ref{Sec:tightdiag}--\ref{Sec:repres-gen}, the tight completeness claims follow.

%

\begin{proposition}
The category $\UD\CCc$ is 
\begin{enumerate}[a)]
\item complete under $\limsupp$ and $\liminff$, and
\item generated by $\limsupp$ and $\liminff$ along $\between\colon\CCc\to \UD\CCc$.
\end{enumerate}
\end{proposition}

\bprf{ \textbf{(a)}} $\limsupp\left<D, \varphi \right>  =   H\left(\limsupp_{\rgt \varphi} \lft D\right)$\ \  and\ \   $\liminff\left<D, \kappa \right>  =    E\left(\liminff_{\lft \kappa}\rgt D\right)$
\paragraph{(b)} 
$\ida \ \  = \ \  H\left(\limsupp_{\between\rgt \ida}\between \lft A\right)\ \  
 = \ \  E\left(\liminff_{\between \lft a} \between\rgt A\right)$
\epr

\begin{theorem}
The embedding $\between\colon \CCc\to \UD\CCc$ preserves any $\limsupp$ and $\liminff$ that may exist in $\CCc$. It is moreover universal among functors from $\CCc$ into $\limsupp$-and-$\liminff$-complete categories. The functor $\UD:\CAT\to\CAT$ is an idempotent pseudomonad whose unit is an equivalence $\CCc\stackrel{\between}\simeq \UD\CCc$ if and only if $\CCc$ is $\limsupp$-and-$\liminff$-complete.
\end{theorem}

\section{Application: Some categorical tights}\label{Sec:real}

\subsection{Tight completions of groups}
Consider an arbitrary group $G$ as a category $\GGg$ with a single object, i.e., $|\GGg| =\{o\}$ and $\GGg(o,o) = G$, with the group operation $\mmult$ as the composition and the group unit $\iota$ as the identity. Starting from the loose completions again, we note (as explained, mutatis mutandis, in Appendix~\ref{Appendix:monmon}) that
\begin{itemize}
\item the category $\Do \GGg$ of left \actions\ $G\times X\tto\ast X$ can be viewed as the (Eilenberg-Moore) algebra category $\Set^{(G\times)}$ for the monad $(G\times) :\Set\to \Set$, whereas
\item the category $\Up\GGg$ of right \actions\ $Y\times G\tto ! Y$ can be viewed as the \emph{opposite}\/ of the algebra category $\Set^{(\times G)}$ for the monad $(\times G):\Set \to \Set$, or equivalently as the  (Eilenberg-Moore) coalgebra category $\left(\Set^\op\right)^{(\times G)}$ for the comonad\footnote{Both underlying functors are written $(\times G)$, without the superscript $o$.}
 $(\times G):\Set^\op\to \Set^{\op}$. 
\end{itemize}
The Isbell adjunction $\Kan : \Up\GGg\to \Do\GGg$ can thus be viewed as running between the categories of algebras and colagebras on the left in the following diagram.
\beq\label{eq:group-isbell}
\begin{tikzar}[row sep=3.3ex,column sep=1.5em]
\&\& \Set^{(G\times)}  \arrow[bend right = 13,thin,shift left=1]{ddddrrrr}
\arrow[bend right = 15]{dddd}[swap]{\Lan} \arrow[phantom]{dddd}{\dashv} 
 \&\&\&\&  \Set_{(G\times)} \arrow[loop, out = 135, in = 45, looseness = 4]{}[description]{\RLan} 
\arrow[hookrightarrow]{llll} 
\arrow[bend right = 15]{dddd}[swap]{\dLan} \arrow[phantom]{dddd}{\dashv}
\\
\\
\GGg \arrow{uurr}{\mnd} \arrow{ddrr}[swap]{\cmn} 
\\
\\ 
\&\& \left(\Set^\op\right)^{(\times G)} \arrow[bend right = 13,crossing over,thin]{uuuurrrr}
\arrow[bend right = 15]{uuuu}[swap]{\Ran}  
\&\&\&\& 
\Set^\op_{(\times G)} 
\arrow[hookrightarrow]{llll} 
\arrow[bend right = 15]{uuuu}[swap]{\dRan}  
\arrow[loop, out = -50, in=-130, looseness = 4]{}[description]{\LRan} 
\end{tikzar}
\eeq
The upshot is that the monad $\RLan=\Ran\Lan$ and the comonad $\LRan=\Lan\Ran$  can be restricted to the Kleisli categories displayed on the right. Let us spell this out.

We saw in Sec.~\ref{Sec:Lambek-problem}  and stated in \eqref{eq:cones-cocones} that the left adjoint $\Lan$ maps any left \action\ $\lft X = \left(G\times X\tto\ast X\right)$ to the right \action\ over the cocones out of it, while the right adjoint $\Ran$ maps any right \action\ $\rgt Y=\left(Y\times G\tto ! Y\right)$ to the left \action\ over the cones into it. For a group, these cones and cocones are just the equivariant homomorphisms $h:\lft X\to G$ and $k:\rgt Y\to G$ in \eqref{eq:groupcon} below, landing in the group itself as the only representable left \action\ $\mnd o$ in the first case, and the representable right \action\ $\cmn o = G$ in the second.
\beq\label{eq:groupcon}
\begin{tikzar}{}
G\times X\ar{d}[description]{\ast} \ar{r}{G\times h} \& G\times G\ar{d}[description]{\mmult} \\
X \ar{r}[swap]{h}\& G
\end{tikzar}
\qquad\qquad\qquad
\begin{tikzar}{}
Y\times G\ar{d}[description]{!} \ar{r}{G\times k} \& G\times G\ar{d}[description]{\mmult} \\
Y \ar{r}[swap]{k}\& G
\end{tikzar}
\eeq
The condition $h(a\ast x) = a\mmult h(x)$ implies that every $x\in X$ and $a\in G$ satisfy
\bea\label{eq:gx}
a\ast x = x & \implies & a\mmult h(x) = h(x)
\eea 
Since any group element $h(x)$ has an inverse, $a\mmult h(x) = h(x)$ implies $a=\iota$, and \eqref{eq:gx} therefore implies that an \action\ $\lft X = \left(G\times X\tto\ast X\right)$ can support a homomorphism (cocone) $h:\lft X\to G$ only if $a\ast x = x$ implies $a= \iota$. If $\lft X$ permits $a\ast x=x$ for $a\neq\iota$, then there are no cocones out of $\lft X$. If there are, then the elements of the orbit $Gx = \{a\ast x\in X\ |\  a\in G\}$ must be in a one-to-one correspondence with the elements of $G$. On the other hand, the orbits are by definition the equivalence classes modulo the relation
\bea\label{eq:partition}
x\approx y & \iff & \exists a\in G. \ a\ast x = y
\eea
The set $X$ is thus a disjoint union of the orbits of the action on it. When each orbit comes with a bijection to  $G$, there is thus a bijection 
\bea\label{eq:decomX}
X  & \cong &  G\times X_{o}
\eea 
where $X_{o} = X/\approx$ is the set of orbits. Since each orbit is obviously fixed under the action, any \action\  $\lft X = \left(G\times X\tto\ast X\right)$ such that there is a cocone $h:\lft X\to G$ must be in the form
\bea\label{eq:frealg}
\lft X = \big(G\times G\times X_{o} & \tto{\ \ \ \ } & G\times X_{o} \big)\\
a\ast <b, \xi> & \longmapsto & <a\mmult b, \xi>
\eea
The \actions\ in this form are called free, and they are also precisely the free algebras for the monad $(G\times):\Set\to\Set$. Going back to \eqref{eq:groupcon}, whenever $\lft X$ is free, any homomorphism $h$ to $G$ satisfies
\bear
h\left(a\ast<b,\xi>\right) & = & a\mmult h(b,\xi)
\eear
Hence $h(a,\xi) = a\mmult h(\iota, \xi)$. As there are no other constraints, for any free \action\  $\lft X$ this gives a one-to-one correspondence between the assignments $\hhh\in G^{X_{o}}$ and the homomorphisms $h:\lft X\to G$. Hence
\bea\label{eq:Langroup}
\Lan \lft X & = & \begin{cases}
G^{X_o} \times G \tto ! G^{X_o} & \mbox{if } \lft X = \big(G\times G\times X_{o} \tto{\mmult\times X_o} G\times X_{o} \big) 
\\
\emptyset \times G \tto{\ \ } \emptyset & \mbox{ otherwise}
\end{cases}
\eea
where the action on $G^{X_o}$ is pointwise
\bea\label{eq:action}
G^{X_o} \times G &\tto{\ \  \mmult\ \ } & G^{X_o}\\
\left<\ggg, a\right> & \longmapsto & \ggg\mmult a = \sseq{g_{\xi}\mmult a}_{\xi \in X_{o}}\notag
\eea
If we think of $\ggg\in G^{X_{o}}$ as an $X_{o}$-dimensional vector $\ggg = \sseq{g_{\xi}}_{\xi \in X_{o}}$, this action is scalar multiplication. To show that the algebra $\Lan\lft X$ is also free and that the functor $\Lan:\Set^{(G\times)}\to \Set^{(\times G)}$ factors through the category  of free algebras $\Set_{(\times G)}$ as claimed in  \eqref{eq:group-isbell}, we need to decompose its underlying set $G^{X_{o}}$ of $\Lan\lft X$ in the same way as the underlying set $X$ of $\lft X$ was decomposed in \eqref{eq:decomX}. The task is thus to determine the set of orbits $\left(G^{X_o}\right)_o$ of $\Lan\lft X$ so that the decomposition
\bea
G^{X_o} & \cong & \big(G^{X_o}\big)_o \times G
\eea
reduces \eqref{eq:action} to the form
\bea\label{eq:action-dec}
\big(G^{X_o}\big)_o \times G \times G &\tto{\ \  !\ \ } & \big(G^{X_o}\big)_o \times G\\
\left<\, \gamma\ ,\  b\ ,\ a\,  \right> & \longmapsto & \left<\, \gamma\ ,\  b\mmult a\, \right>\notag
\eea
Aligning \eqref{eq:action} and  \eqref{eq:action-dec} shows that the orbit set must be the quotient 
\bea\label{eq:Goo}
\big(G^{X_o}\big)_o & = & G^{X_{o}}/ \approx
\eea
modulo the equivalence relation defined for $\ggg, \hhh\in G^{X_o}$ by   
\bea\label{eq:scalmult}
\ggg \approx \hhh &\iff & \exists a\in G.\ \ggg\mmult a = \hhh
\eea
Viewing $\ggg$ and $\hhh$ as vectors makes the elements of $\big(G^{X_o}\big)_o$ into rays of colinear vectors, i.e., \emph{projective lines}. 
If $G^{X_{o}}$ is a vector space presented in  the cartesian coordinates $\sseq{g_{\xi}}_{\xi \in X_{o}}$, then $\big(G^{X_o}\big)_o$ is the corresponding \emph{projective}\/ space presented in the \emph{homogenous}\/ coordinates $\psseq{g_{\xi}}_{\xi \in X_{o}}$, satisfying the usual homogeneity property
\bea\label{eq:homogenous}
\psseq{g_{\xi_{0}}\ ,\  g_{\xi_{1}}\ ,\  g_{\xi_{2}}\ ,\  \ldots}  & = & \psseq{g_{\xi_{0}}\mmult a\ \, ,\  g_{\xi_{1}}\mmult a\, \  ,\   g_{\xi_{2}}\mmult a\ \,  ,\  \ldots} 
\eea
\paragraph{Intuition and notation: The orbit sets are the projective spaces.} To simplify notation and perhaps make use of familiarity with homogenous coordinates, we replace orbit sets $\big(G^{X_o}\big)_o$, defined as the quotients in \eqref{eq:Goo}, by ``projective spaces''
\bea\label{eq:Xcircledast}
X^{\circledast} & = & \Big\{ \psseq{\ggg}\ |\ \ggg \in G^{X_{o}}
\Big\}
\eea
where $\psseq\ggg = \psseq{g_{\xi}}_{\xi\in X_{o}}$ satisfy \eqref{eq:homogenous}. The difference between $X^{\circledast}$ and $\big(G^{X_o}\big)_o$ is, of course, mainly cosmetic, since homogenous coordinates are just a notation for equivalence classes of vectors modulo scalar multiplication. But the notational cosmetics will help us iterate.
%

The analysis of the right \actions\ $\rgt Y= \left(Y\times G\tto ! Y\right)$ is analogous to the left \actions\ and $\Ran \rgt Y$ is in the form $G\times G^{Y_o} \tto \ast G^{Y_o}$ whenever $\rgt Y$ is free and empty otherwise. The $\Ran$-images are  thus free $(G\times )$-algebras in any case, generated by the orbit set $\big(G^{Y_o}\big)_o$ again, which can be presented as the projective space $Y^{\circledast}$ of rays of homogenous  coordinates again. We have thus shown that the functors $\Lan$ and $\Ran$ factor through $\Set^{\op}_{(\times G)}$ and $\Set^{\op}_{(G\times)}$, as claimed in \eqref{eq:group-isbell}. 
%
%
%
%
%

A further claim is that the algebras for the monad $\RLan = \Ran\Lan$ can only be supported by the free $(G\times)$-algebras, and that the coalgebras for the comonad $\LRan = \Lan\Ran$ can only be supported by the cofree $(\times G)$-coalgebras. The reason is that both $\Lan$ and $\Ran$ take the \actions\ that are not free to the empty set, and the algebra structure maps must be surjective. The constructions involving the $\RLan$-algebras can thus be restricted from $\Set^{(G\times)}$ to $\Set_{(G\times)}$ without any loss; and the constructions of the $\LRan$-coalgebras can be restricted from $\left(\Set^{\op}\right)^{(\times G)}$ to $\Set^{\op}_{(\times G)}$. In particular, on the path to the tight completion of the group $G$, the nucleus of the Isbell adjunction $\Kan:\Up\GGg\to \Do\GGg$ can be reduced to the nucleus of the adjunction beween the free algebras $\dKan: \Set^{\op}_{(\times G)}\to \Set_{(G\times)}$. The free algebras are conveniently presented in Kleisli form:
\begin{gather}\label{eq:Kleisli-times}
\Set_{(G\times)}\ \  =\ \  |\Set|\ \ = \ \ \Set_{(\times G)}\\
\Set_{(G\times)}(A,B) \ \ =\ \ \Set(A,G\times B) \qquad \qquad 
\Set_{(\times G)}(A,B) \ \ =\ \ \Set(A,B\times G) \notag
\end{gather}
spelled out in Appendix~\ref{Appendix:monad}. Kleisli composition \eqref{eq:kleislicomp} is in this case
\bea
\left(A\tto{<\varphi, f>}G\times B\right)\boxdot \left(B\tto{<\psi, g>}G\times C\right) & = & \left(A\tto{\left<\varphi\mmult\left(f\bullet \psi\right)\,,\  f\bullet g\right>} G\times C\right)
\eea
The objects $A, B, C$ of $\Set_{(G\times)}$ and $\Set_{(\times G)}$ represent the generator sets of the free algebras. In the Eilenberg-Moore view in $\Set^{(G\times)}$ and $\Set^{(\times G)}$, the generator sets of the free algebras were the sets of orbits $X_o, Y_o$. Now they are the first class citizens. Restricted to free algebras and coalgebras and translated to Kleisli form, the adjoint functors $\Kan$ become
\[\begin{tikzar}[row sep = 2ex]
\dLan\colon \Set_{(G\times)} \ar{r}\& \Set^{\op}_{(\times G)}
\\[-1ex]
\ \ \ \ A \ar[mapsto]{r} \ar[thin,shift left=1ex]{dd}[description]{<\alpha,f>} 
\& A^{\circledast}
\\  
\& A^{\circledast}\times G
\\
\ \ \ \ G\times B
\\
\ \ \ \ B \arrow[mapsto]{r}  \& B^{\circledast} \arrow[thin]{uu}[description]{\left<f^{\circledast},\alpha^{\circledast}\right>}  
\end{tikzar}
\qquad\qquad\qquad\qquad
\begin{tikzar}[row sep = 2ex]
\Set_{(G\times)} \& \Set^{\op}_{(\times G)}\ :\dRan \ar{l}
\\[-1ex]
A^{\circledast}  \ar[thin]{dd}[description]{\left<\beta^{\circledast},t^{\circledast}\right>} 
\& A \ar[mapsto]{l} \ \ \ \ 
\\  
\& A\times G \ \ \ \ 
\\
G\times B^{\circledast}
\\
B^{\circledast}   \& B\ \ \ \  \arrow[mapsto]{l} \arrow[thin,shift left=1ex]{uu}[description]{\left<t,\beta\right>}  
\end{tikzar}
\]
where $A^{\circledast}$ and $B^{\circledast}$ are the ``projective spaces'' \eqref{eq:Xcircledast} of homogenous $A$-tuples and $B$-tuples from $G$. The "projecitivizations" define the object parts of both functors. The arrow parts\footnote{The superscript $\circledast$ in the components is just a convenient way to reuse the names of the input components and does not refer to an operation.} $\dLan(\alpha, f) = \left<f^{\circledast}, \alpha^{\circledast}\right>$ and $\dRan(t,\beta) = \left<\beta^{\circledast}, t^{\circledast}\right>$ are
\begin{align*}
f^{\circledast}\pseq{k_{y}}_{y\in B} & =  \pseq{\alpha_{x}\mmult k_{f(x)}}_{x\in A} & \alpha^{\circledast}\pseq{k_{y}}_{y\in B} & =  \iota \\ 
\beta^{\circledast}\pseq{h_{x}}_{x\in A} & =  \iota  & t^{\circledast}\pseq{h_{x}}_{x\in A} & =  \pseq{h_{t(y)}\mmult \beta_{y}}_{y\in B} 
\end{align*}
The morphisms $\dLan(\alpha,f)$ and $\dRan(t,\beta)$ remain unchanged if any other fixed elements of $G$ are taken to be $\alpha^\circledast$ and $\beta^\circledast$. This is because $f^\circledast$ and $t^\circledast$ are invariant under scalar multiplication. The monad and the comonad are of the form
\[\begin{tikzar}[row sep = 2ex]
\RLan\colon \Set_{(G\times)} \ar{r}\& \Set_{(G\times)}
\\[-1ex]
\ \ \ \ A \ar[mapsto]{r} \ar[thin,shift left=1ex]{dd}[description]{<\alpha,f>} 
\& A^{\circledast\circledast} \arrow[thin]{dd}[description]{\left<\alpha^{\circledast\circledast},f^{\circledast\circledast}\right>}  
\\[3ex]  
\\
\ \ \ \ G\times B \& G\times B^{\circledast\circledast}
\end{tikzar}
\qquad\qquad\qquad\qquad
\begin{tikzar}[row sep = 2ex]
\LRan\colon \Set^{\op}_{(\times G)} \ar{r}\& \Set^{\op}_{(\times G)}
\\[-1ex]
\ \ \ \ A \times G 
\& A^{\circledast\circledast} \times G  
\\[3ex]  
\\
\ \ \ \ B \ar[mapsto]{r}  \ar[thin,shift right=1ex]{uu}[description]{<t,\beta>}   \& B^{\circledast\circledast} \arrow[thin]{uu}[description]{\left<t^{\circledast\circledast}, \beta^{\circledast\circledast}\right>} 
\end{tikzar}
\]
with the components of $\RLan(\alpha, f)$ and $\LRan(t,\beta)>$ in the form
\begin{align*}
\alpha^{\circledast\circledast}\left[\Phi_{\left[\ggg\right]}\right]_{\left[\ggg\right]\in A^\circledast} & =  \iota    
&
f^{\circledast\circledast}\left[\Phi_{\left[\ggg\right]}\right]_{\left[\ggg\right]\in A^\circledast} & =  \left[\Phi_{\Lan(\alpha, f)[\hhh]}\right]_{\left[\hhh\right]\in B^\circledast}
\\ 
t^{\circledast\circledast}\left[\Psi_{\left[\hhh\right]}\right]_{\left[\hhh\right]\in B^\circledast} & =  \left[\Psi_{\Ran(t,\beta)[\ggg]}\right]_{\left[\ggg\right]\in A^\circledast} &
 \beta^{\circledast\circledast}\left[\Psi_{\left[\hhh\right]}\right]_{[\hhh]\in B^\circledast} & =  \iota 
\end{align*}
The $\RLan$-algebras and the $\LRan$-coalgebras provide retractions of the form $A^{\circledast\circledast}\to A$ for the "projective spaces" over a group. However, manipulating the Eilenberg-Moore algebras for  $\RLan$ and $\LRan$ within the Kleisli categories for $(G\times)$ and $(\times G)$ gets quite involved. An explicit presentation of full cuts for a group seems out of reach. Simple and absolute cuts provide a way forward. 
A tight completion $\Labs G$ of a group $G$, along the lines of Def.~\ref{Def:cutabsol}, comprises tuples $<A,B,\ida,\idb,\varphi,\psi>$ as objects, where $A$ and $B$ are sets and $\ida:B^{\circledast}\to G\times B^{\circledast}$ and $\idb:A^{\circledast}\to A^{\circledast} \times G$ are idempotents in $\Set_{(G\times )}$ and $\Set_{(\times G)}$ whose images split on each other:
\bea\label{eq:ida-group}
\begin{tikzar}[column sep=1.8em,row sep = 7ex]
B^{\circledast} \ar[bend right = 15,shift right=2]{dd}[description]{\ltimes}[swap,pos=0.25]{\ida} \ar[two heads]{d}[description]{\ltimes} \& B^{\circledast\circledast} \&\& B^{\circledast\circledast}
\ar{ll}[description]{\rtimes}[swap,pos=0.2]{\psi}
\\
A \ar[tail]{d}[description]{\ltimes} \& A^{\circledast}  \ar[tail]{u}[description]{\rtimes} \& \& A^{\circledast} \ar{ll}[description]{\rtimes}[swap,pos=0.2]{\idb} 
\ar[tail]{u}[description]{\rtimes}  
\\
B^{\circledast}   \& B^{\circledast\circledast} \ar[two heads]{u}[description]{\rtimes} \&\& B^{\circledast\circledast} \ar{ll}[description]{\rtimes}[swap,pos=0.2]{\psi} \ar[two heads]{u}[description]{\rtimes}
 \end{tikzar}
& \qquad &
\begin{tikzar}[column sep=1.8em,row sep = 7ex]
A^{\circledast}   \& A^{\circledast\circledast} \ar{rr}[description]{\ltimes}[pos=0.2]{\varphi} \ar[two heads]{d}[description]{\ltimes} \&\& A^{\circledast\circledast}   \ar[two heads]{d}[description]{\ltimes}
\\
B \ar[tail]{u}[description]{\rtimes}  \&B^{\circledast} \ar{rr}[description]{\ltimes}[pos=0.2]{\ida}  \ar[tail]{d}[description]{\ltimes} \& \& B^{\circledast}  
  \ar[tail]{d}[description]{\ltimes}
\\
A^{\circledast} \ar[bend left = 15,shift left=2]{uu}[description]{\rtimes}[pos=.25]{\idb} \ar[two heads]{u}[description]{\rtimes}  \& A^{\circledast\circledast}\ar{rr}[description]{\ltimes}[pos=0.2]{\varphi}  \&\& A^{\circledast\circledast}  
 \end{tikzar}
\eea
The arrows with the marking $\ltimes$ are in $\Set_{(G\times)}$, whereas the arrows with $\rtimes$ are in $\Set_{(\times G)}$. This means that each marked arrow comes with a $G$-component, as specified in \eqref{eq:Kleisli-times}. The arrows that are left unnamed are the gaps $\gap$ and the intervals $\intv$ from the simple cuts obtained by splitting the displayed idempotents. Switching from the displayed absolute cuts to the underlying simple cut presentation, by erasing the idempotents and naming the gaps and the intervals, displays how $A$ and $B$ determine each other as retracts\footnote{While simple cuts are simpler, we present it here in terms of absolute cuts mainly because diagrams like \eqref{eq:ida} and \eqref{eq:ida-group} display both views, the idempotents and their splittings, which is sometimes helpful.}. The homomorphisms are the same in both cases. A $\Labs G$-morphism from $<A,B,\ida,\idb,\varphi,\psi>$ to $<C,D,\idc,\idd,\xi,\zeta>$ is a pair $<f,g>$ of functions $f:A\to G\times C$ and $g: D\to B\times G$ whose images preserve the idempotents:
\[
\begin{tikzar}[row sep = 10ex,column sep = 10ex]
B^{\circledast} \ar{d}[description]{\ltimes}[swap,pos=0.3]{\Ran g}\ar{r}[description]{\ltimes}[pos=0.3]{\ida} \& B^{\circledast} \ar{d}[description]{\ltimes}[pos=0.3]{\Ran g}
\\
D^{\circledast} \ar{r}[description]{\ltimes}[swap,pos=0.3]{\idc} \& D^{\circledast}
\end{tikzar}\qquad \qquad\qquad\qquad 
\begin{tikzar}[row sep = 10ex,column sep = 10ex]
A^{\circledast} \ar[leftarrow]{d}[description]{\rtimes}[swap,pos=0.7]{\Lan f}\ar[leftarrow]{r}[description]{\rtimes}[pos=0.7]{\idb} \& A^{\circledast}
\ar[leftarrow]{d}[description]{\rtimes}[pos=0.7]{\Lan f}
\\
C^{\circledast} \ar[leftarrow]{r}[description]{\rtimes}[swap,pos=0.7]{\idd} \& C^{\circledast}
\end{tikzar}
\]
Since the idempotents split on each other, $f$ and $g$ determine each other, and specifying either of them suffices.

\paragraph{The tight completion of $\ZZz_{4}$} comprises pairs of sets $A, B$ where $A$ is a retract of $B^{\circledast}$ and $B$ is a retract of $A^{\circledast}$, as displayed in \eqref{eq:ida-group}. When $A$ and $B$ are finite sets, say with $m$ and $n$ elements respectively, then $A^{\circledast}\cong \ZZz_{4}^{m-1}$ and $B^{\circledast}\cong \ZZz_{4}^{n-1}$ means that each constrains the other by $m\leq 4^{n-1}$ and $n\leq 4^{m-1}$. The tight completion of $\ZZz_{4}$ thus comprises retracts of its powers.

\paragraph{The tight completion of $\QQq$} is similarly comprised of pairs of retractions $A$ of $B^{\circledast}$ and $B$ of $A^{\circledast}$, where $A^{\circledast}$ and $B^{\circledast}$ are familiar either as the 1-dimensional Grasmanians of the spaces $\QQq^{A}$ and $\QQq^{B}$, or as the rational subspaces of real projective spaces of dimensions $A$ and $B$. They are also the projective objects in the category of $\QQq$'s \actions\ on its powers. Projective spaces gave name to projective modules, which gave name to the abstract notion of projectivity in categories, which in tight completions refers back to projective spaces.

\subsection{Tight completions of monoids and categories}
Viewing a monoid $(M, \mmult, \iota)$ as a category $\MMm$ with a single object, i.e., $|\MMm| = \{o\}$ and $\MMm(o,o)=M$ leads to the loose completions $\Do\MMm$ and $\Up\MMm$. The correspondence of $\MMm$-\actions\ and the algebras for the induced monads, described in Appendix~\ref{Appendix:monmon}, still admits interpreting  $\Do\MMm$ as the category of algebras $\Set^{(M\times)}$, and the $\Up\MMm$ as the category of coalgebras $\left(\Set^\op\right)^{(\times M)}$. The adjunction as in \eqref{eq:group-isbell} on the left is again realized by homming into the sole representable, like in \eqref{eq:groupcon}, with $G$ replaced by $M$. A cocone $h:\lft X\to M$ still satisfies $h(a\ast x) = a\mmult h(x)$, but this is as far as replacing groups by monoids goes. Since the group elements have the inverses, homming into a group induced a free \action, and the adjunction in \eqref{eq:group-isbell} factored through free algebras and coalgebras. Homming into a monoid maps the orbits $Mx = \{a\ast x\in X\ |\  a\in M\}$ into monoid ideals, and the orbits are generally not as big as $M$, and may not partition $X$. The relation
\bea\label{eq:preord}
x\prec y & \iff & \exists a\in M. \ a\ast x = y\ \ \iff\ \ Mx\subseteq My
\eea
is now a mere preorder.

\begin{lemma} \label{lemma:freemon} An orbit $Mx$ is in a one-to-one correspondence with $M$ if and only if 
\bea\label{eq:cancellable}
 a\ast x=b\ast x &\implies & a=b
\eea
The preorder $\prec$ is symmetric if and only if every pair that has a lower bound also has an upper bound:
\bea\label{eq:gcd}
z\prec x, y & \implies & \exists u.\ x,y \prec u
\eea
\end{lemma}
The monoid \action\ extensions $\Kan: \left(\Set^\op\right)^{(\times M)} \tto{\ \ \ \ }\Set^{(M\times)}$ are thus in the form \eqref{eq:Langroup} and factor through free algebras just when the monoid $M$ itself satisfies (\ref{eq:cancellable}--\ref{eq:gcd}). In general, the induced algebras are still projective, but not free:
\bea\label{eq:Lanmonoid}
\Lan \lft X & = & \left(X^\circledast \times M \tto ! X^\circledast\right) \mbox{ where}\\
&& X^\circledast \ =\ \left\{ \hhh\in M^X\ |\ \hhh_{a\ast x} = a\mmult \hhh_x\right\}
\eea
with the pointwise action $\hhh \mmult b = \sseq{ h_x\mmult b}_{x\in X}$. The right extension $\Ran\rgt Y$ is analogous. The $\RLan$-algebras and the $\LRan$-coalgebras for $\RLan=\Ran\Lan$ and $\LRan=\Lan\Ran$ in the standard (Eilenberg-Moore) form carry a lot of structure, but using the simple nucleus presentation from \cite[Sec.~8]{PavlovicD:nucleus} simplifies the task. This construction is, however, nearly as general as it gets, as noted in Appendix~\ref{Appendix:catspan}.

\bibliographystyle{plain}
\bibliography{ZZZZ-ref-ded,CT,PavlovicD,math,logic}

\appendix
\section{Appendix: Categories and \actions}\label{Appendix:basic}

While it is unlikely that a reader completely uninitiated in categories will be able to follow all  \emph{technical}\/ details of tight completions by reading any single paper at this moment, the summary of the basic categorical concepts in this appendix will hopefully enable a sufficiently curious reader to follow the \emph{conceptual}\/ developments in the present paper.

Even the experienced reader may find the following alignments useful, since categories are many things to many people, and some of our developments have been simplified by nonstandard approaches and context-dependent notation.

\subsection{Categories}\label{Appendix:category}
\subsubsection{Monoids, preorders}\label{Appendix:mon-preord}
Categories are a common generalization of preorders and monoids. The common denominator of preorders and monoids is the pattern echoed in the preorder axioms and monoid operations:
\begin{align*}
\top &\implies x \leq x & x\leq y \wedge y \leq z &\implies x\leq z\\
1&\tto{\ \ \iota\ \  }\ \  M & M\ \ \times\ \  M\ \ \,  & \tto{\ \ \mmult\ \ } M
\end{align*}
The monoid operations are required to be associative and unitary, which means that the following diagram commutes
\beq\label{eq:monoid}
\begin{tikzar}[row sep = 6ex,column sep = 1.8em]
M\times M\times M\ar{rr}{M\times \mmult} \ar{dd}[swap]{\mmult \times M}
\&\& M\times M\ar{dd}{\mmult}
\\
\& M\ \ \ar[equals]{dr}\ar{ur}[description]{M\times \iota} \ar{dl}[description]{\iota \times M} \\
M\times M\ar{rr}[swap]{\mmult} \&\& M
\end{tikzar}
\eeq
A (small) category $\CCcc$ comprises
\begin{itemize}
\item a set of \emph{objects}\/ $\lvert\CCcc\rvert$,
\item for any pair\footnote{We often write $x,y\in \CCcc$ instead of $x,y\in |\CCcc|$.} $x, y \in |\CCcc|$ a \emph{hom-set} of morphisms $\CCcc(x,y)$, and
\item for all $x, y, z \in |\CCcc|$ the operations 
\begin{align}\label{eq:catsig}
1&\tto{\ \ \iota_{x}\ \  } \CCcc(x,x) & \CCcc(x,y) \times \CCcc(y,z)  & \tto{\ \ \bullet\ \ }  \CCcc(x,z)
\end{align} 
\end{itemize}
making the following diagram commute
{\footnotesize
\beq\label{eq:category}
\begin{split}
\begin{tikzar}[row sep = 8ex,column sep = 1.3em]
\CCcc(x,y)\times \CCcc(y,z)\times \CCcc(z,a)\ar{rr}{\id \times \bullet} \ar{dd}[swap]{\bullet \times \id}
\&\& \CCcc(x,y)\times \CCcc(y,a)\ar{dd}{\bullet}
\\
\mbox{\hspace{4em}}
\\
\CCcc(x,z)\times \CCcc(z,a)\ar{rr}[swap]{\bullet} \&\& \CCcc(x,a)
\end{tikzar}
\quad 
\begin{tikzar}[row sep = 7ex,column sep = .4em]
\&\& \CCcc(x,a)\times \CCcc(a,a)\ar{dd}{\bullet}
\\
\& \CCcc(x,a)\ar{ur}[description]{\id\times \iota_{a}} \ar{dl}[description]{\iota_{x} \times \id} \ar[equals]{dr} 
\\
\CCcc(x,x)\times \CCcc(x,a)\ar{rr}[swap]{\bullet} \&\& \CCcc(x,a)
\end{tikzar}
\end{split}
\eeq}%
A monoid is thus a category with one object, while a preorder is a category with at most one morphism in each hom-set $\CCcc(x,y)$, and $x\leq y$ that there is a morphism from $x$ to $y$. A category is thus a family $\CCcc:\lvert\CCcc\rvert \times \lvert\CCcc\rvert \to \Set$ of hom-sets indexed over the objects, and a preorder is such a family of sets with at most one element, i.e., a relation $(\leq):\lvert\CCcc\rvert \times \lvert\CCcc\rvert \to \{0,1\}$. For a preorder, all diagrams like \eqref{eq:category} obviously commute.

\paragraph{The composition notations.} When morphisms are drawn and composed left to right, we write composition as above, in  "geometric" order:
\bear
\left(x\tto f y \tto g z\right) & = & \left(x\tto{f\bullet g} z\right)
\eear
In algebra, however, the arguments of functions are written on the right of the function symbols, we write composition in "algebraic" order:
\bea\label{eq:comp}
g\left(f(a)\right) & = & \left(g\circ f\right)(a)
\eea
The sequential composition in \eqref{eq:catsig} would be written in the form
\bea\label{eq:catsig-rl}
\CCcc(y,z) \times \CCcc(x,y) &\tto{\ \ \circ\ \ } & \CCcc(x,z)
\eea
and the diagrams in \eqref{eq:category} would need to be adjusted. While most people get used to \eqref{eq:comp} early on, prefer \eqref{eq:catsig-rl} to \eqref{eq:catsig}, and adjust \eqref{eq:category} accordingly, the notations for \actions, starting from \eqref{eq:cataction}, become cumbersome and error-prone. For simplicity, we keep in use both \eqref{eq:catsig-rl} and \eqref{eq:catsig}.

\subsubsection{Categories are monoids of spans}\label{Appendix:catspan}
If a monoid is summarized as a diagram in the form
\bea\label{eq:monspan}
\MMm & = & \left(\begin{tikzar}[row sep = 8ex,column sep = 7ex] 1 \ar{r}{\iota} \& M \& M\times M \ar{l}[swap]{\mmult} \end{tikzar} \right)
\eea
satisfying \eqref{eq:monoid}, then a category can be summarized as a diagram in the form
\bea\label{eq:catspan}
\CCcc & = & \left(\begin{tikzar}[row sep = 8ex,column sep = 8ex] C_0 \ar{r}{\iota} \ar{dr}[description,pos=0.35]{\Delta} \& C_1 \ar{d}[description,pos=0.35]{<\delta_0,\delta_1>}\& C_2 \ar{l}[swap]{\bullet} \ar{dl}[description,pos=0.35]{<d_0,d_2>}\\
\& C_0\times C_0\end{tikzar} \right)
\eea
where $C_0 = \lvert\CCcc\rvert$ is the set of objects, $\Delta(x) = <x,x>$ is the diagonal, $C_1 = \coprod_{x,y\in C_0} \CCcc(x,y)$ is the set of morphisms with the projections $\delta_0, \delta_1 : C_1\to C_0$ mapping each  $x\to y$ to $x$ and $y$ respectively, and $C_2$ is obtained as the pullback
\beq\label{eq:catspan-gen}
\begin{tikzar}[column sep = 1cm,row sep = 1cm]
\&\&C_2 \ar[bend right=40]{ddll}[swap]{d_0} \ar{dr}[description]{p_1} \ar{dl}[description]{p_0} \ar[bend left=40]{ddrr}{d_2} \ar[phantom]{dd}[description,pos=0.1]{\pbdown}
\\
\& C_1 \ar{dl}[description]{\delta_0}  \ar{dr}[description]{\delta_1} \&\& C_1 \ar{dl}[description]{\delta_0}  \ar{dr}[description]{\delta_1}
\\
C_0 \&\& C_0 \&\& C_0
\end{tikzar}
\eeq 
%
The set $C_2$ thus consists of all pairs of morphisms in the form $<x\to y, y\to z> \in C_1\times C_1$. The projections $d_0, d_2:C_2\to C_0$, induced in \eqref{eq:catspan-gen}, map each pair $<x\to y, y\to z>$ to $x$ and $z$, whereas the diagonal of the pullback $d_1:C_2\to C_0$ maps it to $y$. Note that the set $C_2$ is often written in the form $C_1\underset{C_0}{\times} C_1$, which is an imprecise but intuitive reminder that the pullback of $\delta_1$ and $\delta_0$ is their product in the slice category\footnote{See Appendix~\ref{Appendix:comma}.} $\comm{\Set}{C_0}$. When $C_0 = 1$, diagram \eqref{eq:catspan} boils down to \eqref{eq:monspan}, which just says that monoids are categories with one object. The other way around, categories can also be viewed as indexed monoids. This was already suggested by diagram \eqref{eq:category} as an indexed version of \eqref{eq:monoid}. But now a much bigger picture emerges.  Any monoid $M$ can be viewed as a diagram \eqref{eq:monspan} in the category $\Set$ of sets. Any small category $\CCcc$ can be viewed as monoid in the category $\comm{\Set}{C_0\times C_0}$ over the set of objects $C_0$ \cite[Sec.~5.4.3]{BenabouJ:bicats}. 

Note, however, that the categories $\Set$ and $\comm{\Set}{C_0\times C_0}$ are not small.

\subsubsection{Large categories}\label{Appendix:size}
Taking 
\begin{itemize}
\item sets as the objects $A, B,\ldots \in |\Set|$ and 
\item functions $A\to B$ as the morphisms $f,g,\ldots \in \Set(A,B)$
\end{itemize}
yields the category $\Set$ of sets. It carries all structure described in Sec.~\ref{Appendix:category}, but it does not satisfy the requirement that the objects form a set. Assuming that the collection of all sets is a set leads to logical paradoxes. The need to structure the universe of mathematical discourse to avoid the paradoxes propelled the development of set theory \cite{vanHeijenoortJ}. The simplest solution, due to von Neumann, Bernays, and G\"odel, was to distinguish (small) sets, comprehended by predicates over other sets, from (large) classes, definable by the predicates alone \cite{BernaysP:axiomatic,KanamoriA:bernays}. This solution was adopted in categorical practice early on, sometimes conveniently but conservatively iterated \cite[Expos\'e i, Appendice: Univers]{GrothendieckA:SGA4}. The issue of the size of categories also inspired interesting categorical treatments  \cite{StreetR:cosmoi,StreetR:size}.  
For the present work, it should only be borne in mind that
\begin{itemize}
\item the completions that we study are under \emph{small}\/ limits and colimits: our diagrams are always indexed by small categories, but
\item a completion under \emph{all}\/ small limits and colimits cannot in general be small itself: completing under operations indexed over all sets generally leads to \emph{large}\/ categories.
\end{itemize}
%
%
%
%
The issue is settled by distinguishing two universes: the universe $\Set$ of (small) sets and the universe of (large) classes, to which we do not assign a name, since we will never need to refer to it explicitly again. 

%

{\footnotesize 
\paragraph{\footnotesize Full disclosure: small = finite.} Our own thinking, motivated by the applications in data analysis, is to think of small diagrams as finite, and of large categories as infinite. Our experience from \cite{PavlovicD:moncom} suggests that the infinitary constructions of large categories are at least to some extent an artifact of uncontrolled abstraction. Most of the large categories one encounters are specified by small descriptions. For example, the category of groups is specified by two operations and a few equations, over the category of sets, which is specified by a couple of rules for generating and manipulating sets. Some theories require infinite sets of axioms, but these infinite sets have finite specifications themselves. The whole process resembles programming, where infinite computable functions and processes arise from finite programs \cite[computable categories chapter]{PavlovicD:moncom}. This is not a peculiarity of programming, but the essence of language, described by Wilhelm von Humboldt as \emph{"the infinite use of finite means"}. The ongoing stream of conversations and texts, human and mechanical, is generated from finite dictionaries by finite sets of rules. We need just two words to describe the "Infinite Universe", and just one word to describe "Everything". We sometimes use transfinite constructions to build large categories, but we use finite descriptions to specify the transfinite constructions. That is why the problems arising from size distinctions might disappear if we do not abstract away the way in which we abstract away the towers of abstractions that we build to describe the descriptions of categories. }

\subsubsection{Functors, natural transformations}\label{Appendix:functor}
The fact that categories can be studied within categories propelled categories into a foundational role \cite{BenabouJ:naive,LawvereFW:sets,LawvereFW:Foundation,hottbook}. Categories form a category $\Cat$, with the small categories as objects and functors as morphisms. As the definition of small categories in Sec.~\ref{Appendix:category} suggests, a functor $F:\AAa\to \BBb$, as a morphism from a category $\AAa$ to a category $\BBb$, should consist of an object-function $F\colon\lvert\AAa\rvert\to\lvert\BBb\rvert$ and a family of morphism-functions $F_{xy}\colon \AAa(x,y)\to \BBb(Fx,Fy)$ such that $F_{xx}{\iota_x} = \iota_{Fx}$ and $F_{xy}f \bullet F_{yz}g = F_{xz}(f\bullet g)$  for all $x,y, z\in\lvert \AAa\rvert$. More precisely, $\Cat$ is not just a category but also a \emph{2-category}, since there are also morphisms between functors, called \emph{natural transformations}. In the seminal paper where functors and categories were introduced as mathematical structures, natural transformations were the main purpose \cite{Eilenberg-MacLane:natural}. A natural transformation $\varphi:F\to G$, where $F,G:\AAa\to\BBb$ are functors, is a family of morphisms $\varphi x: Fx\to Gx$ such that $\varphi x\bullet G f = Ff \bullet \varphi y$ holds for all $f\in \AAa(x,y)$. The theory of categories, functors, natural transformations, and higher categories is described in textbooks and handbooks \cite{AdamekJ:CT,AwodeyS:CT,BorceuxF:handbook,LeinsterT:operads,MacLaneS:CWM}.

\subsection{\Actions} \label{Appendix:action}
A vector space can be defined as a multiplicative monoid of a field acting on an   abelian group. A module is a multiplicative monoid of a ring acting on an abelian group. But monoid \actions\ do not necessitate the abelian group structure, and can also be studied over sets irrespective of any underlying structure.   
%
%
%
 A monoid $\MMm$ with operation $\mmult$ acts on a set $X$ on the left or on a set $Y$ on the right along the functions in the form
\begin{align*}
M\times X&\tto{\ \ \ast\ \  }\ \  X & Y\times M & \tto{\ \ !\ \ } Y
\end{align*}
satifying the equations that make the following diagrams commute
\beq\label{eq:monaction}
\begin{tikzar}[row sep = 6ex,column sep = 1.8em]
M\times M\times X\ar{rr}{M\times \ast} \ar{dd}[swap]{\mmult \times X}
\&\& M\times X\ar{dd}{\ast}
\\
\& X\ \ \ar[equals]{dr} \ar{dl}[description]{\iota \times X} \\
M\times X\ar{rr}[swap]{\ast} \&\& X
\end{tikzar}
\qquad\qquad
\begin{tikzar}[row sep = 6ex,column sep = 1.8em]
Y\times M\times M\ar{rr}{Y\times \mmult} \ar{dd}[swap]{\ast \times M}
\&\& Y\times M\ar{dd}{!}
\\
\& Y\ \ \ar[equals]{dr}\ar{ur}[description]{Y\times \iota} \\
Y\times M\ar{rr}[swap]{!} \&\& Y
\end{tikzar}
\eeq
Diagram \eqref{eq:monoid} thus says that every monoid acts on itself. When the monoid $M$ admits all inverses, i.e., when it is a group, the monoid \actions\ are the group actions, a powerful mathematical tool going back to Cayley and Klein. When the underlying set of the monoid is a group and the monoid operation is a group homomorphism, the monoid is a ring, it acts on groups, and its \actions\ are ring modules. When the ring is a field, its actions are vector spaces. 

Stepping up to a category $\CCcc$, left and the right \actions\ become families of functions
\begin{align*}
\CCcc(x,y)\times A{y}&\tto{\ \ \ast\ \  }\ \  A{x} & B{x}\times \CCcc(x,y) & \tto{\ \ !\ \ } B{y}
\end{align*}
where $A, B:|\CCcc| \to \Set$ are families of sets indexed over the objects of $\CCcc$. This time the \actions\ are required to satisfy
{\small
\beq\label{eq:cataction}
\begin{split}
\begin{tikzar}[row sep = 7ex,column sep = 1.2em]
\CCcc(x,y)\times \CCcc(y,z)\times A{z}\ar{rr}{\id\times \ast} \ar{dd}[swap]{\bullet \times \id}
\&\& \CCcc(x,y)\times A{y}\ar{dd}{\ast}
\\
\mbox{\hspace{4em}}
\\
\CCcc(x,z)\times A{z}\ar{rr}[swap]{\ast} \&\& A{x}
\end{tikzar}
\\[3ex]
\begin{tikzar}[row sep = 9ex,column sep = 2.2em]
\&\& 
\\
\& Ax
\ar{dl}[description]{\iota_{x} \times \id} \ar[equals]{dr} 
\\
\CCcc(x,x)\times Ax\ar{rr}[swap]{\ast} \&\& Ax
\end{tikzar}
\end{split}
\qquad
\begin{split}
\begin{tikzar}[row sep = 7ex,column sep = 1.2em]
B{x}\times \CCcc(x,y)\times \CCcc(y,z) \ar{rr}{Bx\times \bullet} \ar{dd}[swap]{! \times \id}
\&\& Bx\times \CCcc(x,z)\ar{dd}{!}
\\
\mbox{\hspace{4em}}
\\
By\times \CCcc(y,z) \ar{rr}[swap]{!} \&\& Bz
\end{tikzar}
\\[3ex]
\begin{tikzar}[row sep = 7ex,column sep = 1.2em]
\&\& Bz\times \CCcc(z,z)\ar{dd}{!}
\\
\& Bz \ar{ur}[description]{\id \times \iota_{z}} 
\ar[equals]{dr} 
\\
\&\&Bz
\end{tikzar}
\end{split}
\eeq}
These diagrams express the  \emph{functoriality}\/ conditions for $A$ and $B$, since they say that the assignments
\[
\begin{tikzar}[row sep = 2ex]
A \colon \CCcc^{\op} \ar{r}\& \Set
\\[-1ex]
\ \ \ \ y \ar[mapsto]{r} 
\& Ay \ar[thin]{dd}[description]{f{\ast}} 
\\  \\
\ \ \ \ x \arrow[mapsto]{r} \arrow[thin,shift right=1ex]{uu}[description]{f} \& Ax  
\end{tikzar}
\qquad\qquad\qquad\qquad
\begin{tikzar}[row sep = 2ex]
B \colon \CCcc \arrow{r}\&  \Set
\\[-1ex]
\ \ \ \ x\ar[thin,shift left=1ex]{dd}[description]{f} \ar[mapsto]{r} 
\&  Bx\arrow[thin]{dd}[description]{!f} 
\\ \\
\ \ \ \ y  \ar[mapsto]{r} 
\& By \end{tikzar}
\]
are functors, in the sense that 
\begin{align*}
f\ast(g\ast a) & = (f\bullet g)\ast a\qquad \qquad & (b!f)! g = b!(f\bullet g)\\
\iota \ast a & = a & b! \iota = b
\end{align*}
where the functor notations $Ag(a)$ and $Bf(b)$ are replaced with $g\ast a$ and $b!f$.

\subsubsection{\Actions\ as functions and functors}
Cayley representations of a monoid $\MMm$ are obtained by transposing\footnote{Computer scientists would sey ``by currying''.} its actions
\[
\prooftree
M\times X\to X
\justifies
M\to \Set(X, X)
\endprooftree
\qquad
\qquad
\prooftree
Y\times M \to Y
\justifies
M\to \Set(Y,Y)
\endprooftree
\]
In this way, the \actions\ of the monoid $M$ induce its representations in terms of monoids of functions on sets. The action of the monoid on its own elements embeds it into the monoid of functions on its underlying set:
\[
\prooftree
M\times M\tto\mmult M
\justifies
M\mono \Set(M,M)
\endprooftree
\]
This is the representation that Cayley spelled out for groups. Since groups are monoids where every element has an inverse, the same representation goes through for monoids. Each element of a monoid is thus represented by the endofunction that it induces by acting on the set of all elements of the same monoid, be it on the left, or on the right. For a group, the presence of the inverses constrains the representing functions to bijections. 

Stepping up to categories again, left and the right \actions\ represent morphisms as functions
\beq\label{eq:catactions}
\prooftree
\CCcc(x,y)\times  Ay\to Ax
\justifies
\CCcc(x,y) \to \Set(Ay, Ax)
\endprooftree
\qquad
\qquad
\prooftree
Bx\times \CCcc(x,y)\to By
\justifies
\CCcc(x,y)\to \Set(Bx,By)
\endprooftree
\eeq
In particular, a category $\CCcc$ acts by the compositions on its own hom-sets
\beq\label{eq:cat-act-homsets}
\prooftree
\CCcc(x,y)\times \CCcc(y,z)\tto{\ \bullet\ } \CCcc(x,z)
\justifies
\CCcc(x,y) \to \prod_{z\in \CCcc}\Set\big(\CCcc(y,z), \CCcc(x,z)\big)\qquad \quad\Big|\quad\qquad \CCcc(y,z) \to \prod_{x\in \CCcc}\Set\big(\CCcc(x,y), \CCcc(x,z)\big)
\endprooftree
\eeq
But conditions \eqref{eq:cataction} say that the representations of objects as sets  
\begin{align}
\mndu z \colon \CCcc^{\op} &\to \Set & \cmnu x\colon \CCcc & \to \Set\\
x& \mapsto \CCcc(x,z) & z & \mapsto \CCcc(x,z)\notag
\end{align}
are functorial. The functors representing the objects in this way are said to be \emph{representable}\/ by the objects.

\subsubsection{Cayley representation and Yoneda embedding}
Collecting the $x,y$-indexed families of representations brings about the \emph{Yoneda embeddings}: 
\beq\label{eq:yoneda}
\prooftree
\CCcc(x,y) \to \Set\big(\mndu y, \mndu x\big)
\justifies
\mndu\colon \CCcc \to \Set^{\CCcc^{\op}}
\endprooftree
\qquad \qquad 
\prooftree
\CCcc(y,z) \to \Set\big(\cmnu y, \cmnu z\big)
\justifies
\cmnu \colon \CCcc \to \left(\Set^{\CCcc}\right)^{\op}
\endprooftree
\eeq
The embedding $\mndu$ represents the objects $x\in \CCcc$ as contravariant functors $\mndu x\in \Set^{{\CCcc^{op}}}$, while $\cmnu$ represents them as covariant functors $\cmnu x\in \Set^{\CCcc}$. It is often convenient to view the former as left \actions, and the latter as right \actions: 
\beq\label{eq:representable-actions}
\CCcc(x,y)\times  \mndu z (y)\to \mndu z(x)
\qquad 
\qquad
\qquad
\cmnu x (y) \times \CCcc(y,z)\to \cmnu x(z)
\eeq

\subsection{Monads, algebras}\label{Appendix:monad}
Monads are useful not only for tightening the completions in the present paper, but also for encapsulating the datatype constructors and effects in computation, and as a succinct and uniform view of algebraic theories. A monad is a functorial presentation of the free algebra construction for a given algebraic theory. E.g., the theory of groups can be viewed in terms of a functor $\Gamma:\Set\to \Set$ which maps an arbitrary set $X$ to the free group $\Gamma X$ generated by $X$. Any function $f:X\to Y$ then induces a homomorphism $\Gamma f:\Gamma X\to \Gamma Y$ between the free groups. The structure of any given  group $G$ can then be captured in a canonical way by a function $\gamma:\Gamma G\to G$, where $\Gamma G$ is the free group generated over the underlying set of the group $G$, forgetting its structure. This structure is recovered by evaluating every expression $q\in \Gamma G$, formed from the elements of $G$ by the group operations, into a group element $\gamma(q)\in G$ using the group operations of $G$. 
Note that the group structure of the free group $\Gamma X$ for any set $X$ is captured by an evaluation $\mu_X:\Gamma\Gamma X\to \Gamma X$. The family of such evaluations turns out to form a natural transformation $\mu:\Gamma\Gamma\to \Gamma$, which together with a natural transformation $\eta:\Id\to \Gamma$, whose components $\eta_X:X\to \Gamma X$ embed the generator set $X$ into the free group $\Gamma X$, forms the \emph{group monad}\/ $(\Gamma,\mu,\eta)$. Along the same lines, all other algebraic theories, for example, monoids and semigroups on the one hand, and abelian groups, rings, or lattices on the other hand, can be captured as such triples\footnote{For about 30 years, monads were called \emph{"triples"}\/ by a significant fraction, perhaps a majority, of researchers \cite{BarrM:ttt}. Before that, they were recognized in sheaf theory, albeit in dual form, as a \emph{"standard construction"}\/ \cite{Godement}, and even that name stuck  for a while \cite{HuberP:monad}, so that the \emph{"triples"}\/ proposal \cite{Eilenberg-Moore} came as a relief of sorts.}, providing a uniform framework for studying their algebras.

A monad is an endofunctor $T:\SSS\to \SSS$ with a pair of natural transformations    
\bea\label{eq:monadspan}
\TTt & = & \left(\begin{tikzar}[row sep = 8ex,column sep = 7ex] \Id \ar{r}{\eta} \& T \& TT \ar{l}[swap]{\mu} \end{tikzar} \right)
\eea
making the following diagram commute 
\beq\label{eq:monad}
\begin{tikzar}[row sep = 6ex,column sep = 1.8em]
TTT\ar{rr}{T\mu} \ar{dd}[swap]{\mu T}
\&\& TT\ar{dd}{\mu}
\\
\& T\ \ \ar[equals]{dr}\ar{ur}[description]{T\eta} \ar{dl}[description]{\eta T} \\
TT \ar{rr}[swap]{\mu} \&\& T
\end{tikzar}
\eeq
Comparing (\ref{eq:monadspan}--\ref{eq:monad}) with \eqref{eq:monspan} and \eqref{eq:monoid} makes it clear that a monad is a monoid in a functor category. A reader who made their way through any of the category theory textbooks will easily find a way to formalize  this observation. The conceptual consequence of interest here is that the algebras for the algebraic signature that induces the monad precisely correspond to the morphisms $\alpha:TA\to A$ in $\SSS$ that make the following diagram commute
\beq\label{eq:algebra}
\begin{tikzar}[row sep = 6ex,column sep = 1.8em]
TTA\ar{rr}{T\alpha} \ar{dd}[swap]{\mu_A}
\&\& TA\ar{dd}{\alpha}
\\
\& A\ \ \ar[equals]{dr}\ar{ur}[description]{\eta_A} \ar{dl}[description]{\eta_A} \\
TA \ar{rr}[swap]{\alpha} \&\& A
\end{tikzar}
\eeq
One of the $\eta_{A}$s in \eqref{eq:algebra} is redundant, but we keep both to emphasize the symmetry echoing \eqref{eq:monad}. Formally, \eqref{eq:monad} is a diagram of natural transformations, but the diagram of their components, obtained by appending on the right a set of $X$ to each sequence of $T$s, says that  \eqref{eq:monad} is a special case of the algebra diagram in \eqref{eq:algebra}, and that $\mu:TTX\to TX$ is the structure map of the free $T$-algebra $TX$. The \emph{category of $T$-algebras}\/  $\SSS^{T}$, as defined by Eilenberg and Moore in \cite{Eilenberg-Moore}, has the $\SSS$-morphisms $\alpha:TA\to A$ satisfying \eqref{eq:algebra} as its objects, and the $\SSS$-morphisms $f:A\to B$ such that $\alpha\bullet f = Tf\bullet \beta$ are the $\SSS^{T}$-morphisms, where $\beta:TB\to B$ is another $T$-algebra. The subcategory $\SSS_{T}\inclusion \SSS^{T}$ spanned by the free algebras $\mu:TTX\to TX$ can be equivalently but more succinctly presented in the form spelled out by Heinrich Kleisli \cite{KleisliH}  with the object class and the hom-sets 
\beq\label{eq:kleisli}
|\SSS_{T}| = |\SSS| \qquad \qquad\qquad \qquad \SSS_{T}(X,Y) = \SSS(X,TY)
\eeq
composed by
\bea\label{eq:kleislicomp}
\SSS_{T}(X,Y) \times \SSS_{T}(Y,Z) & \tto{\ \ \boxdot\ \ } & \SSS_{T}(X,Z)\\
\left<\Big(X\tto \varphi TY\Big), \Big(Y\tto\psi TZ\Big)  \right>& \longmapsto & \left(X\tto \varphi TY \tto{T\psi} TTZ\tto\mu TZ\right)\notag
\eea
with $\eta_{X}\in \SSS_{T}(X,X)$ as the identity morphism. Categories of algebras and of free algebras are studied in many textbooks and monographs \cite{Adamek-Lawvere,BarrM:ttt,MacLaneS:CWM,ManesE:algt}.

\subsubsection{Monoids induce monads}\label{Appendix:monmon}
In appendices \ref{Appendix:category} and  \ref{Appendix:action}, we considered preorders, monoids, and categories as structures within the category $\Set$ of sets and functions. The same structures can be defined and studied in any category $\SSS$ with cartesian products, represented by a functor $(\times )\colon\SSS\times \SSS\to \SSS$ with a natural bijection 
\bear
\SSS(X,A)\times \SSS(X,B) & \begin{tikzar}[row sep = 4em]\hspace{.1ex} \ar[bend left]{r}{<-,->} \ar[phantom]{r}[description]{\cong} 
\& \hspace{.1ex}
\ar[bend left]{l}{\left<\pi_{A}, \pi_{b}\right>}
\end{tikzar}
&\SSS(X, A\times B)
\eear
while $\times$ on the right and the pairing $<-,->$ pointed at it are in $\SSS$, whereas $\times$ on the left and the corresponding pairing are in $\Set$. Moreover, any monad $\MMm = \left(1\tto\iota M\oot{\mmult}M\times M\right)$ in $\SSS$ induces two monads
\begin{align}
(M\times) \colon \SSS& \to \SSS & ( \times M) \colon \SSS& \to \SSS\\
X & \mapsto M\times X & Y &\mapsto Y\times M 
\notag
\end{align}
where the monoid structures induce the monad structures
\[
X \ \tto{\eta = (\iota\times \id)}\  M\times X\ \oot{\mu = (\mmult\times \id)} M\times M\times X\qquad\qquad
Y \ \tto{\eta = (\id\times \iota)}\  Y\times M\ \oot{\mu = (\id\times \mmult)} Y\times M\times M
\]
Instantiated to these monads, the conditions for monad algebras in \eqref{eq:algebra} are precisely the requirements for monoid-actions in \eqref{Appendix:action}. Hence
\[
\Do\MMm = \SSS^{(M\times)}\qquad\qquad\mbox{ and }\qquad\qquad \Up\MMm = \left(\SSS^{(\times M)}\right)^\op
\]
This is just a change of angle, since corresponding categories are identical, i.e., comprising the same objects, the same morphisms, with the same structures, just satisfying different definitions, one in terms of monoids, the other in terms of monads. This change of angle, however, makes a difference in Sec.~\ref{Sec:real}
%


\section{Appendix: Limits}\label{Appendix:lim}

While preorders are a special case of categories, posetal infima and suprema are a special case of the categorical limits and colimits. Indeed, infima and suprema can be defined by the adjunctions
\beq\label{eq:inf-sup}
\begin{tikzar}{}
\PPp \arrow[bend right=10]{rr}{\scriptstyle \bot}[swap]{\mnd} \&\& \Do \PPp \arrow[bend right=10]{ll}[swap]{\vee}
 \&\&\& \PPp \arrow[bend right=10]{rr}{\scriptstyle \top}[swap]{\cmn} \&\& \Up \PPp  \arrow[bend right=10]{ll}[swap]{\wedge}
\end{tikzar}
\eeq
which means 
\[
\vee L \leq x \iff L\subseteq \mnd x \qquad \quad\qquad  x\leq \wedge U \iff \cmn x \supseteq U
\]
limits and colimits can be defined using the adjunctions
\beq\label{eq:lim-colim}
\begin{tikzar}{}
\CCc \arrow[bend right=10]{rr}{\scriptstyle \bot}[swap]{\mnd} \&\& \Do \CCc \arrow[bend right=10]{ll}[swap]{\supp}
 \&\&\& \CCc \arrow[bend right=10]{rr}{\scriptstyle \top}[swap]{\cmn} \&\& \Up \CCc  \arrow[bend right=10]{ll}[swap]{\inff}
\end{tikzar}
\eeq
which correspond to the natural families of bijections
\[
\CCc(\supp \lft A, x) \cong \Do\CCc(\lft A, \mnd x) \qquad\qquad \quad  \CCc(x,\inff \rgt B) \cong \Up\CCc(\cmn x, \rgt B)
\]
The lattices $\Do\PPp$ and $\Up\PPp$ are defined in \eqref{eq:douppos}, and the categories $\Do\CCc$  and $\Up\CCc$ in \eqref{eq:doupcat}. 

\section{Appendix: Arrow and comma categories}\label{Appendix:comma}

The \emph{arrow category}\/ $\arrow\CCc$ \cite[II.4]{MacLaneS:CWM} and the \emph{twisted arrow category} $\tarrow\CCc$ \cite[IX.6.Ex.3]{MacLaneS:CWM} of a given category $\CCc$ both have the $\CCc$-arrows as the objects, and the pairs of $\CCc$-arrows as the morphisms:
\bear
\arrow\CCc (x, y) &=& \left\{\ <f,f'>\ \ \  \Bigg|\ \ \ \ \ 
\begin{tikzar}[row sep=1.8em,column sep=1.8em]
X \ar{r}{f} \ar{d}[description]{x} \& Y \ar{d}[description]{y}  \\ 
X' \ar{r}[swap]{f'}\& Y'
\end{tikzar}\  \right\}\\[2ex]
\tarrow\CCc (x, y) &= & \left\{\ <t,t'>\ \ \  \Bigg|\ \ \ \ \ 
\begin{tikzar}[row sep=1.8em,column sep=1.8em]
X \ar{r}{t} \ar{d}[description]{x} \& Y \ar{d}[description]{y}  \\ 
X' \& Y' \ar{l}{t'}
\end{tikzar}\  \right\}
\eear
Note the opposite orientations of the second components $f'$ and $t'$. This difference accounts for two different ways to project the two categories to their base category:
\[
\begin{tikzar}[row sep = 6ex]
\CCc/\CCc \ar[shift right=2ex]{d}[swap]{\scriptstyle \Dom} 
\ar[shift left=2ex]{d}{\scriptstyle \Cod} 
\\
\CCc \ar{u}[description]{\scriptstyle \Ids}
\end{tikzar}
\qquad \qquad\qquad \qquad
\begin{tikzar}[row sep = 6ex]
\tarrow{\CCc} \ar{d}[description]{\scriptstyle <\Dom,\Cod>} \\
\CCc \times \CCc^\op
\end{tikzar}
\]
While the projection on the left presents the category $\CCc$ as a truncated simplicial set, the projection on the right is the discrete fibration corresponding to the presheaf $hom: \CCc^\op\times \CCc\to \Set$, which maps any pair of objects $x,y$ of $\CCc$ to the hom-set $\CCc(x,y)$.

For any given pair of functors $F:\AAa\to \CCc$ and $G:\BBb\to \CCc$, the \emph{comma}\/ construction \cite[II.6]{MacLaneS:CWM} is defined  by
\bea
\lvert F/G\rvert & = &  \coprod_{\substack{a\in\AAa\\ b\in \BBb}} \CCc(Fa, Gb) \\
F/G(\varphi^{Fa}_{Gb}, \psi^{Fc}_{Gd}) & = &\left\{<f,g>\in \AAa(a,c)\times \BBb(b,d)\ \Big|\ \begin{tikzar}[column sep = 1.5em]
Fa\ar{r}{Ff} \ar{d}[description]{\varphi}\& Fc\ar{d}[description]{\psi}\\
Gb\ar{r}[swap]{Gg} \& Gd
\end{tikzar}
\right\}\notag
\eea
They come with the functors $F/G \tto{<\Dom,\Cod>} \AAa\times \BBb$ which map the arrows to the sources of their domains and codomains. 

When $\AAa=\CCc$ and $F$ is the identity functor, the comma notation $\Id_\CCc/G$ is often simplified to $\CCc/G$. Similarly when $\BBb=\CCc$, one often writes $F/\CCc$ instead of $F/\Id_\CCc$.

When $\BBb = 1$ and $G$ is the constant functor to $x\in \CCc$, the comma construction gives the \emph{slice}\/ category $\CCc/x$ of arrows into $x$. Similarly when $\AAa = 1$ and $F$ is the constant functor to $x\in \CCc$, the comma construction gives the \emph{coslice}\/ category $x/\CCc$ of arrows from $x$.

\section{Appendix: Comprehending diagrams}\label{Appendix:diags}

Limits and colimits were originally defined for diagrams as graph morphisms into categories, and settled in  \cite[Ch.~III, Sections 3--4]{MacLaneS:CWM} for diagrams as arbitrary functors from small categories. But a supremum and an infimum of an arbitrary poset are equal, respectively, to the supremum of its lower closure and to the infimum of its upper closure. The definitions streamlined in terms of the lower and the upper sets, given in \eqref{eq:inf-sup}, are slightly more convenient for posets and preorders but much more convenient when lifted to \eqref{eq:lim-colim} for categories. Here we spell out the reduction of colimits and limits over arbitrary diagrams to colimits and limits over discrete fibrations and opfibrations, i.e., left and right \actions. For preorders, the reductions of joins and meets of an arbitrary subset $S$ to its lower and upper closures was
\[\prooftree
S\subseteq \PPp
\justifies
\rule{0ex}{2.8ex}
\lft S = \bigcup_{x\in S} \mnd x\qquad \qquad\qquad\qquad \rgt S = \bigcup_{x\in S} \cmn x
\endprooftree\]
giving 
\[
\bigvee S = \bigvee \lft S \qquad \qquad \qquad \qquad  \bigwedge S = \bigwedge \rgt S
\]
The analogous reductions for the colimit and the limit of an arbitrary diagram or functor $\DDd\tto D \CCc$ are based on the comma categories
\[\prooftree\prooftree
\DDd \tto D \CCc
\justifies
\comm\CCc D \tto{\Dom} \CCc\qquad \qquad \qquad \qquad \comm D\CCc \tto{\Cod} \CCc 
\endprooftree
\justifies
\lft \DDd\tto{\lft D}\CCc \qquad \qquad \quad\qquad \qquad \rgt\DDd \tto{\rgt D} \CCc 
\endprooftree
\]
where the objects of $\lft\DDd$ and $\rgt\DDd$ are connected components\footnote{Intuitively, a connected component of a space is a subspace where every two points can be connected by a path. Analogously, a connected component of a category is a subcategory where every two objects can be connected by a path through its arrows, ignoring the arrowheads.} of $\comm D\CCc$ and $\comm\CCc D$ 
\begin{align*}
\lvert \lft \DDd\rvert & = \lvert \comm \CCc D\rvert \Big/ \underset\DDd{\approx} &  \lvert \rgt \DDd\rvert & = \lvert \comm D\CCc\rvert \Big/ \overset\DDd\approx
\end{align*}
and $\underset\DDd{\approx}$ and $\overset\DDd\approx$ are the symmetric transitive closures of the relations $\underset\DDd{\sim}$ and $\overset\DDd\sim$ defined by
\bear
\alpha_{Dx}^a \underset\DDd{\sim} \beta^b_{Dy} & \iff & \begin{tikzar}[column sep = 1.5em]
\&a=b\ar{dl}[description]{\alpha}\ar{dr}[description]{\beta}\\
Dx \ar{rr}{\exists u\in \DDd(x,y)}[swap]{Du} \&\& Dy
\end{tikzar}\\[3ex]
\alpha^{Dx}_a \overset\DDd{\sim} \beta_b^{Dy} & \iff & \begin{tikzar}[column sep = 1.5em]
Dx \ar{rr}{\exists u\in \DDd(x,y)}[swap]{Du} \ar{dr}[description]{\alpha} \&\& Dy \ar{dl}[description]{\beta}\\
\&a=b
\end{tikzar}
\eear
Since an object $\left[\beta^b_{Dy}\right]\in \lvert\lft \DDd\rvert$ is thus a family of arrows with the same domain, the domain functor $\comm \CCc D\tto\Dom \CCc$ induces the projection $|\lft D|\colon |\lft \DDd| \to |\CCc|$ with $\lft D\left[\beta^b_{Dy}\right] = b$. Moreover, any $\CCc$-morphism $f:d\to b$ lifts the relation $\alpha^a_{Dx}\sim\beta^b_{Dy}$ to $f\bullet \alpha^a_{Dx}\sim f\bullet\beta^b_{Dy}$. The $\CCc$-arrows thus induce the inverse images which map connected components to connected components in $\comm \CCc D$. The $\DDd$-morphisms can therefore be defined by setting $\comm{\lft \DDd} {\left[\beta^v\right]} = \comm \CCc v$, which makes $\lft D$ into a discrete fibration. 
The discrete opfibration $\rgt D:\rgt\DDd\to \CCc$ is defined dually. See \cite{Street-Walters:comprehensive} for a full account of such \emph{comprehensive factorizations}\/ of functors.

\begin{proposition} \label{Prop:ess-small}
The following statements are true
\begin{enumerate}[a)]
\item $\supp\left( \DDd \tto D \CCc\right)\ =\ \supp \left(\CCc/D \tto{\Dom} \CCc\right) \ =\ \supp \left(\lft\DDd \tto{\lft{D}} \CCc\right)$
\item $\inff\left( \DDd \tto D \CCc\right)\ =\ \inff \left(D/\CCc \tto{\Cod} \CCc\right) \ =\ \inff \left(\rgt{\DDd} \tto{\rgt D} \CCc\right)$
\end{enumerate}
\end{proposition}
 
 \bpr 
Consider case (b). A family of arrows $\gamma = \{\gamma_x\in \CCc(c,Dx)\ |\ x\in \DDd\}$ is a cone to the diagram $\DDd\tto D\CCc$ if and only if all $u\in \DDd(x,y)$ satisfy $Du\circ \gamma_x = \gamma_y$. The following diagram 
\beq\label{eq:redcones}\begin{tikzar}[column sep = 1.5em]
\& c \ar{dl}[description]{\gamma_x}\ar{dr}[description]{\gamma_y}\\
Dx \ar{rr}[swap]{Du} \ar{dr}[description]{\alpha} \&\& Dy \ar{dl}[description]{\beta}\\
\&a=b
\end{tikzar}\eeq
shows that the family $\widehat \gamma = \{\widehat \gamma_\alpha \in \CCc\big(c,\Cod( \alpha) \big)\ |\ \alpha \in D/\CCc\}$ defined by 
\bear
\widehat \gamma_\alpha & = & \alpha\circ\gamma_x
\eear
is a cone to the diagram $D/\CCc\tto{\Cod} \CCc$. 
Conversely, given a cone $\widehat{\gamma} = \{\widehat{\gamma}_\alpha \in \CCc\big(c,\Cod( \alpha) \big)\ |\ \alpha \in D/\CCc\}$ to $D/\CCc\tto{\Cod} \CCc$, the induced cone $\gamma$ to $\DDd\tto D\CCc$ is defined by
\bear \gamma_x &= & \widehat{\gamma}_{\id_{Dx}}
\eear 
Diagram \eqref{eq:redcones} also shows that cones $\widehat{\gamma}$ to $D/\CCc\tto{\Cod} \CCc$ must satisfy 
\bear
\alpha  \overset\DDd{\sim} \beta & \implies & \widehat{\gamma}_\alpha = \widehat{\gamma}_\beta
\eear
so that a cone to $\DDd\tto D\CCc$ is always also a cone to $\lft \DDd\tto{\lft D} \CCc$, and vice versa. This proves claim (b). The argument for (a) is symmetric.
 \epr
 
\paragraph{Essentially small diagrams.} 
 A \presheaf\ $\lft{\mathsf{D}}:\lft \DDD\to \CCc$ (resp.\ \postsheaf\  $\rgt{\mathsf{D}}: \rgt \DDD\to \CCc$) is  \emph{essentially small} when there is a small diagram $D:\DDd\to \CCc$ such that the \presheaf\ $\lft D:\lft \DDd\to \CCc$ (resp.\ \postsheaf\ $\rgt D:\rgt \DDd\to \CCc$) constructed as in Prop.~\ref{Prop:ess-small} is equivalent to $\lft{\mathsf{D}}:\lft \DDD\to \CCc$ (resp.\ $\rgt {\mathsf{D}}: \rgt \DDD\to \CCc$). When confusion is unlikely, all \presheaf s written in the form $\lft D:\lft \DDd\to \CCc$ and \postsheaf s in the form $\rgt D:\rgt \DDd\to \CCc$ are tacitly assumed to be essentially small.

\section{Appendix: Splittings}\label{Appendix:split}

A categorical structure is called \emph{absolute}\/ if it is preserved by all functors. Since functors are defined to only preserve identity maps and compositions, absolute constructions must be characterized in terms of those. Here we are only concerned with absolute limits and colimits. They were studied by Bob Par\'e \cite{PareR:absolute-coeq,PareR:absolute}. There are only a few of them, and we call them \emph{splittings}. 

\subsection{Idempotents and retractions}
A \emph{retraction}\/ is a pair of maps $A\overset{q}{\underset{i}{\rightleftarrows}} B$ such that $q\circ i = \id_B$.  The type $B$ is called a \emph{retract} of $A$ when there is such a pair. The composite $\varphi = i\circ q$ is then idempotent, in the sense that $\varphi\circ\varphi = \varphi$. The retraction $A\overset{q}{\underset{i}{\rightleftarrows}} B$ is called the \emph{splitting}\/ of the idempotent $\varphi$. The following diagram summarizes a retraction 
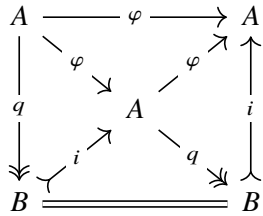
\[\begin{tikzar}{}
A \ar[two heads]{dd}[description]{q} \ar{rr}[description]{\varphi} \ar{dr}[description]{\varphi} \&\& A\\
\& A\ar{ur}[description]{\varphi}\ar[two heads]{dr}[description]{q}\\
B \ar[tail]{ur}[description]{i} \ar[equal]{rr} \&\& B \ar[tail]{uu}[description]{i}
\end{tikzar}
\]

\subsection{Split pairs and their equalizers}
The splitting of an idempotent $\varphi$ is easily seen to be an equalizer, or a coequalizer, of $\varphi$ and the identity. A \emph{split pair} is usually drawn as a mnemonic diagram
\beq\label{eq:splitpair}
\begin{tikzar}
A \arrow[shift left=1ex]{rr}[description]{f}
\arrow[tail,shift right=1ex]{rr}[description]{i}
\&\&
B\arrow[shift left = 1.5ex,bend left = 25,two heads]{ll}[description]{r} 
\end{tikzar}
\eeq
which is meant to encode two equations:
\[ r\circ i = \id_a \qquad\qquad\qquad \qquad i\circ r\circ f = f\circ r \circ f\]
It follows that $r\circ f \circ r\circ f = r\circ i\circ r\circ f = r\circ f$ is an idempotent. A split pair $f,i$ thus has $r$ as the joint splitting of the identity $r\circ i$ and of the idempotent $r\circ f$. An equalizer of a split pair is also split, and it is also drawn mnemonically
\beq\label{eq:splitpair-equalizer}
\begin{tikzar}
E\arrow[tail]{rr}[description]{e}\&\&A \arrow[shift left=1ex]{rr}[description]{f}
\arrow[tail,shift right=1ex]{rr}[description]{i}
\arrow[shift right = .75ex,bend right = 25,two heads]{ll}[swap]{q}
\&\&
B\arrow[shift left = 1.5ex,bend left = 25,two heads]{ll}{r} 
\end{tikzar}
\eeq
which is meant to encode the equations of the following less mnemonic diagram
\beq\label{eq:splitpair-equalizer-big}
\begin{tikzar}
E\arrow[tail]{rr}[description]{e}\ar[equals]{dd}\&\&A 
\\ \\
E\arrow[tail]{rr}[description]{e} \&\&
A  \ar[two heads]{uull}[description]{q} \ar[equals]{dd}
\arrow[shift left=1ex]{rr}[description]{f}
\arrow[tail,shift right=1ex]{rr}[description]{i} 
\&\&
B\arrow[two heads,shift left=.75ex]{ddll}[description]{r} \arrow[two heads,shift right=.75ex]{uull}[description]{r} 
\\ \\
\&\& A
\end{tikzar}
\eeq
Showing that $e$ is an equalizer of a pair $f$ and $i$ with a splitting $r$ if and only if $e$ has a splitting $q$ with $e\circ q = r\circ f$ is an instructive exercise in function algebra. It is an absolute equalizer because it is characterized by these equations. The fact that \emph{all}\/ absolute limits must be in this form, and the absolute colimits in the dual form of split coequalizers, has broad repercussions through mathematics.  The role of idempotents is also hard to overestimate. Note that a splitting of an idempotent $\varphi$ can be viewed as an equalizer of the split pair where $f=\varphi$ and $i = r = \id$.

\end{document}